\documentclass{article}  

\usepackage[english]{babel}
\usepackage[utf8]{inputenc}

\usepackage[left=1in, right=1in, top=1in, bottom=1in]{geometry}
\usepackage{setspace}  
\usepackage{graphicx}
\usepackage{subcaption}
\usepackage{amsmath, amssymb, amsfonts}
\usepackage{booktabs}
\usepackage{appendix}
\usepackage{verbatim}
\usepackage{multirow, multicol}
\usepackage{array}

\usepackage{hyperref}

\usepackage[authoryear]{natbib}  
\usepackage{xcolor}

\definecolor{BlueGreen}{RGB}{0,50,125}

\definecolor{mysilver}{RGB}{220,220,220}

\definecolor{mygreen}{RGB}{0,128,0}

\newcommand{\distance}{\textbf{Distance}}
\newcommand{\mnl}{\textbf{MNL}}
\newcommand{\mxlmean}{\textbf{MXL-Mean}}
\newcommand{\mxltwentyfive}{\textbf{MXL-25}}

\title{What makes a good public EV charging station? A revealed preference study}
\author{
    Steven Lamontagne$^{1}$ \quad
    Margarida Carvalho$^{1}$ 
    \quad
    Emma Frejinger$^{1}$ \quad
    Ribal Atallah$^{2}$\\[1ex]
    \small $^{1}$CIRRELT and Département d'informatique et de recherche opérationnelle, Université de Montréal\\
    \small $^{2}$Institut de Recherche d'Hydro-Québec
}
\date{}

\begin{document}
\maketitle

\begin{abstract}
To determine the optimal locations for electric vehicle charging stations, optimisation models need to predict which charging stations users will select. We estimate discrete choice models to predict the usage of charging stations using only readily available information for charging network operators. Our parameter values are estimated from a unique, revealed preferences dataset of charging sessions in Montreal, Quebec. We find that user distance to stations, proximity to home areas, and the number of outlets at each station are significant factors for predicting station usage. Additionally, amenities near charging stations have a neutral effect overall, with some users demonstrating strong preference or aversion for these locations. High variability among the preferences of users highlight the importance of models which incorporate panel effects. Moreover, integrating mixed logit models within the optimization of charging station network design yields high-quality solutions, even when evaluated under other model specifications.
\end{abstract}

\noindent\textbf{Keywords:} Electric vehicles, Discrete choice models, Mixed logit, Revealed preferences

\section{Introduction}


Whether operated by a private company or a governmental entity, electric vehicle (EV) public charging station network operators want their charging stations to have a high usage rate.  Each charging station is expensive to install, and as such, private companies require many customers to recuperate this cost, while high usage of charging stations demonstrates a wise deployment of resources for governmental entities.  To aid in this, a tool for predicting aggregate usage of public charging stations can inform on the potential impacts of candidate locations for new stations. Moreover, such a tool can then be embedded within optimisation models for charging network planning, with the goal of selecting future public charging station locations in a long-term context. 

A crucial component of this optimisation is the prediction of demand for each station. In charging station placement models, the demand for a charging station is based on either the location of the users and charging stations on the network (so called \emph{node-based models}), or the path of a user between origin and destination on a trip (\emph{flow-based models}, or \emph{activity-based models} if considering a sequence of trips)~\citep{Kchaou2021, Metais2022}.  Since it only requires the locations of users and charging stations, the node-based approach is particularly well-suited to intracity charging network optimisation, which may have limited data on the daily activities of their users. However, the level of demand in these models is typically determined via distance-based coverage or p-median models, which do not allow for more complex behaviours such as cannibalisation of demand between stations or limited charging capacity~\citep{Metais2022}.  


By contrast, the state-of-the-art in charging station demand modeling is dominated by discrete choice models, where the alternatives are the charging stations available to the driver~\citep{Potoglou2023}. The flexibility and interpretability of discrete choice models make it an ideal choice for predicting existing behaviour and extending these predictions for candidate charging stations. In addition, the underlying behavioral assumptions in these models can allow for the more complex behaviours described above. However, the integration of existing discrete choice demand models to charging network operation is non-trivial. 
Additionally, while discrete choice models can be estimated using stated preference (SP) or revealed preference (RP) data, the vast majority are from SP data. More specifically, the surveys for these SP-driven works present hypothetical trip and charging stations to respondent in order to simulate real-time decision making. As a consequence, these experiments typically include real-time operation attributes  such as remaining vehicle range, parking time at destination, or sociodemographic characteristics~\citep[e.g.][]{Wang2021, Visaria2022, Anderson2023}. Since these attributes cannot be known in advance by the charging network operator, their inclusion in the demand models renders them ill-suited for strategic optimisation.

In this work, we estimate discrete choice models specifically tailored for integration within intracity, node-based charging network optimisation. More concretely, the models predict charging station choices based on vehicle-agnostic and trip-agnostic characteristics, for example, the home area of the user 
and station characteristics. Such information is readily available by charging network operators, thus allowing for \emph{a priori} estimation of demand of both existing and candidate stations.  Recent advances in competitive facility location problems then describe the integration of these discrete choice models directly within optimisation frameworks~\citep[e.g.,][]{Mai2020, Paneque2020, Lamontagne2023a, Legault2023,shen2024sequential}. These frameworks embed the evaluation of charging station demand directly into the optimisation model, preserving the underlying behavioural assumptions of the demand models, and thus enabling the selection of optimal charging station locations. 
Public charging infrastructure is available with different power outputs, referred to as level 2 and level 3 charging, and which exhibit drastically different usage and behaviour~\citep{Hardman2018, Figenbaum2019, Tal2020}. As such, we estimate separate parameter values for each of these charging levels. 


\paragraph{Contributions}
Our work presents several contributions to EV demand modeling. 
\emph{Firstly}, we specify models tailored for charging station operators, 
ensuring their high interpretability and rendering seamless their integration into network optimisation models. 
More specifically, we estimate multinomial logit (MNL) and mixed logit (MXL), with the latter including the panel effect of repeated observations from individual users. 
\emph{Secondly}, we estimate our models from RP data rather than SP data, thus avoiding the well-studied hypothetical bias~\citep{Haghani2021a, Haghani2021b}. To the best of our knowledge, only \cite{Sun2016} uses RP data for charging station demand estimation, and only in the context of intercity EV charging. In addition, the dataset we use for estimation is unique in that the charging network represents 90\% of charging stations within the province~\citep{ChargingNetworks}, thus accounting for nearly all available alternatives. By point of comparison, the charging network operator in \cite{Pevec2018} is reported to only account for 15\% of the total network. 
\emph{Thirdly}, we complement this unique usage data with Geographical Information System (GIS) data from OpenStreetMap~\citep{OpenStreetMap}, allowing us to examine the travel network and amenities in proximity to the stations when each session took place. Notably, the presence of amenities has been observed to affect EV charging decisions in \cite{Philipsen2016, Anderson2018, Sheldon2019} but, to the best of our knowledge, has not yet been included in RP studies with discrete choice models. We note that RP data and GIS information have been used in machine learning models for charging station prediction~\citep{Pevec2018, Straka2020}. However the resulting machine learning models are not well-suited for selecting multiple new charging stations simultaneously, thus making for a difficult integration into charging network optimisation. 
\emph{Fourthly}, we present an application where we use our estimated demand model (utility functions) to select the placement of level 2 charging stations across the city of Montreal as to provide the best service to all users. This application illustrates how the solution of an optimisation model can vary based on the assumptions for the utility functions. 
Our findings indicate, unsurprisingly, that the distance between the charging stations is the most significant factor for predicting charging station usage. Additionally, the number of outlets at each station and the charging station being within a short walk of the home area of the user are also found to increase usage of a station. In the MNL models, the average effects of some amenities are significant. However, the MXL models indicate an overall neutral average effect and significant variance among users. This suggests that the significant average effect observed in the MNL models may be caused by a high number of sessions from users with strong preferences or aversions, rather than an overall trend within the population. As a consequence, optimisation models which employ MNL formulations for demand modeling may incorrectly relate the usage of stations with amenities.

\paragraph{Paper Organisation}
In Section~\ref{ThirdArticle:SectionLiteratureReview}, we present a review of the relevant literature on charging preferences and charging station selection. Section~\ref{ThirdArticle:SectionData} presents general characteristics about EV charging in Quebec, along with the process for generating the value of the attributes in our models. In Section~\ref{ThirdArticle:SectionMethods}, we discuss the specifications of our MNL and MXL models. The results of the estimation process are presented in Section~\ref{ThirdArticle:SectionResults}, while they are discussed in Section~\ref{ThirdArticle:SectionDiscussion}. Section~\ref{ThirdArticle:SectionDiscussion:opt} compares charging station network design solutions for optimisation problems using different estimated utility functions. Finally, we conclude our work in Section~\ref{ThirdArticle:SectionConclusion}.

\section{Literature Review}
\label{ThirdArticle:SectionLiteratureReview}

There is a vast literature on EV charging habits. Thus, for the sake of brevity, we present here works relating to the use or selection of charging stations by private EV owners. We start with a focus toward the data used, followed by an exploration of literature addressing charging preferences. This specifically includes preferences related to the public charging network. Lastly, we delve into the methodologies employed for predicting charging station selection and we discuss the relation to our work. For recent reviews of these and other aspects of charging activity, we refer to \cite{Hardman2018} and \cite{Potoglou2023}. 
\subsection{Charging Behaviour and Requirements}

We start by noting the works detailing the charging behaviour for private vehicle users. These allow us to identify and validate charging behaviour of real users, as well as isolate irregular behaviour. Collectively, these works report on millions of charging sessions and thousands of EV owners across the Netherlands, British Columbia (Canada), Norway, Ireland, California (United States), and New Zealand~\citep{vanDenHoed2013, Axsen2015, Figenbaum2016, Morrissey2016, Figenbaum2019, Nicholas2017a, Tal2018, Tal2020, EECA2021}. Additionally, while other works focus on a specific research question, they may also present survey results or RP data on charging behaviour. Examples include the surveys of \cite{Franke2013, Chakraborty2019, Lee2020, Visaria2022}, and \cite{Anderson2023}. By combining this vast literature, we obtain a diverse portfolio of preferences of EV drivers across multiple years and geographical regions, thus rendering it possible to compare against existing behaviour and identify inconsistencies.


The dependency on the public charging network varies greatly depending on the user. When charging their vehicle, users may have access to private charging outlets located at their home or workplace. Indeed, currently, the majority of users charge primarily at home, with between 82\% and 93\% of EV owners recharging frequently at home~\citep{Figenbaum2019, Lee2020, Tal2020, Visaria2022, Anderson2023}.  For workplace charging, there is more variability in terms of availability and frequency of use, with the percentage of frequent users ranging from 19\% to 40\%~(\emph{ibid.}). Notably, in \cite{Helmus2020} they find that a minority of sessions in the Netherlands are associated with workplace charging, at around 14\%. This compares with 35\% associated with non-workplace daytime charging and 50\% associated with overnight public charging. The attributes contributing to the frequency of each type of charging location were examined in more detail in \cite{Chakraborty2019} and \cite{Lee2020}, which included income, dwelling type, access to level 2 charging at home versus level 1 charging, and the electric range of the vehicle. The access to home or workplace charging can decrease the reliance on the public charging network, and can thus lead to a disproportionate representation of public charging sessions by those without access to these other types of charging. 

\subsection{Predicting Charging Station Selection}

Once a user has opted to use the public charging network, they must then select a charging station to use. Depending on its characteristics, each station may successfully draw the demand of that user. In this sense, this viewpoint is similar to the approach of node-based optimisation models, with the demand of each charging station depending on multiple attributes rather than simply the distance. In the demand modeling literature, researchers have then tried to characterise or predict this selection. In many cases, they use a discrete choice model to predict a charging station choice, with significant attributes including cost, charger availability or waiting time, distance to home or detour time, location type, proximity to amenities, number of chargers at the station, electricity obtained from a renewable source, and the rating of the charging station on a mobile phone application~\citep{Luo2015, Cui2018, Moon2018, Sheldon2019, Wang2021, Lamontagne2023a, Ma2022, Visaria2022, Anderson2023}. In the case of \cite{Luo2015, Cui2018, Lamontagne2023a}, charging demand was predicted using discrete choice models, and embedded within optimisation models. No parameter values or estimation results are presented in \cite{Luo2015, Cui2018}, while parameter values for a MNL model are presented in \cite{Lamontagne2023a} based on real charging session data. In all other cases, the models are estimated based on a curated survey, where respondents are presented a set of charging options with different characteristics and must select an option. 

Of particular note for our study is the work in \cite{Sheldon2019}, involving a choice experiment with 1,261 drivers (not necessarily EV) in California, as it includes the proximity to amenities as an attribute. More specifically, the proximity of a station to amenities is treated as a categorical characteristic which classifies charging stations into one of the following location types: workplace, grocery stores, shopping malls, public transit, sports facilities, schools, entertainment venues, level 3 near home, and level 3 near the highway.
They find that respondents exhibited a preference for locations identified as near grocery stores, shopping malls, and having level 3 charging stations near their homes and close to the highway network, as opposed to the baseline location near entertainment locations. 
Moreover, respondents demonstrated indifference or a tendency to disfavor locations identified as being near transit, sports facilities, and schools (even students). In all cases, these public charging locations were favoured less than workplace charging. However, the authors find that these preferences vary significantly among the respondents.

Under a similar node-based viewpoint but not using discrete choice models, in \cite{Philipsen2015}, they ask a discussion group of 15 non-EV drivers from Germany about their preferences for charging. As part of this, a series of evaluation criteria for charging stations were brought up: combining the charging session with everyday activities (dual use), compatibility with existing habits, avoiding detours or added delays (accessibility), easy to find and see (visibility), availability of chargers when needed (reliability), safe to leave vehicle or to stay for extended periods of time, allowing for longer trips that could otherwise not be done, and connection to the public transportation network. The participants also proposed several locations which fit the above criteria, such as supermarkets, public authorities,  gas stations, motorway service stations, medical centres, recreational facilities, and sports venues. 

These findings were validated in \cite{Philipsen2016} through a survey of 252 respondents in Germany, asking each to rank the importance of the evaluation criteria and locations proposed in \cite{Philipsen2015}.  For the ranking of criteria, reliability and dual use were most important, followed by accessibility and visibility. For location type, motorway service stations were deemed most important, followed by workplace, gas stations, and shopping facilities (at roughly equal importance), then leisure facilities and educational facilities. These preferences for evaluation criteria and location were found to vary based on sociodemographic (notably, based on gender) and if the respondent was presently an EV driver. 

In a similar approach, in \cite{Anderson2018}, they conduct a survey of 761 EV owners in Germany in which respondents were asked where to place additional charging stations based on their personal needs. For each station, they could indicate the power level (3.7 kW, 22 kW, 50 kW), the projected frequency of use, and the activity while charging at that location (such as work, stop-to-charge, or shopping). Overall, users placed charging stations within cities, preferred 22 kW charging, and wanted activities close to their charging stations (with only 29\% of the stations indicating stop-to-charge). Additionally, the activities selected varied considerably between the different power levels, highlighting different usage.

Across the node-based viewpoint, we note the lack of RP data, relying exclusively on choice experiments or survey results. Additionally, while the distance or accessibility of a station is crucial, there is support for the idea that users prefer activities near their public charging stations (an idea which is implicit in the activity-based optimisation). As such, this suggests that the inclusion of such attributes in node-based and activity-based optimisation would be beneficial.


Rather than focusing on the charging stations in isolation, other works consider the station as one part of a user's trip between their origin and destination. In this context, the selection of a charging station can thus be viewed as the deviation of a path for the purposes of charging. As such, this viewpoint is akin to the path-based optimisation models. In all of these works, discrete choice models are used to predict which charging station is selected by the users, however the exact decision being modelled can vary. In both \cite{Sun2016} and \cite{Yang2016}, some routes can be completed without charging, so users must decide if they select a route with charging at all and, if applicable, which route with charging. They find that the state of charge and sociodemographic characteristics were significant attributes for selecting a route with charging. In terms of route deviation for charging, \cite{Sun2016} noted that maximum acceptable deviations from the shortest path varied between 500 metres and 1,750 metres depending on vehicle type (private versus commercial) and day (weekday versus weekend). Similarly to \cite{Sun2016, Yang2016}, in \cite{Ashkrof2020}, users must select a route between their origin and destination, with some routes including fast charging stations along the way and some destinations including a slow charger.  They find that availability of a slow charging station at the destination decreased the chances of selecting a route with fast charging, while female drivers and younger drivers were more likely to select routes with fast charging. In general, for shorter trips, local streets are preferred while for longer trips freeways and arterial ways are preferred.

Instead of selecting a route, in \cite{Ge2022}, users are presented with a route and vehicle characteristics. The process then simulates a real trip, where users arrive at a charging station with given vehicle and trip characteristics (e.g. remaining state of charge and distance to next charging station), and then must decide whether to charge at this charging station or to continue on the route. They find that the remaining state of charge, charging cost, access time, and having access to full amenities (restroom, Wifi, and a restaurant) were significant attributes for users deciding to charge. In terms of data for discrete choice model estimation, we note that only \cite{Sun2016} uses revealed-preference data, while all others use survey data. 

Overall, these path-based viewpoints consider analogous attributes to those of node-based, such as the distance to the charging station being replaced by path deviation. However, a major difference is in the data assumption about users, where the researcher is assumed to be aware of the path (or at least the origin and destination) of the user in the path-based viewpoint.



\subsection{Relation to Our Work}

As with many of the prior articles, our work uses discrete choice models to predict the choice of charging station by users, more specifically following the node-based viewpoint for intracity charging. However, our work expands upon the existing literature by estimating with RP data, previously only used in the intercity model of \cite{Sun2016}. Additionally, by combining this RP data with GIS information from OpenStreetMaps, we can include amenity information, previously only considered in discrete choice estimation with SP data such as in \cite{Sheldon2019, Visaria2022}. By comparing the available amenities near a charging station with the usage of the stations, we can then empirically verify the value of the location types proposed in \cite{Philipsen2015, Philipsen2016, Anderson2018}. Finally, we assess the impact of different estimated models when used to optimize a charging network.


\section{Data}
\label{ThirdArticle:SectionData}

In our work, we use a dataset of charging sessions from the charging network operator \emph{Circuit Électrique}. We begin, in Section~\ref{ThirdArticle:SectionDataDescription}, with a brief discussion about the charging network and the charging process in Quebec, while more details can be found in~\ref{ThirdArticle:AppendixDataShared}.  As we are interested in intracity charging, in Section~\ref{ThirdArticle:SectionDataProcessing}, we isolate the data for the city of Montreal, Quebec. Finally, in Section~\ref{ThirdArticle:SectionAttributeEnconding}, we describe the procedure for obtaining attribute values which are not found in the charging session data.

\subsection{Data Description}
\label{ThirdArticle:SectionDataDescription}

 A key component of our charging session data, and a notable difference with other charging networks, is the membership card. Indeed, in addition to managing operations of the charging stations, \emph{Circuit Électrique} also has a membership program to which users can register~\citep{Membership}. Registered users can receive by mail one or more cards with Radio Frequency Identification (RFID) capabilities, which are connected to their account information; All sessions initiated by a membership card (including at some charging stations not operated by \emph{Circuit Électrique}) are then linked and charged to the associated account. However, a membership card is not required at all charging outlets operated by \emph{Circuit Électrique}, as sessions may be initiated by guest accounts and charged to a credit or debit card. As such, recorded sessions contained within the dataset include those completed by registered and guest members at  \emph{Circuit Électrique} charging, as well as some sessions completed by registered members at third-party charging stations. 


Since we cannot accurately identify the preferences and habits of guest accounts (which represent roughly 47\% of users and 11\% of sessions), they are excluded from the analysis and estimation processes.
 Even for non-guest members, reliable vehicle information is not available for member accounts and, in addition, multiple vehicles may be associated with the same account. As such, it is not directly known if the account corresponds to a private user, or a fleet of commercial vehicles. Since our aim is to evaluate public charging stations from the perspective of private vehicle owners, we must identify and isolate the corresponding member accounts. To do this, we have combined the expertise of our industrial partner, as well as the descriptions of typical charging behaviour from the literature, to create a series of filters for categorisation. For the purposes of these filters, we do not take into account charging behaviour or sessions which occurred between February 1st 2020 and June 30th 2021, as the COVID-19 pandemic had noticeable effects on charging habits which could interfere with the filters.
 
\paragraph{Unplugged} Member accounts which have no charging sessions associated with them. These members may correspond to users who only use home charging~\citep[which correspond to 53\% of respondents in ][]{Lee2020}, or PHEV owners who never charge their vehicles~\citep{Chakraborty2020}. Our industrial partner has also suggested that these could correspond to duplicate accounts, where users have forgotten account sign-in information, and have simply created a new account. Since the classification of users excludes the period of February 1st 2020 to June 30th 2021, users are classified as unplugged if their only charging sessions are during this time. 
    
\paragraph{Shared} Member accounts presumed to correspond to fleets or commercial vehicles. Several filters identify this case:
    \begin{itemize}
        \item Member accounts which are associated with multiple sessions occurring at the same time. 
        \item Member accounts which have been charged four or more times in a single day. This threshold value is much higher than the average charging rates of 18 to 20 sessions per year in \cite{vanDenHoed2013}, average 1.47 sessions per day in \cite{Axsen2015} (median of 1), average 1.17 sessions per day in \cite{Tal2018}, median 0.71 sessions per day in \cite{Visaria2022} (mean of 0.75, maximum of 4.29),  and average 1.1 sessions per day in \cite{Tal2020} (highest rate among all vehicle types). Additionally, while less commonly found in the literature,  in \cite{Tal2020} the average number of charging sessions on days in which charging occurred was reported to be below 1.66 sessions per day for all vehicle types. 
        \item Member accounts which have an average recharge rate over 500 kWh per month, or 6,000 kWh per year. This threshold is similar to the highest annual energy requirements reported in \cite{Tal2020} (6,565 kWh, followed by 5,043 kWh). While the energy efficiency varies by vehicle and climate, 6,000 kWh per year corresponds to over 30,000 kilometres of travel~\citep{vanDenHoed2013, Hardman2018}. By comparison, the average annual driving distances for EVs was 20,150 kilometres in \cite{Axsen2015}, 15,563 kilometres for non-Tesla EVs in \cite{Figenbaum2016}, and 23,367 kilometres for Tesla EVs in \cite{Figenbaum2016}.
    \end{itemize}
    
\paragraph{Rental} Member accounts which have less than 14 days between the creation of the account and the last charging session, and have at least three charging sessions during that time. While this corresponds to a reasonable charging rate, the short duration of the account suggests that this may be a temporary vehicle. As such, the charging behaviour may not correspond to general vehicle ownership, particularly if the rental vehicle is used for long-distance trips to cottages or other vacation destinations~\citep{Figenbaum2016, Figenbaum2019}. 

\paragraph{Private} All the accounts that are not labelled as unplugged, shared or rental are considered private.\\

\begin{figure}
    \centering
\includegraphics[width = 0.75\textwidth]{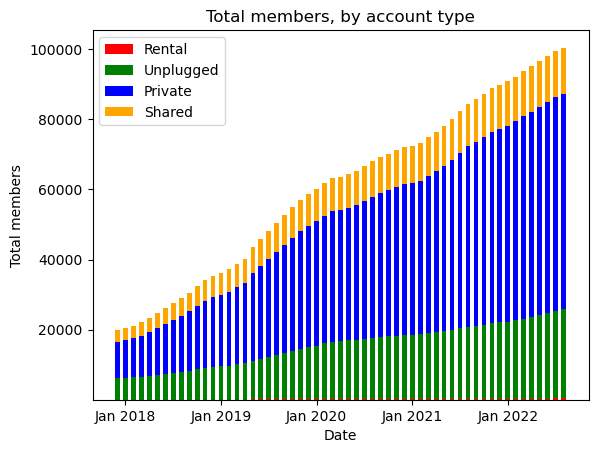}
    \caption{Number of total members, by date and account type.}
    \label{ThirdArticle:Data:MembersSpread}
\end{figure}

We note that a different classification of users was proposed in \cite{Helmus2020}, partly based on the types of charging sessions initiated by that user. However, this relies on (unavailable to us) vehicle information, and the resulting classifications do not isolate car sharing and taxis from private vehicles. As such, they are not of practical interest for our study.  The total number of members in each month and of each type is given in Figure~\ref{ThirdArticle:Data:MembersSpread}.

\begin{figure}
    \centering
\includegraphics[width = 0.75\textwidth]{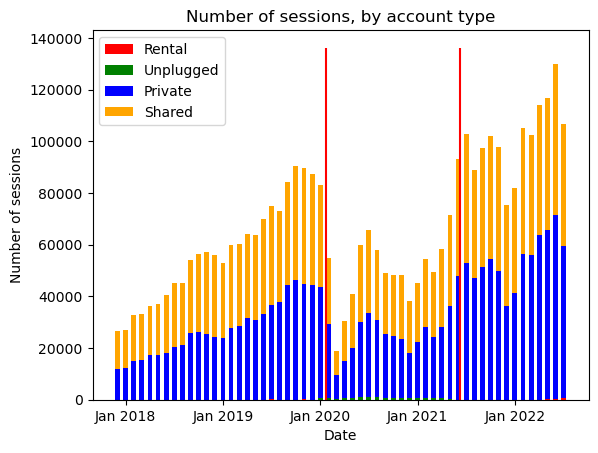}
    \caption{Number of sessions, by date and account type. The narrow red lines indicate the excluded period for the COVID 19 pandemic. }
    \label{ThirdArticle:Data:SessionSpread}
\end{figure}

In total, there are 2,873,345 valid recorded sessions by 120,952 non-guest member accounts at 81,239 outlets and 43,056 stations, which take place between January 2018 to August 2022. For each charging session (of level 2 or level 3), the unique identifiers for the charging outlet and the member account are recorded, along with the time that the charging connection started, the time the connection ended, and the amount of energy charged.  Additionally, for some level 3 charging outlets, the starting and ending state of charge are also recorded. However, we note that the price of charging is not available, nor is waiting time before charging.  The number of sessions in each month and by each type of member is given in Figure~\ref{ThirdArticle:Data:SessionSpread}.

\subsection{Data Processing}
\label{ThirdArticle:SectionDataProcessing}

For discrete choice model estimation and validation, we use the public charging network on the island of Montreal. As we consider intracity charging by private vehicle owners in this work, we only select members who have provided a postal code, and whose postal code lies within the island. Furthermore, since our industrial partner has advised us that there are many taxi accounts in Montreal, we exclude users if they have ever charged three or more times within the island in a day. These two conditions --postal code within Montreal, and maximum of two or fewer daily charging sessions on the island-- are applied in addition to the filters for private vehicle owners in Section~\ref{ThirdArticle:SectionDataDescription}. 

\begin{figure}
    \centering
    \includegraphics[scale = 0.5]{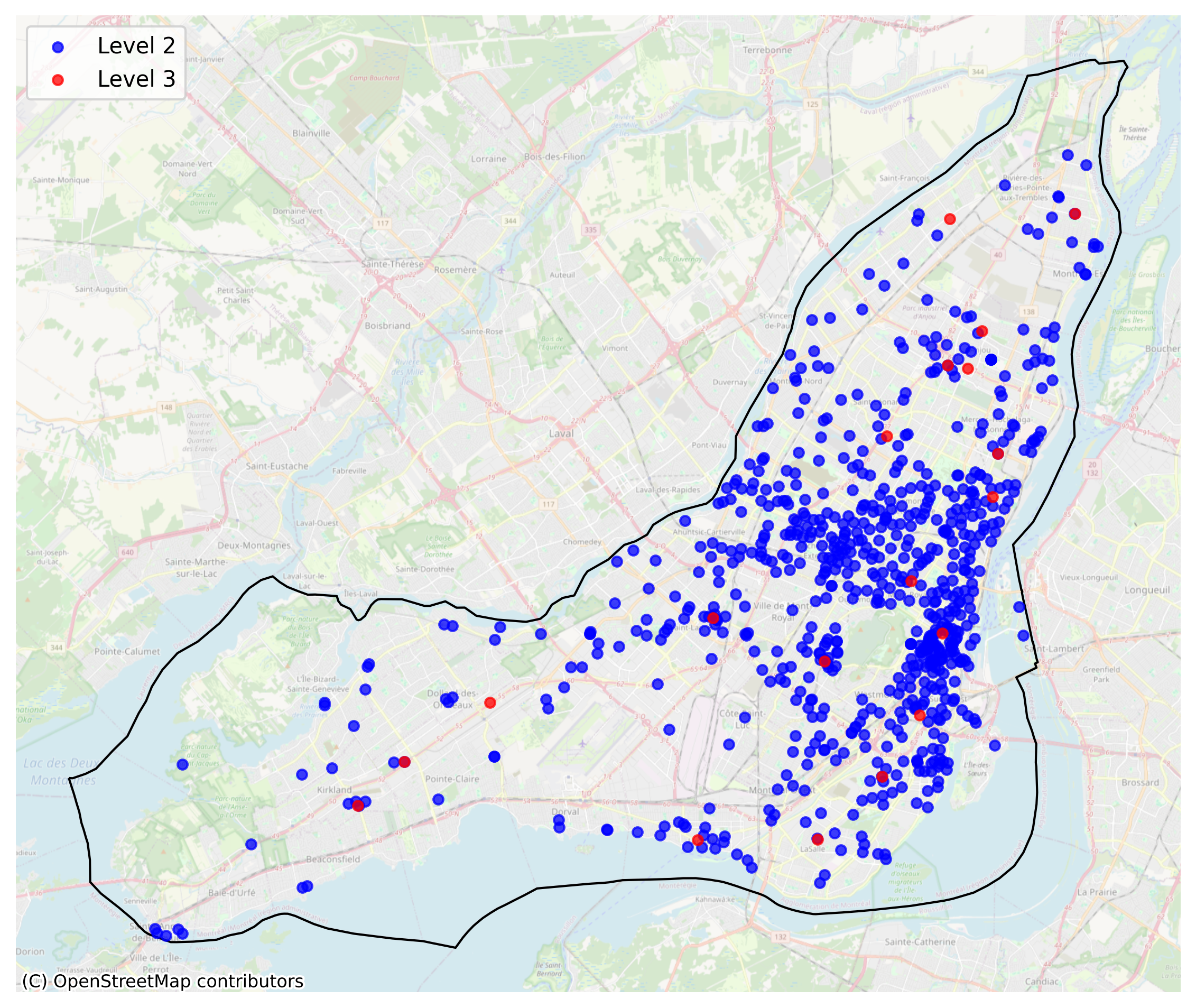}
    \caption{Public charging stations within the Island of Montreal.}
    \label{ThirdArticle:Data:SiteMap}
\end{figure}

As for the charging stations, they are excluded if they are private or semi-private stations~\citep[similar to ][]{Sun2016} 
, or if they were permanently closed before January 1st 2018. Additionally, while members can use their \emph{Circuit Électrique} membership cards at charging stations operated by \emph{ChargePoint}, this incurs an additional fee~\citep{FindAStation}. These stations have few recorded sessions, likely due to users using a \emph{ChargePoint} rather than a \emph{Circuit Électrique} membership card to avoid the additional fee. As such, to avoid incorrectly biasing the results, we exclude the charging stations operated by \emph{ChargePoint}. After filtering, there are 736 charging stations of level 2 charging stations, and 19 charging stations of level 3. The stations are displayed in Figure~\ref{ThirdArticle:Data:SiteMap}.

\subsection{Attribute Encoding}
\label{ThirdArticle:SectionAttributeEnconding}

To estimate model parameter values, we first determine values for each attribute. 
As the session data spans from January 2018 to August 2022, it is important to get historical data. For charging station attributes, we use the data directly from \emph{Circuit Électrique} as well as publicly available charging station information~\citep{FindAStation}. For network information and available amenities, we use OpenStreetMaps~\citep{OpenStreetMap} data, queried at the start of every month. We use attributes which have been examined, or for which similar attributes have been examined, previously in the literature. These attributes are:
\begin{description}
    \item[Distance:] The distance between the user and the charging station (or the deviation from their path between home and destination) is the most common attribute in charging station choice modelling, and has consistently been found to be a significant attribute for the user's choices~\citep[e.g. ][]{Yang2016, Lamontagne2023a, Ma2022, Visaria2022}. A shorter distance between the user's home and the charging station improves the \emph{accessibility} of the station as described in \cite{Philipsen2015, Philipsen2016}. We use two measures for these distances; First, we use OpenStreetMaps to calculate the road network distance between the user's home and the charging station. Second, we include a binary flag indicating whether the Euclidean (i.e. walking) distance between the charging station and the user's home is less than 400 metres, which indicates that the charging station may be convenient for charging near home~\citep{Lee2020}. Though \cite{Lee2020} proposes a distance of 300 metres, we adjust this threshold slightly to account for the distance errors due to having postal codes rather than precise addresses.  
    \item[Number of outlets:] Included in \cite{Lamontagne2023a, Visaria2022}, this attribute indicates the total number of outlets at each charging station. A higher number of outlets at the station increases the likelihood that one is available when needed, thus improving the \emph{reliability} of the station as described in \cite{Philipsen2015, Philipsen2016}. The number of outlets in each station is calculated from the \emph{Circuit Électrique} data based on the installation and closure dates of each outlet.
    \item[Is at a gas station:] The charging station being installed at a gas station has been found to be beneficial for users~\citep{Philipsen2015, Sun2016, Philipsen2016}. This is attributed in \cite{Philipsen2015} to habit, where the gas stations are those previously used by the users for their ICEVs. In our case, we use a binary variable to indicate whether or not the charging station is located at a gas station, and is determined by the presence of a gas station within 50 metres of the charging station in OpenStreetMaps.
    \item[Proximity to amenities:] In \cite{Anderson2018, Figenbaum2019}, it was reported that users prefer to do activities while their vehicle is charging, with activities differing between users and charging level. Various location types and activities have been proposed across the literature, which we discuss in detail below. For each amenity type, we calculate density measures by taking the logarithm of the number of amenities of that type in OpenStreetMaps within various threshold distances around the charging station. 

    \begin{description}
        \item[Restaurants, fast food:] The presence of restaurants was included as an attribute in the demand models of \cite{Luo2015, Cui2018}, though the value was not estimated. The highest tier of amenities in \cite{Ge2022, Visaria2022, Anderson2023} included access to the combination of washrooms, a restaurant, and free Wifi. As such it is unclear if the presence of a restaurant alone is sufficient. Restaurants are included as part of an ``Other" category in \cite{Anderson2018, Anderson2023}, while the location is not included in \cite{Philipsen2016}. Additionally, while less than $1$\% of respondents in \cite{EECA2021} indicated that they currently used charging stations near restaurants, caf\'es, or bars, 76\% of respondents indicated that they prefer to eat and drink while charging. In summary, while restaurants have been included in many surveys and reports, there are conflicting assessments as to the value of such amenities. For our estimation, we separate fast food restaurants from other types; Due to the difference in preparation time, fast food restaurants may be more beneficial for level 3 charging sessions. 
        \item[Shopping, supermarkets, shopping malls:] At an aggregate level, the presence of general shopping facilities was found to increase charging station usage for some users in \cite{Anderson2023}, while the demand models in \cite{Luo2015, Cui2018, Sheldon2019} includes separate attributes for supermarkets and shopping malls. Notably, in \cite{Sheldon2019}, both supermarkets and shopping malls were found to increase the likelihood that users select a station in comparison to locations near entertainment venues, with supermarkets having a larger impact than malls. The overall benefit of shopping is consistent with the findings in \cite{Anderson2018}, where 18.4\% of the placed charging outlets were designated for shopping purposes, and in \cite{Philipsen2016}, where shopping was rated as important for users. In practice, \cite{Anderson2023} report that around 20-32\% of respondents charged at least once a month near shopping locations. In \cite{Axsen2015}, shopping malls were reported as the most frequent location type for public charging (though they note that public charging was overall infrequently used), while \cite{Figenbaum2016} report that around 7\% of BEVs recharge at least weekly at shopping centres and similar commercial locations. In \cite{EECA2021}, only 3\% and 2\% of respondents indicated that they recharged near supermarkets and shopping malls respectively, while 55\% of respondents said they liked to go shopping while recharging their vehicle. In summary, the presence of shopping facilities (and, where applicable, supermarkets and shopping centres) are generally seen as important by users, though they may be used in practice less often than anticipated. For our estimation, we consider supermarkets and shopping malls separately from other types of shopping. For supermarkets, this is due to the habitual and more frequent usage~\citep[as discussed in ][]{Philipsen2015}. On the contrary, for shopping malls (a category which also includes the public marketplaces in Montreal), this is due to the presence of multiple types of different shopping in one location, often not individually identified in OpenStreetMaps. As such, using a separate attribute for shopping malls can better capture the more diverse shopping opportunities available.
        \item[Leisure, sports:] Combined into one category, sport locations (such as tracks or arenas) and leisure locations (such as cinemas or dog parks) were considered in \cite{Anderson2023} and found to be important for some users. In \cite{Sheldon2019}, entertainment venues were used as the baseline location type for other types, and found that gyms and sports facilities were dispreferred compared with this baseline, with high variability among users.  In \cite{Anderson2018}, around 18.3\% of the placed charging outlets were designated for leisure purposes, almost identical to those for shopping purposes. By contrast, in \cite{Philipsen2016},  leisure locations were deemed less important than shopping. Likewise, in \cite{Anderson2023} report that only around 4\% of respondents charged at least once a month near leisure locations.   This is consistent with the findings in \cite{EECA2021}, who report that 1\% of respondents charged in the ``Other'' category (which includes entertainment venues and sports facilities), while 13\% of respondents indicated that they would like to use such facilities while charging. In summary, the use of charging stations near entertainment venues or sports facilities is typically lower than other amenities, though is important for some users. As with fast food restaurants, we separate sports and leisure locations due to potential variations in the duration of these activities. 
    \end{description}
\end{description}
The maximum value 
for each of the attributes for each charging level are presented in Table~\ref{ThirdArticle:TableAttributeMaximums}. Additionally, we present the associated notation $x_{k}$ for each attribute $k$ in the model. 

We note that some common attributes from the literature have not been included, either because they are not appropriate for our case study or because accurate values can not be determined: 
\begin{description}
    \item[Charging power and duration:] The power supplied by the charging station or, similarly, the duration of charging has been shown to have an impact on the selection of charging stations, with higher power or lower duration being preferred~\citep[e.g.][]{Yang2016, Ge2022, Visaria2022}. However, in Quebec, the charging power can be set independently for each charging outlet, and as such may not be consistent within each charging station.  Additionally, while over 95\% of the level 2 outlets have a charging power of 7.2 kW~\citep{FindAStation}, the charging power for level 3 outlets varies between 24-100 kW with different power levels deliberately placed at the same stations~\citep{CostOfCharging}. 
    \item[Charging price:] As expected, users have been found to avoid charging stations which have higher costs compared to others~\citep{Ge2022, Visaria2022, Wang2021}.  However, while the cost of charging in Quebec can vary slightly by location and provider, it is subject to governmental regulation, and is determined based on the charging power of the outlet, the charging power drawn by the vehicle, and state of charge of the vehicle~\citep{CostOfCharging}. As a consequence, the charging price can vary not only within a charging station but based on the charging profile of each individual vehicle. In Montreal, over 95\% of the level 2 charging stations have a cost of 1\$ per hour which, given the median duration of a level 2 session, results in a per-session cost of around 1.93\$~\citep{FindAStation}. By contrast, the cost for level 3 charging varies between 6.39-15.73\$ per hour, with the majority having a cost of 12.77\$ per hour (resulting in a per-session median cost of around 4.90\$). 
\end{description}

Charging stations are considered to be available for users if at least one charging outlet is installed at the time of the charging, and if a path is found between the users' home and the charging station.

\begin{table}[]
    \centering
\begin{tabular}{l l r r}
\toprule
Attribute name & Notation & Level 2 & Level 3
\\
\midrule
Network distance to station (km)  & $x_\text{dist}$ & 98.658 & 49.668 \\
Is the station within a short walk & $x_\text{isWalkHome}$ & 1 & 1\\
Number of outlets & $x_\text{outlets}$ & 13 & 5 \\
Is at a gas station & $x_\text{isGas}$ & 1 & 1\\
Restaurant density & $x_\text{rest}$ & 4.635 & 3.638 \\
Fast food density & $x_\text{ff}$ & 4.060 & 2.079 \\
Shopping facility density & $x_\text{shop}$ & 5.257 & 4.007 \\
Supermarket density & $x_\text{sm}$ &  2.303 & 0.693\\
Shopping mall density & $x_\text{mall}$ & 1.099 & 0.693 \\
Leisure density & $x_\text{leis}$ & 3.526 & 3.091\\
Sports facility density & $x_\text{sport}$ & 2.833 & 2.565 \\
\bottomrule
\end{tabular}
    \caption{Maximum value of attributes for level 2 and level 3 charging sessions.}
    \label{ThirdArticle:TableAttributeMaximums}
\end{table}
\section{Discrete Choice Modeling}
\label{ThirdArticle:SectionMethods}

Similar to many of the works described in Section~\ref{ThirdArticle:SectionLiteratureReview}, we use Random Utility Maximisation (RUM) discrete choice models for predicting the selection of charging station by users. As mentioned previously, the reasons why we use these models are twofold: \emph{(i)} the disaggregate nature of our data and \emph{(ii)} the existing literature permitting their direct integration into optimisation models. 

In Section~\ref{ThirdArticle:SectionModelSpecification}, we provide the specifications of the models as to the interactions of the attributes and attribute levels. While in Section~\ref{ThirdArticle:SectionMethodsValidation}, we describe the validation process used to compare the various models. 

\subsection{Model Specifications}
\label{ThirdArticle:SectionModelSpecification}

For the sake of simplicity, in this section, we only detail aspects of our models which are relevant for the discussion of results. For the model properties and derivations, we refer to \cite{McFadden1974, Revelt1998, McFadden2000} and \cite{Train2002}. We provide the specifications for both MNL and MXL models, with the latter considering panel effects. As noted in, e.g., \cite{Philipsen2016, Sheldon2019, EECA2021}, the evaluation of location types and amenities vary significantly among individuals, making the MXL well-suited for this application. The high number of alternatives (i.e. charging stations), the total number of users, and the small number of observations for many users make it computationally infeasible to compute alternative-specific or user-specific constants. Additionally, while such constants are useful in predicting individual behaviour, they are not practical for large-scale network operation. As such, we adopt a random parameters specification for our MXL model. 

In our application, each user $i$ in a set of users $N$ selects a charging station alternative $j$ in a set of alternatives $M$ based on the attributes $k = 1, \dotsc, K$ given in Section~\ref{ThirdArticle:SectionAttributeEnconding}. A \emph{utility function} 
$u_{ij} = V_{ij} + \epsilon_{ij}$ is associated with each alternative $j\in M$ and user $i\in N$, which combines a function of \emph{observable} attributes $V_{ij}$ known by the modeller and \emph{unobservable} attributes $\epsilon$ which are known only to the user. In both the MNL and MXL models, we assume that our observable utility takes the standard form which is linear-in-parameters
$$
V_{ij} = \sum_{k =1}^K \beta_{k} x_{ijk}.
$$
We recall that $x_{ijk}$ denotes the value of attribute $k$ for user $i$ and charging station $j$. We note that the parameters $\beta_{k}$ are not indexed by $i$ in this formulation as the dependence of the observable utility on the user appears via the attributes $x_{ijk}$.
In the case of our MXL model, the mean $\beta_k^{\mu}$ and standard deviation $\beta_k^{\sigma}$ of Normal distributions are estimated for each attribute $k$, while only the mean value is estimated for MNL. 

We assume that users are either charging near their home or charging while performing daily activities~\citep{Hardman2018}, with different attributes considered for each case. We note that \cite{Hardman2018} also list workplace charging and charging during long-distance travel, which do not apply in our situation. 
Since the charging session dataset does not include the paths travelled by the user~\citep[contrary to ][]{Sun2016}, 
it is not possible to determine whether the user is charging near their home or while performing daily activities (most notably, their activity may be near their home). As such, we use a cutoff threshold based on the distance from home, with charging sessions within 1.5 km (network distance) being considered as charging near home and charging sessions outside of 1.5 km being considered as charging near an activity.  This threshold distance was selected based on the attributes for walking proximity to home, and is also similar to the maximum path deviation distance of 1.75 km for charging by private vehicle users in \cite{Sun2016}.

More specifically, the observable utility within 1.5 km is given by 
\begin{equation*}
    V = \beta_\text{distNear} x_\text{dist} + \beta_\text{outletsNear} x_\text{outlets} + \beta_\text{isWalkHome} x_\text{isWalkHome},
\end{equation*}
where the indices for user $i$ and station $j$ have been omitted for ease of reading. The parameters $\beta_\text{distNear}, \beta_\text{outletsNear}$, and $\beta_\text{isWalkHome}$ are respectively those for the distance to home, the number of outlets, and the station being a short walk from home. 

When outside of 1.5 km, the observable utility is given by
\begin{align*}
    V = \beta_\text{distFar} &x_\text{dist} + \beta_\text{outletsFar} x_\text{outlets} + \beta_\text{isGas} x_\text{isGas} + \beta_\text{leis} x_\text{leis} + \beta_\text{sport} x_\text{sport} \\
    &+ \beta_\text{sm} x_\text{sm} + \beta_\text{shop} x_\text{ishop} + \beta_\text{mall} x_\text{mall} + \beta_\text{rest} x_\text{rest} + \beta_\text{ff} x_\text{f},
\end{align*}
where the indices $i$ and $j$ have been omitted as before. The parameters $\beta_\text{distFar}$ and $\beta_\text{outletsFar}$, as well as the attribute levels  $x_\text{dist}$ and $x_\text{outlets}$ are equivalent to the prior case. The parameter and attribute level $\beta_\text{isGas}$  and $x_\text{isGas}$ are those for whether station $j$ is located at a gas station. The parameters $\beta_\text{leis}, \beta_\text{sport}, \beta_\text{sm}, \beta_\text{shop}, \beta_\text{mall}, \beta_\text{rest}$, and $\beta_\text{ff}$ are those for the proximity to amenities, respectively for leisure facilities, sports facilities, supermarkets, shopping facilities, shopping malls, restaurants, and fast food locations.

\subsection{Estimation and Validation}
\label{ThirdArticle:SectionMethodsValidation}
To estimate, validate, and compare the MNL and MXL models, we employ internal five-fold cross-validation~\citep{Parady2021}. More specifically, due to the presence of panel data, we use a grouped sampling approach based on the users to prevent \emph{data leakage}~\citep[see, e.g.,][for discussions about data leakage related to panel data]{Hillel2020, Hillel2021} 
In this approach, users are randomly divided into five groups, while keeping a similar number of observations in each group. The observations from some groups are used for estimation, while the others are kept for validation (the specifics of which observations are different for level 2 and level 3 charging, and are given below). For validation, only the last observation (i.e. the latest observation in time) is used. 
Additionally, in original testing, we observed that outliers in the distance attribute (that is, sessions with much higher distance between the user's home and the charging station) negatively impacted the estimation process. As such, we have removed the sessions with the distance attribute in the top 2\% (over 28.991 km for level 2 and 27.838 km for level 3). The original distributions for the distance can be found in~\ref{ThirdArticle:AppendixAttributes} as Figure~\ref{ThirdArticle:FigureDistributionDistanceToHomeSelected}.

For level 3 charging, each group contains 214 users and approximately 830 observations. For each fold, the estimation set is composed of four out of the five groups, with the last observation of the users in the remaining group reserved for validation. This results in a roughly 94-6 split of observations between estimation and validation 
For level 2 charging, each group contains approximately 1,180 users and 11,000 observations. However, the high number of observations and alternatives make it intractable to estimate a model even for one group. As such, for each fold, the estimation set samples 5,000 observations from one group, with the last observation of the users in the remaining groups reserved for validation. This results in a roughly 51-49 split of observations between estimation and validation.

In addition to parameter values, we report a series of standard performance indicators about the results in both the estimation and validation sets, such as the log-likelihood (of both the final parameter values and with all parameter values set to 0), $\rho$, $\bar{\rho}^2$, Akaike Information Criterion, and Bayesian Information Criterion~\citep{Parady2021, Bierlaire2023}.  Moreover, for the validation set, we present some statistics relating to the distribution of the probability of the chosen alternatives for the estimated model, denoted ``DPSA, final'' in the tables below. The mean of this distribution -- given by 
the estimated probability of the observed alternatives averaged over individuals
-- corresponds to the \emph{fitting factor} described in \cite{Parady2021}. Since the number of available alternatives is not constant throughout the observations, we also present the distribution of the choice probabilities corresponding to the uniform distribution, denoted ``DPSA, null'' in the tables below. Higher values are preferred for the log-likelihoods, $\rho$, $\bar{\rho}^2$, and DPSA, while lower values are desirable for the Akaike Information Criterion and the Bayesian Information Criterion.

\section{Estimation and Validation Results}
\label{ThirdArticle:SectionResults}
In this section, we present the results of the estimation and validation processes. Model estimation is performed with Biogeme version 3.2.12~\citep{Bierlaire2023}. A total of 1,000 draws are used for the simulations, where applicable.  The results for level 2 and level 3 charging are presented, respectively, in Sections~\ref{ThirdArticle:ResultsL2} and~\ref{ThirdArticle:ResultsL3}, while a discussion of these results is reserved for Section~\ref{ThirdArticle:SectionDiscussion}.
We refer to \ref{ThirdArticle:AppendixAdditionalResults} for more detailed results.

\subsection{Level 2}
\label{ThirdArticle:ResultsL2}

Parameter ratio 
values (dividing with the estimate of $\beta_{\text{isWalkHome}}^\mu$) for both MNL and MXL models in each fold are presented in Table~\ref{ThirdArticle:TableParametersL2Combined}, while the averages of these ratios across folds are presented in Table~\ref{ThirdArticle:TableL2RatioTest}. Performance indicators for MNL and MXL models in the estimation sets are, respectively, presented in Tables~\ref{ThirdArticle:TableL2EstimationMNL} and~\ref{ThirdArticle:TableL2EstimationMXL}. We recall that the estimation results are not comparable between the MNL and MXL models, due to the panel effect in the MXL model. In Table~\ref{ThirdArticle:TableL2ValidationCombined}, we present the performance indicators of both models in the validation sets, which can be compared. 

\begin{table}
\caption[Parameter ratio values across models and folds for level 2 charging]{Parameter ratio values across models and folds for level 2 charging. *** indicates significance at 1\% level, ** significance at 5\% level, and * significance at 10\% level.}
\label{ThirdArticle:TableParametersL2Combined}
\centering
\resizebox{1.0\textwidth}{!}{
\begin{tabular}{lllllllllll}
\toprule
Fold & \multicolumn{2}{c}{0} & \multicolumn{2}{c}{1} & \multicolumn{2}{c}{2} & \multicolumn{2}{c}{3} & \multicolumn{2}{c}{4} \\
Model & MNL & MXL & MNL & MXL & MNL & MXL & MNL & MXL & MNL & MXL \\
\midrule
$\beta_\text{distNear}^{\mu}$ & 0.1885*** & -0.2592*** & 0.3046*** & -0.1756 & 0.1822*** & -0.5973*** & 0.1790*** & -0.1857*** & 0.2378*** & -0.7890*** \\
$\beta_\text{distNear}^{\sigma}$ & - & 0.7842*** & - & 0.7804*** & - & 1.0277*** & - & 0.5848*** & - & 1.3037*** \\
$\beta_\text{distFar}^{\mu}$ & -0.0612*** & -0.1688*** & -0.0549*** & -0.1091*** & -0.0459*** & -0.1019*** & -0.0546*** & -0.0862*** & -0.0558*** & -0.1376*** \\
$\beta_\text{distFar}^{\sigma}$ & - & 0.1202*** & - & 0.0871*** & - & 0.0811*** & - & 0.0655*** & - & 0.1036*** \\
$\beta_\text{ff}^{\mu}$ & -0.0216*** & -0.1107*** & -0.0527*** & -0.2021*** & -0.0218*** & -0.1551*** & -0.0421*** & -0.0918*** & -0.0368*** & -0.1642*** \\
$\beta_\text{ff}^{\sigma}$ & - & 0.24*** & - & 0.3404*** & - & 0.2791*** & - & 0.1468*** & - & 0.2486** \\
$\beta_\text{leis}^{\mu}$ & -0.0684*** & -0.1353*** & -0.0434*** & -0.0684* & -0.078*** & -0.08*** & -0.0584*** & -0.0505*** & -0.0596*** & -0.0239 \\
$\beta_\text{leis}^{\sigma}$ & - & 0.5218*** & - & 0.3816*** & - & 0.3428*** & - & 0.2461*** & - & 0.4789*** \\
$\beta_\text{mall}^{\mu}$ & 0.1322*** & -0.0437 & 0.0837*** & -0.0415 & 0.0627*** & -0.1288* & 0.2276*** & 0.1169*** & 0.1826*** & 0.0970 \\
$\beta_\text{mall}^{\sigma}$ & - & 0.8658*** & - & 0.6641*** & - & 0.6802*** & - & 0.3395*** & - & 0.7178*** \\
$\beta_\text{outletsNear}^{\mu}$ & 0.0938*** & 0.0542 & 0.0876*** & 0.0527 & 0.0773*** & 0.1293*** & 0.0687*** & 0.0273 & 0.0640*** & 0.1002** \\
$\beta_\text{outletsNear}^{\sigma}$ & - & 0.3015*** & - & 0.4774*** & - & 0.3153*** & - & 0.2549*** & - & 0.5203*** \\
$\beta_\text{outletsFar}^{\mu}$ & 0.0384*** & 0.0510*** & 0.031*** & 0.0283*** & 0.0245*** & 0.0088 & 0.0253*** & 0.0289*** & 0.0238*** & 0.0340** \\
$\beta_\text{outletsFar}^{\sigma}$ & - & 0.0639*** & - & 0.0811*** & - & 0.1002*** & - & 0.0374*** & - & 0.0862*** \\
$\beta_\text{rest}^{\mu}$ & 0.1225*** & 0.2428*** & 0.1348*** & 0.2692*** & 0.0964*** & 0.2297*** & 0.1039*** & 0.1200*** & 0.1264*** & 0.3159*** \\
$\beta_\text{rest}^{\sigma}$ & - & 0.3552*** & - & 0.2810*** & - & 0.3143*** & - & 0.1653*** & - & 0.3541*** \\
$\beta_\text{shop}^{\mu}$ & -0.0479*** & -0.1008*** & -0.0362*** & -0.1041*** & -0.0052 & -0.061* & -0.0293*** & -0.0496*** & -0.0502*** & -0.1146*** \\
$\beta_\text{shop}^{\sigma}$ & - & 0.2383*** & - & 0.1725*** & - & 0.2529*** & - & 0.1276*** & - & 0.1332*** \\
$\beta_\text{sport}^{\mu}$ & 0.0315*** & -0.0083 & 0.0153* & -0.0162 & 0.0355*** & 0.0621* & 0.0422*** & 0.0156 & 0.0197*** & 0.0190 \\
$\beta_\text{sport}^{\sigma}$ & - & 0.4469*** & - & 0.4508*** & - & 0.4082*** & - & 0.2299*** & - & 0.4444*** \\
$\beta_\text{sm}^{\mu}$ & 0.0422*** & 0.0607 & 0.0411*** & 0.0036 & 0.0072 & -0.0055 & 0.0356*** & -0.0063 & 0.0742*** & 0.1036** \\
$\beta_\text{sm}^{\sigma}$ & - & 0.4553*** & - & 0.4723*** & - & 0.3752*** & - & 0.2977*** & - & 0.5632*** \\
$\beta_\text{isGas}^{\mu}$ & -0.0931*** & -0.3351*** & -0.1072*** & -0.2207 & 0.0749*** & -0.6976*** & -0.1207*** & -0.4598*** & -0.1128*** & -0.3961 \\
$\beta_\text{isGas}^{\sigma}$ & - & 0.3803*** & - & 0.3635** & - & 0.8264*** & - & 0.5106*** & - & 0.5407 \\
$\beta_\text{isWalkHome}^{\mu}$ & 1.0000*** & 1.0000*** & 1.0000*** & 1.0000 & 1.0*** & 1.0000*** & 1.0000*** & 1.0000*** & 1.0000*** & 1.0000*** \\
$\beta_\text{isWalkHome}^{\sigma}$ & - & 1.5693*** & - & 1.4397 & - & 1.2911*** & - & 0.8739*** & - & 1.7154*** \\
\bottomrule
\end{tabular}}
\end{table}

\begin{table}[]
    \centering
\begin{tabular}{lrr}
\toprule
 & MNL & MXL \\
\midrule
$\beta_\text{distNear}^{\mu}$ & 0.2184 & -0.4014 \\
$\beta_\text{distNear}^{\sigma}$ & - & 0.8962 \\
$\beta_\text{distFar}^{\mu}$ & -0.0545 & -0.1207 \\
$\beta_\text{distFar}^{\sigma}$ & - & 0.0591 \\
$\beta_\text{ff}^{\mu}$ & -0.0350 & -0.1448 \\
$\beta_\text{ff}^{\sigma}$ & - & 0.2510 \\
$\beta_\text{leis}^{\mu}$ & -0.0616 & -0.0716 \\
$\beta_\text{leis}^{\sigma}$ & - & 0.3942 \\
$\beta_\text{mall}^{\mu}$ & 0.1378 & -0.0000 \\
$\beta_\text{mall}^{\sigma}$ & - & 0.6535 \\
$\beta_\text{outletsNear}^{\mu}$ & 0.0783 & 0.0727 \\
$\beta_\text{outletsNear}^{\sigma}$ & - & 0.3739 \\
$\beta_\text{outletsFar}^{\mu}$ & 0.0286 & 0.0302 \\
$\beta_\text{outletsFar}^{\sigma}$ & - & 0.0337 \\
$\beta_\text{rest}^{\mu}$ & 0.1168 & 0.2355 \\
$\beta_\text{rest}^{\sigma}$ & - & 0.2940 \\
$\beta_\text{shop}^{\mu}$ & -0.0338 & -0.0860 \\
$\beta_\text{shop}^{\sigma}$ & - & 0.1849 \\
$\beta_\text{sport}^{\mu}$ & 0.0288 & 0.0144 \\
$\beta_\text{sport}^{\sigma}$ & - & 0.3960 \\
$\beta_\text{sm}^{\mu}$ & 0.0401 & 0.0312 \\
$\beta_\text{sm}^{\sigma}$ & - & 0.4327 \\
$\beta_\text{isGas}^{\mu}$ & -0.0718 & -0.4219 \\
$\beta_\text{isGas}^{\sigma}$ & - & 0.5243 \\
$\beta_\text{isWalkHome}^{\mu}$ & 1.0000 & 1.0000 \\
$\beta_\text{isWalkHome}^{\sigma}$ & - & 1.3779 \\
\bottomrule
\end{tabular}
    \caption[Average parameter ratios for level 2 charging]{Parameter ratios for MNL and MXL models for level 2 charging, average across folds. }
    \label{ThirdArticle:TableL2RatioTest}
\end{table}

\begin{table}[]
    \centering
        \resizebox{1.0\textwidth}{!}{
\begin{tabular}{lrrrrr}
\toprule
 & 0 & 1 & 2 & 3 & 4 \\
\midrule
Null log-likelihood & -31030.6143 & -30924.5374 & -30876.1183 & -30911.5930 & -30881.7735 \\
Final log-likelihood & -26141.5274 & -26029.6193 & -25712.7924 & -25368.9420 & -25445.0283 \\
$\rho$ & 0.1576 & 0.1583 & 0.1672 & 0.1793 & 0.1761 \\
$\bar{\rho}^2$ & 0.1571 & 0.1579 & 0.1668 & 0.1789 & 0.1756 \\
Akaike Information Criterion & 52309.0547 & 52085.2387 & 51451.5848 & 50763.8841 & 50916.0566 \\
Bayesian Information Criterion & 52393.7783 & 52169.9622 & 51536.3083 & 50848.6076 & 51000.7801 \\
\bottomrule
\end{tabular}
}
    \caption{Performance indicators for level 2 charging and the MNL model, estimation sets.}
    \label{ThirdArticle:TableL2EstimationMNL}
\end{table}

\begin{table}[]
    \centering
        \resizebox{1.0\textwidth}{!}{
\begin{tabular}{lrrrrr}
\toprule
 & 0 & 1 & 2 & 3 & 4 \\
\midrule
Null log-likelihood & -31030.6143 & -30924.5374 & -30876.1183 & -30911.5930 & -30881.7735 \\
Final log-likelihood & -22557.9877 & -22763.1175 & -22747.6526 & -22622.8858 & -22534.8722 \\
$\rho$ & 0.2730 & 0.2639 & 0.2633 & 0.2681 & 0.2703 \\
$\bar{\rho}^2$ & 0.2722 & 0.2631 & 0.2624 & 0.2673 & 0.2694 \\
Akaike Information Criterion & 45167.9754 & 45578.2350 & 45547.3051 & 45297.7715 & 45121.7444 \\
Bayesian Information Criterion & 45337.4224 & 45747.6820 & 45716.7521 & 45467.2186 & 45291.1914 \\
\bottomrule
\end{tabular}
}
    \caption{Performance indicators for level 2 charging and the MXL model, estimation sets.}
    \label{ThirdArticle:TableL2EstimationMXL}
\end{table}

\begin{table}[]
    \centering
        \resizebox{1.0\textwidth}{!}{
\begin{tabular}{lrrrrrrrrrr}
\toprule
Fold & \multicolumn{2}{c}{0} & \multicolumn{2}{c}{1} & \multicolumn{2}{c}{2} & \multicolumn{2}{c}{3} & \multicolumn{2}{c}{4} \\
Model & MNL & MXL & MNL & MXL & MNL & MXL & MNL & MXL & MNL & MXL \\
\midrule
Null log-likelihood & -29429.3387 & -29429.3387 & -29375.5664 & -29375.5664 & -29382.9599 & -29382.9599 & -29446.8929 & -29446.8929 & -29377.0021 & -29377.0021 \\
Final log-likelihood & -25743.6379 & -26146.6463 & -25479.0207 & -25764.3990 & -25597.1986 & -25814.6837 & -25711.0459 & -25976.8103 & -25528.0620 & -25680.6044 \\
$\rho$ & 0.1252 & 0.1115 & 0.1326 & 0.1229 & 0.1288 & 0.1214 & 0.1269 & 0.1178 & 0.1310 & 0.1258 \\
$\bar{\rho}^2$ & 0.1248 & 0.1107 & 0.1322 & 0.1220 & 0.1284 & 0.1206 & 0.1264 & 0.1170 & 0.1306 & 0.1249 \\
Akaike Information Criterion & 51513.2757 & 52345.2926 & 50984.0415 & 51580.7980 & 51220.3973 & 51681.3674 & 51448.0919 & 52005.6207 & 51082.1239 & 51413.2089 \\
Bayesian Information Criterion & 51596.9237 & 52512.5887 & 51067.6640 & 51748.0431 & 51304.0227 & 51848.6182 & 51531.7399 & 52172.9167 & 51165.7437 & 51580.4484 \\
DPSA null, min & 0.0014 & 0.0014 & 0.0014 & 0.0014 & 0.0014 & 0.0014 & 0.0014 & 0.0014 & 0.0014 & 0.0014 \\
DPSA null, mean & 0.0017 & 0.0017 & 0.0017 & 0.0017 & 0.0017 & 0.0017 & 0.0017 & 0.0017 & 0.0017 & 0.0017 \\
DPSA null, max & 0.0037 & 0.0037 & 0.0038 & 0.0038 & 0.0038 & 0.0038 & 0.0038 & 0.0038 & 0.0038 & 0.0038 \\
DPSA final, min & 0.0000 & 0.0001 & 0.0000 & 0.0001 & 0.0000 & 0.0001 & 0.0000 & 0.0001 & 0.0000 & 0.0001 \\
DPSA final, mean & 0.0209 & 0.0245 & 0.0216 & 0.0228 & 0.0236 & 0.0204 & 0.0222 & 0.0279 & 0.0231 & 0.0196 \\
DPSA final, max & 0.5853 & 0.4673 & 0.6238 & 0.3494 & 0.5847 & 0.3881 & 0.6147 & 0.5312 & 0.6432 & 0.3690 \\
DPSA average, min & 0.0000 & 0.0000 & 0.0000 & 0.0000 & 0.0000 & 0.0000 & 0.0000 & 0.0000 & 0.0000 & 0.0000 \\
DPSA average, mean & 0.0028 & 0.0028 & 0.0029 & 0.0029 & 0.0030 & 0.0030 & 0.0032 & 0.0032 & 0.0031 & 0.0031 \\
DPSA average, max & 0.0216 & 0.0216 & 0.0132 & 0.0132 & 0.0126 & 0.0126 & 0.0184 & 0.0184 & 0.0184 & 0.0184 \\
DPSA closest, min & 0.0000 & 0.0000 & 0.0000 & 0.0000 & 0.0000 & 0.0000 & 0.0000 & 0.0000 & 0.0000 & 0.0000 \\
DPSA closest, mean & 0.0185 & 0.0185 & 0.0172 & 0.0172 & 0.0183 & 0.0183 & 0.0167 & 0.0167 & 0.0179 & 0.0179 \\
DPSA closest, max & 1.0000 & 1.0000 & 1.0000 & 1.0000 & 1.0000 & 1.0000 & 1.0000 & 1.0000 & 1.0000 & 1.0000 \\
\bottomrule
\end{tabular}}
    \caption{Performance indicators for level 2 charging, validation sets.}
    \label{ThirdArticle:TableL2ValidationCombined}
\end{table}

\subsection{Level 3}
\label{ThirdArticle:ResultsL3}


Parameter values for both MNL and MXL models are presented in Table~\ref{ThirdArticle:TableParametersL3Combined}, while the average of the parameter ratios are presented in Table~\ref{ThirdArticle:TableL3RatioTest}. Performance indicators for MNL and MXL models in the estimation sets are, respectively, presented in Tables~\ref{ThirdArticle:TableEstimationL3MNL} and~\ref{ThirdArticle:TableEstimationL3MXL}. As with level 2 charging, the estimation results are not comparable between the MNL and MXL models. In Table~\ref{ThirdArticle:TableValidationL3Combined}, we present the performance indicators of both models in the validation sets, which can be compared.

\begin{table}
\centering
\resizebox{1.0\textwidth}{!}{
\begin{tabular}{lllllllllll}
\toprule
Fold & \multicolumn{2}{c}{0} & \multicolumn{2}{c}{1} & \multicolumn{2}{c}{2} & \multicolumn{2}{c}{3} & \multicolumn{2}{c}{4} \\
Model & MNL & MXL & MNL & MXL & MNL & MXL & MNL & MXL & MNL & MXL \\
\midrule
$\beta_\text{distNear}^{\mu}$ & 0.5640*** & 0.1820 & 0.8865*** & -0.0371 & 0.5504*** & -0.0004 & 0.7996*** & 0.1813 & 0.8894*** & 0.3881** \\
$\beta_\text{distNear}^{\sigma}$ & - & 0.2096*** & - & 1.0174*** & - & 0.5111*** & - & 0.5584*** & - & 0.6135*** \\
$\beta_\text{distFar}^{\mu}$ & -0.1089*** & -0.1457*** & -0.1346*** & -0.2025*** & -0.1231*** & -0.1065*** & -0.1710*** & -0.2813*** & -0.1109*** & -0.1517*** \\
$\beta_\text{distFar}^{\sigma}$ & - & 0.1070*** & - & 0.1307*** & - & 0.0764*** & - & 0.1928*** & - & 0.1092*** \\
$\beta_\text{ff}^{\mu}$ & -0.0363 & -0.0022 & -0.116*** & -0.2568*** & 0.0676 & 0.0196 & 0.0315 & 0.0220 & 0.0213 & -0.0827 \\
$\beta_\text{ff}^{\sigma}$ & - & 0.4933*** & - & 0.7711*** & - & 0.2917*** & - & 0.7814*** & - & 0.6147*** \\
$\beta_\text{leis}^{\mu}$ & 0.1844*** & 0.1872*** & 0.3106*** & 0.3934*** & 0.1695*** & 0.1608*** & 0.2783*** & 0.2795*** & 0.1877*** & 0.2279*** \\
$\beta_\text{leis}^{\sigma}$ & - & 0.3820*** & - & 0.6995*** & - & 0.4456*** & - & 0.9465*** & - & 0.5517*** \\
$\beta_\text{mall}^{\mu}$ & 0.3301*** & 0.2343*** & 0.4297*** & 0.4599*** & 0.5037*** & 0.2706*** & 0.6088*** & 0.3935** & 0.4088*** & 0.3852*** \\
$\beta_\text{mall}^{\sigma}$ & - & 0.8221*** & - & 1.3356*** & - & 0.5464*** & - & 1.568*** & - & 0.8198*** \\
$\beta_\text{outletsNear}^{\mu}$ & -0.0247 & 0.0181 & -0.0948 & 0.2328* & 0.0689 & 0.3015*** & 0.0802 & 0.5312*** & -0.1816*** & -0.0252 \\
$\beta_\text{outletsNear}^{\sigma}$ & - & 0.4024*** & - & 0.6039*** & - & 0.5134*** & - & 1.2447*** & - & 0.3516*** \\
$\beta_\text{outletsFar}^{\mu}$ & 0.1019*** & 0.004 & 0.1085*** & -0.0104 & 0.1077*** & -0.0078 & 0.1246*** & 0.0115 & 0.1188*** & -0.0006 \\
$\beta_\text{outletsFar}^{\sigma}$ & - & 0.1854*** & - & 0.4369*** & - & 0.1934*** & - & 0.4083*** & - & 0.2527*** \\
$\beta_\text{rest}^{\mu}$ & -0.0748*** & -0.0931** & -0.0674** & -0.0283 & -0.1263*** & -0.0156 & -0.1837*** & -0.1596 & -0.1535*** & -0.0714 \\
$\beta_\text{rest}^{\sigma}$ & - & 0.1955*** & - & 0.4363*** & - & 0.1593*** & - & 0.2502*** & - & 0.1856*** \\
$\beta_\text{shop}^{\mu}$ & -0.1224*** & -0.1744*** & -0.2051*** & -0.3363*** & -0.1592*** & -0.1612*** & -0.2046*** & -0.2775*** & -0.1206*** & -0.1907*** \\
$\beta_\text{shop}^{\sigma}$ & - & 0.0566*** & - & 0.0909 & - & 0.1259*** & - & 0.3556*** & - & 0.0973*** \\
$\beta_\text{sport}^{\mu}$ & -0.1464*** & -0.1479*** & -0.2871*** & -0.3261*** & -0.1841*** & -0.2025*** & -0.2402*** & -0.3621*** & -0.1836*** & -0.2118*** \\
$\beta_\text{sport}^{\sigma}$ & - & 0.1003*** & - & 0.4266*** & - & 0.3092*** & - & 0.6418*** & - & 0.2614*** \\
$\beta_\text{sm}^{\mu}$ & 0.1213*** & 0.1252* & 0.3174*** & 0.2888*** & 0.2236*** & 0.1563*** & 0.3764*** & 0.4199*** & 0.1776*** & 0.1434*** \\
$\beta_\text{sm}^{\sigma}$ & - & 0.6107*** & - & 0.4681*** & - & 0.4054*** & - & 1.0503*** & - & 0.6399*** \\
$\beta_\text{isGas}^{\mu}$ & -0.1778*** & -0.3027*** & -0.2509*** & -0.4688*** & -0.3397*** & -0.3126*** & -0.3396*** & -0.4871*** & -0.2418*** & -0.2344*** \\
$\beta_\text{isGas}^{\sigma}$ & - & 0.5651*** & - & 0.4866*** & - & 0.4613*** & - & 0.6096*** & - & 0.1547 \\
$\beta_\text{isWalkHome}^{\mu}$ & 1.0000*** & 1.0000*** & 1.0000*** & 1.0000*** & 1.0000*** & 1.0000*** & 1.0000*** & 1.0000*** & 1.0000*** & 1.0000*** \\
$\beta_\text{isWalkHome}^{\sigma}$ & - & 0.1364 & - & 0.1221 & - & 0.5924*** & - & 0.5810*** & - & 0.9382*** \\
\bottomrule
\end{tabular}
}
\caption[Parameter ratio values across models and folds for level 3 charging]{Parameter ratio values across models and folds for level 3 charging. *** indicates significance at 1\% level, ** significance at 5\% level, and * significance at 10\% level}
\label{ThirdArticle:TableParametersL3Combined}
\end{table}

\begin{table}[]
    \centering
        \resizebox{1.0\textwidth}{!}{
\begin{tabular}{lrrrrr}
\toprule
 & 0 & 1 & 2 & 3 & 4 \\
\midrule
Null log-likelihood & -8238.9696 & -8193.4340 & -8213.6351 & -8458.5849 & -8432.0811 \\
Final log-likelihood & -5002.5323 & -4898.4266 & -4931.0994 & -5005.9935 & -4950.2417 \\
$\rho$ & 0.3928 & 0.4022 & 0.3996 & 0.4082 & 0.4129 \\
$\bar{\rho}^2$ & 0.3912 & 0.4006 & 0.3981 & 0.4066 & 0.4114 \\
Akaike Information Criterion & 10031.0646 & 9822.8532 & 9888.1987 & 10037.9871 & 9926.4834 \\
Bayesian Information Criterion & 10110.3706 & 9902.1986 & 9967.6382 & 10117.6548 & 10006.1742 \\
\bottomrule
\end{tabular}
}
    \caption{Performance indicators for level 3 charging and the MNL model, estimation sets.}
    \label{ThirdArticle:TableEstimationL3MNL}
\end{table}

\begin{table}[]
    \centering
        \resizebox{1.0\textwidth}{!}{
\begin{tabular}{lrrrrr}
\toprule
 & 0 & 1 & 2 & 3 & 4 \\
\midrule
Null log-likelihood & -8238.9696 & -8193.4340 & -8213.6351 & -8458.5849 & -8432.0811 \\
Final log-likelihood & -3499.2768 & -3334.5446 & -3341.5290 & -3402.0914 & -3483.0145 \\
$\rho$ & 0.5753 & 0.5930 & 0.5932 & 0.5978 & 0.5869 \\
$\bar{\rho}^2$ & 0.5721 & 0.5898 & 0.5900 & 0.5947 & 0.5838 \\
Akaike Information Criterion & 7050.5537 & 6721.0892 & 6735.0580 & 6856.1827 & 7018.0291 \\
Bayesian Information Criterion & 7209.1657 & 6879.7801 & 6893.9369 & 7015.5183 & 7177.4106 \\
\bottomrule
\end{tabular}

}
    \caption{Performance indicators for level 3 charging and the MXL model, estimation sets.}
    \label{ThirdArticle:TableEstimationL3MXL}
\end{table}

\begin{table}[]
    \centering
    \resizebox{1.0\textwidth}{!}{
\begin{tabular}{lrrrrrrrrrr}
\toprule
Fold & \multicolumn{2}{c}{0} & \multicolumn{2}{c}{1} & \multicolumn{2}{c}{2} & \multicolumn{2}{c}{3} & \multicolumn{2}{c}{4} \\
Model & MNL & MXL & MNL & MXL & MNL & MXL & MNL & MXL & MNL & MXL \\
\midrule
Null log-likelihood & -540.6550 & -540.6550 & -539.4718 & -539.4718 & -521.6798 & -521.6798 & -527.4577 & -527.4577 & -534.0180 & -534.0180 \\
Final log-likelihood & -339.6777 & -341.4395 & -385.9410 & -383.0379 & -339.2753 & -341.6490 & -387.8558 & -376.3883 & -367.5855 & -350.2515 \\
$\rho$ & 0.3717 & 0.3685 & 0.2846 & 0.2900 & 0.3496 & 0.3451 & 0.2647 & 0.2864 & 0.3117 & 0.3441 \\
$\bar{\rho}^2$ & 0.3477 & 0.3204 & 0.2605 & 0.2418 & 0.3247 & 0.2953 & 0.2400 & 0.2371 & 0.2873 & 0.2954 \\
Akaike Information Criterion & 705.3554 & 734.8789 & 797.8819 & 818.0758 & 704.5506 & 735.2980 & 801.7116 & 804.7766 & 761.1709 & 752.5030 \\
Bayesian Information Criterion & 748.8678 & 821.9037 & 841.3323 & 904.9765 & 747.7497 & 821.6962 & 844.8471 & 891.0477 & 804.5589 & 839.2790 \\
DPSA null, min & 0.0556 & 0.0556 & 0.0556 & 0.0556 & 0.0556 & 0.0556 & 0.0556 & 0.0556 & 0.0556 & 0.0556 \\
DPSA null, mean & 0.0825 & 0.0825 & 0.0825 & 0.0825 & 0.0864 & 0.0864 & 0.0801 & 0.0801 & 0.0836 & 0.0836 \\
DPSA null, max & 0.2500 & 0.2500 & 0.2500 & 0.2500 & 0.2500 & 0.2500 & 0.2500 & 0.2500 & 0.2500 & 0.2500 \\
DPSA final, min & 0.0021 & 0.0049 & 0.0012 & 0.0035 & 0.0032 & 0.0079 & 0.0010 & 0.0039 & 0.0019 & 0.0043 \\
DPSA final, mean & 0.3268 & 0.3292 & 0.3084 & 0.2916 & 0.3392 & 0.3138 & 0.2827 & 0.2761 & 0.3119 & 0.3051 \\
DPSA final, max & 0.9729 & 0.8713 & 0.9509 & 0.8211 & 0.9710 & 0.8426 & 0.9606 & 0.8110 & 0.9784 & 0.8211 \\
DPSA average, min & 0.0076 & 0.0076 & 0.0085 & 0.0085 & 0.0009 & 0.0009 & 0.0009 & 0.0009 & 0.0006 & 0.0006 \\
DPSA average, mean & 0.0867 & 0.0867 & 0.0924 & 0.0924 & 0.0840 & 0.0840 & 0.0912 & 0.0912 & 0.0887 & 0.0887 \\
DPSA average, max & 0.2175 & 0.2175 & 0.2275 & 0.2275 & 0.2568 & 0.2568 & 0.2414 & 0.2414 & 0.2398 & 0.2398 \\
DPSA closest, min & 0.0000 & 0.0000 & 0.0000 & 0.0000 & 0.0000 & 0.0000 & 0.0000 & 0.0000 & 0.0000 & 0.0000 \\
DPSA closest, mean & 0.4143 & 0.4143 & 0.4163 & 0.4163 & 0.3415 & 0.3415 & 0.3039 & 0.3039 & 0.4183 & 0.4183 \\
DPSA closest, max & 1.0000 & 1.0000 & 1.0000 & 1.0000 & 1.0000 & 1.0000 & 1.0000 & 1.0000 & 1.0000 & 1.0000 \\
\bottomrule
\end{tabular}
}
    \caption{Performance indicators for level 3 charging, validation sets.}
    \label{ThirdArticle:TableValidationL3Combined}
\end{table}

\begin{table}[]
    \centering
\begin{tabular}{lrr}
\toprule
 & MNL & MXL \\
\midrule
$\beta_\text{distNear}^{\mu}$ & 0.7380 & 0.1428 \\
$\beta_\text{distNear}^{\sigma}$ & - & 0.2937 \\
$\beta_\text{distFar}^{\mu}$ & -0.1313 & -0.1775 \\
$\beta_\text{distFar}^{\sigma}$ & - & 0.0804 \\
$\beta_\text{ff}^{\mu}$ & -0.0064 & -0.0600 \\
$\beta_\text{ff}^{\sigma}$ & - & 0.5904 \\
$\beta_\text{leis}^{\mu}$ & 0.2261 & 0.2498 \\
$\beta_\text{leis}^{\sigma}$ & - & 0.6051 \\
$\beta_\text{mall}^{\mu}$ & 0.4562 & 0.3487 \\
$\beta_\text{mall}^{\sigma}$ & - & 1.0184 \\
$\beta_\text{outletsNear}^{\mu}$ & -0.0304 & 0.2117 \\
$\beta_\text{outletsNear}^{\sigma}$ & - & 0.6232 \\
$\beta_\text{outletsFar}^{\mu}$ & 0.1123 & -0.0007 \\
$\beta_\text{outletsFar}^{\sigma}$ & - & 0.2953 \\
$\beta_\text{rest}^{\mu}$ & -0.1211 & -0.0736 \\
$\beta_\text{rest}^{\sigma}$ & - & 0.1711 \\
$\beta_\text{shop}^{\mu}$ & -0.1624 & -0.2280 \\
$\beta_\text{shop}^{\sigma}$ & - & 0.1453 \\
$\beta_\text{sport}^{\mu}$ & -0.2083 & -0.2501 \\
$\beta_\text{sport}^{\sigma}$ & - & 0.3479 \\
$\beta_\text{sm}^{\mu}$ & 0.2433 & 0.2267 \\
$\beta_\text{sm}^{\sigma}$ & - & 0.6349 \\
$\beta_\text{isGas}^{\mu}$ & -0.2700 & -0.3611 \\
$\beta_\text{isGas}^{\sigma}$ & - & 0.4555 \\
$\beta_\text{isWalkHome}^{\mu}$ & 1.0000 & 1.0000 \\
$\beta_\text{isWalkHome}^{\sigma}$ & - & 0.4740 \\
\bottomrule
\end{tabular}
    \caption[Average parameter ratios for level 3 charging]{Parameter ratios for MNL and MXL models for level 3 charging, average across folds.}
    \label{ThirdArticle:TableL3RatioTest}
\end{table}

\section{Discussion}
\label{ThirdArticle:SectionDiscussion}


We begin by discussing aspects which apply to both the level 2 and level 3 models. Firstly, we note that both the MNL and MXL models have a DPSA (notably, a fitting factor) that is better than the null model (which reflects random assignment). Combined with the $\rho$ and $\bar{\rho}^2$, this indicates an increased predictive power in comparison with random assignment. 

Secondly, consistent with the literature, we observe a high level of heterogeneity among the population. Notably, the standard deviation terms in the MXL models are nearly always highly statistically significant, and larger than their mean counterpart.  In addition, we note that there are some differences in sign or significance between the two models, such as the network distance when near home for level 2 charging. For the average parameter ratios for level 2 in Table~\ref{ThirdArticle:TableL2RatioTest}, we also note that the parameter value for the network distance in the MXL model is more than double that of the MNL model. 
In terms of performance on the validation sets, we note that the MNL model consistently has a better $\rho$, $\bar{\rho}^2$, AIC, and BIC than the MXL model, while the better DPSA (including, notably, the fitting factor) depends on the fold and charging level. We note that the MXL model nevertheless displays strong performance in an optimization application (see Section~\ref{ThirdArticle:SectionDiscussion:opt}).

Thirdly, when considering charging close to home, the effects of the network distance and the charging station being a short walk from home are comparable (since, by assumption, the distance is bound by 1.5 km in this case), while the number of outlets has a lesser (though still positive and statistically significant) effect. However, when considering charging far from home, the most important attribute is certainly the distance. More specifically, there is a higher bound on the attribute values for the distance compared with the other parameters, with a maximum of around 50 km for level 3 and 100 km for level 2. By contrast, the density measures for the amenities have a maximum of around 5, while the number of outlets has a maximum of  5 outlets for level 3 and 13 outlets for level 2. Fourthly, we observe different preferences between level 2 and level 3 charging stations, particularly as it pertains to the amenities. This highlights the different use cases for the types of charging, consistent with \cite{Anderson2018}.

As it pertains to the impact of amenities for level 2 charging in Table~\ref{ThirdArticle:TableL2RatioTest}, the only overall positive and statistically significant effects were for the density of restaurants. In addition to this, the effects of the density of leisure facilities, shopping malls, sports facilities, and supermarkets were overall neutral, with a standard deviation term at least an order of magnitude higher than the mean terms. 
Despite the positive impact of restaurants, the density of fast food locations had a negative and statistically significant effect, suggesting a preference for longer duration or possible other confounding factors. Likewise, the charging station being located at a gas station had a negative and statistically significant effect, as did the densities of shopping facilities. As the main value of gas stations was associated with habit in \cite{Philipsen2015, Philipsen2016}, this suggests that the long-term habits of EV users have adapted around the new refueling requirements. Given their respective bounds, the effects of the amenities near the charging station can have a larger impact than the number of outlets. 

For level 3 charging, we note that the average parameter ratios in Table~\ref{ThirdArticle:TableL3RatioTest} indicates a positive effect for the network distance when charging close to home. As discussed previously, a possible explanation for the positive effect of the distance is that users with access to home charging may prefer to use that rather than a public charger when near home. However, in over 75\% of the observations, there are no charging stations within 1.5 km of the user~(for more details, see Table~\ref{ThirdArticle:TableThresholdSplitL3} in \ref{ThirdArticle:AppendixAttributes}). As such, unexpected phenomena may be due to a lack of observations. For the amenities, the density of leisure locations, shopping malls, and supermarkets have a positive and statistically significant effect, restaurant density exhibits a neutral and highly variable effect,  while shopping facilities, sports facilities, and gas stations had negative and statistically significant effects. The aversion to sports facilities for level 3 charging may be attributable to longer duration of these activities in comparison with the charging time. However, while the preference for leisure locations, shopping malls, and supermarkets are consistent with the discussion in Section~\ref{ThirdArticle:SectionAttributeEnconding}, the aversion to shopping facilities is somewhat surprising. Since the ``shopping'' category in OpenStreetMaps considers a diverse portfolio of shop types (including specialty stores), a further refinement in terms of types of shopping locations may be necessary for improving the accuracy.

\section{Application in Optimisation}\label{ThirdArticle:SectionDiscussion:opt}

In this section, we examine if the differences between the MNL and MXL models lead to different solution recommendations in an optimisation context. For this purpose, we propose two different optimisation models based on improving the service level for level 2 stations of the charging network for a set of customers. Model~\eqref{ModelPMedian} is formulated as a p-median problem, where the goal is to maximise the total utility from assigning each customer to their preferred station. Specifically, we have 
\begin{subequations}
\begin{alignat}{3}
\operatorname{Maximise} \, & \sum_{i \in N} \sum_{j \in M} u_{ij} y_{ij}, 
\label{ModelPMedian:Objective} 
\\
\text{subject to }
& \sum_{j \in M} y_{ij} = 1, \quad &&   i \in N, 
\label{ModelPMedian:Assignment} 
\\
& y_{ij} \leq z_j, \quad && i \in N, j \in M, 
\label{ModelPMedian:Assignment_y_z} 
\\
& \sum_{j \in M} z_j = p,
\label{ModelPMedian:OpenP} 
\\
& y_{ij} \geq 0, \quad && i \in N, j \in M,
\label{ModelPMedian:domain_y} 
\\
& z_j \in \{0,1\}, \quad && j \in M,
\label{ModelPMedian:domain_z} 
\end{alignat}
\label{ModelPMedian}
\end{subequations}
where $z_j$ indicates whether station $j$ is open or closed, and $y_{ij}$ specifies whether customer $i$ is assigned or not to station $j$. Objective function~\eqref{ModelPMedian:Objective} maximises the sum of the utility across all customers achieved from assigning them to their preferred station. Constraints~\eqref{ModelPMedian:Assignment} enforce that each customer may only be assigned to one station, while Constraints~\eqref{ModelPMedian:Assignment_y_z} impose that customers may only be assigned to open stations. Finally, Constraint~\eqref{ModelPMedian:OpenP} indicates that exactly $p$ stations may be opened, and the remaining are domain constraints.

By contrast, Model~\eqref{ModelMaxMin} maximises the utility of the customer with the worst service level (in other words, it maximises the minimum utility across customers). More concretely, we have
\begin{subequations}
\begin{alignat}{3}
\operatorname{Maximise} \ \ \, & \alpha, 
\label{ModelMaxMin:Objective} 
\\
\text{subject to } \ \ 
& \alpha \leq \sum_{j \in M} u_{ij} y_{ij}, \quad &&   i \in N, 
\label{ModelMaxMin:ObjectiveBound} 
\\
& \alpha \in \mathbb{R}, \quad && i \in N, 
\notag
\\
& \text{Constraints}~\eqref{ModelPMedian:Assignment}-~\eqref{ModelPMedian:domain_z}. \notag
\end{alignat}
\label{ModelMaxMin}
\end{subequations}
Objective function~\eqref{ModelMaxMin:Objective} maximises our auxiliary variable $\alpha$, which is upper-bounded by the utility of each customer through Constraints~\eqref{ModelMaxMin:ObjectiveBound}. The remaining decision variables and constraints are as in Model~\eqref{ModelPMedian}.

For this experiment, we use the 2021 Statistics Canada~\citep{StatsCan2021} census of Montreal to generate the set of customers $N$, where each customer home is assumed to be located at the centroid of the dissemination areas (the smallest geographical unit available). Distances to each station are determined via the driving network in Open Street Maps~\citep{OpenStreetMap} using these centroids. Meanwhile, 200 new locations for charging stations were randomly generated at locations within the city, ensuring they are at least 200 meters from other charging stations (both existing and new). 

For each customer $i$ and charging station $j$, we generate four different sets of utilities for comparison: 
\begin{description}
    \item[\distance:] The utility is given simply by the negative distance, without considering the amenities. In this case, the resulting model reduces to a classic p-median formulation. 
    \item[\mnl:] The utility is calculated following the average parameter ratios for the MNL model in Table~\ref{ThirdArticle:TableL2RatioTest}. 
    \item[\mxlmean:] The utility is calculated following the average parameter ratios for the MXL model in Table~\ref{ThirdArticle:TableL2RatioTest}. 
    \item[MXL-25:] This set of utilities is generated via the same distribution of utility values as the mean case discussed above. However, instead of setting the utility value of station $j$ to customer $i$ as the mean of their distribution, we instead set the value as the 25th percentile.
\end{description}
We emphasise that in each case, a single utility value $u_{ij}$ is generated for each customer $i$ and each station $j$. As a consequence, the instance sizes are identical across cases.  See~\ref{ThirdArticle:AppendixSimulUtilities} for details on the simulation of utilities for the models using the multinomial and mixed logit.

We consider a series of experiments, with different combinations of various hyperparameters: (i)~the p-median Model~\eqref{ModelPMedian} or the max-min Model~\eqref{ModelMaxMin}, (ii)~the utility values calculated via \distance, \mnl, \mxlmean, or \mxltwentyfive, (iii)~whether or not to consider existing charging stations as available alternatives for the demand points, and (iv)~the number of stations to open, with $p \in \{ 10, 25, 50, 75, 100\}$. 

A valid comparison method across the different utility specifications is not trivial, as the utility values themselves are not comparable. As such, we compare the resulting optimal solutions obtained from solving the models to optimality. 

In Table~\ref{TableExperiments}, we present descriptive statistics across the set of experiments, examining the solutions themselves. The `Active' column indicates the percentage of the $p$ stations that have users assigned to them. For example, a value of $0.75$ indicates that three-quarters of the stations have at least one user assigned to them. In other words, one quarter of the stations are not used at all, and could be substituted with other stations without modifying the objective function value. Therefore, lower values indicate higher levels of symmetry within the problem. 

The columns for `Similarity' examines how often the same set of solutions are selected by multiple different utility specification. More precisely, for each set of hyperparameter values, we examine the overlap in the set of stations which are opened in the optimal solution with each utility specification. For example, a value of $0.75$ indicates that, on average, three-quarters of the stations are identical across the utility specifications. We examine the similarity considering either all of the utility specifications or comparing {\mxlmean} and {\mxltwentyfive} versions. 

\begin{figure}[p]
    \centering
        \begin{subfigure}[b]{0.475\textwidth}
            \centering
            \includegraphics[width=\textwidth]{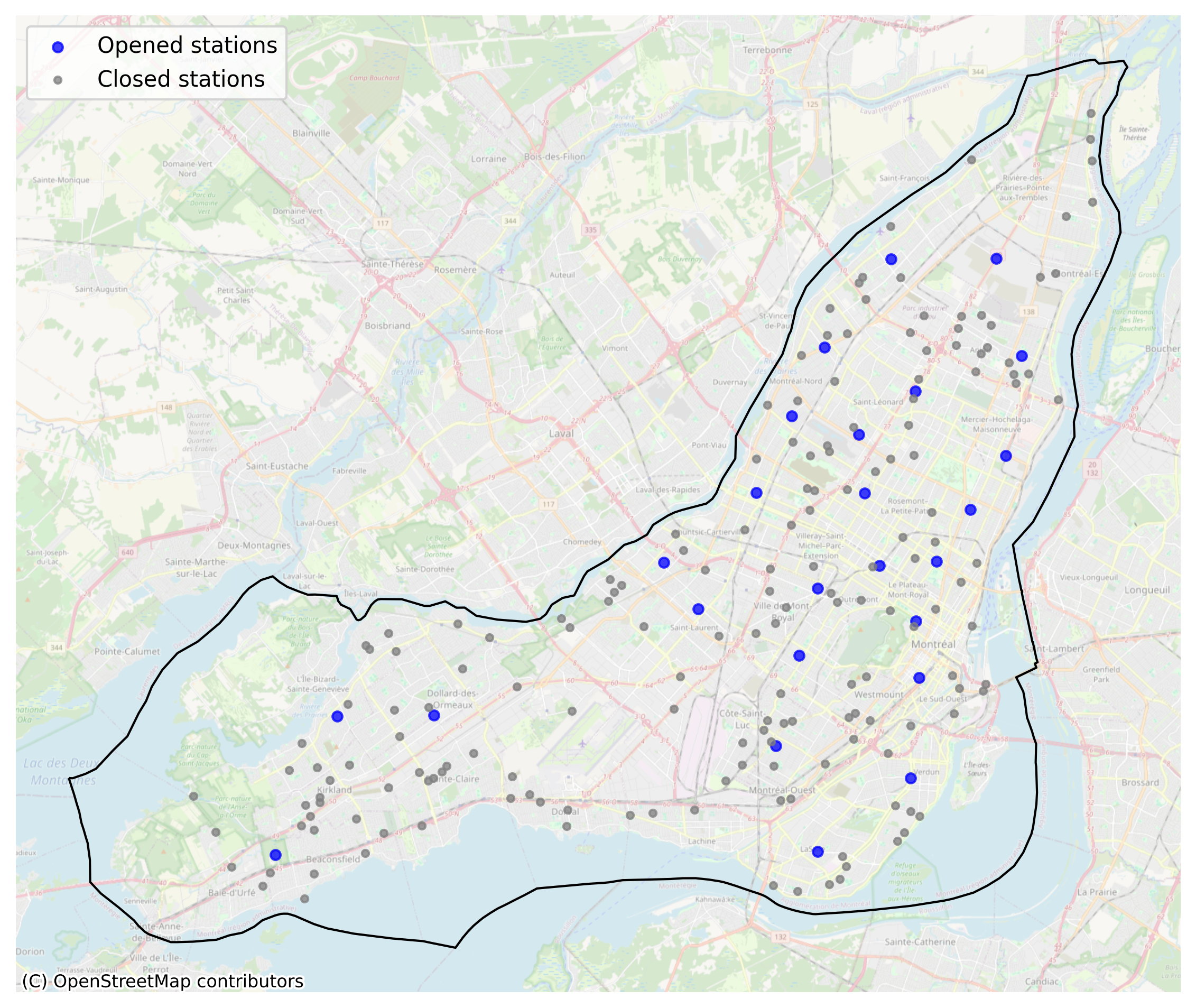}
            \caption
            {{\small \distance}}    
            \label{FigureDistanceSolution}
        \end{subfigure}
        \hfill
        \begin{subfigure}[b]{0.475\textwidth}
            \centering
            \includegraphics[width=\textwidth]{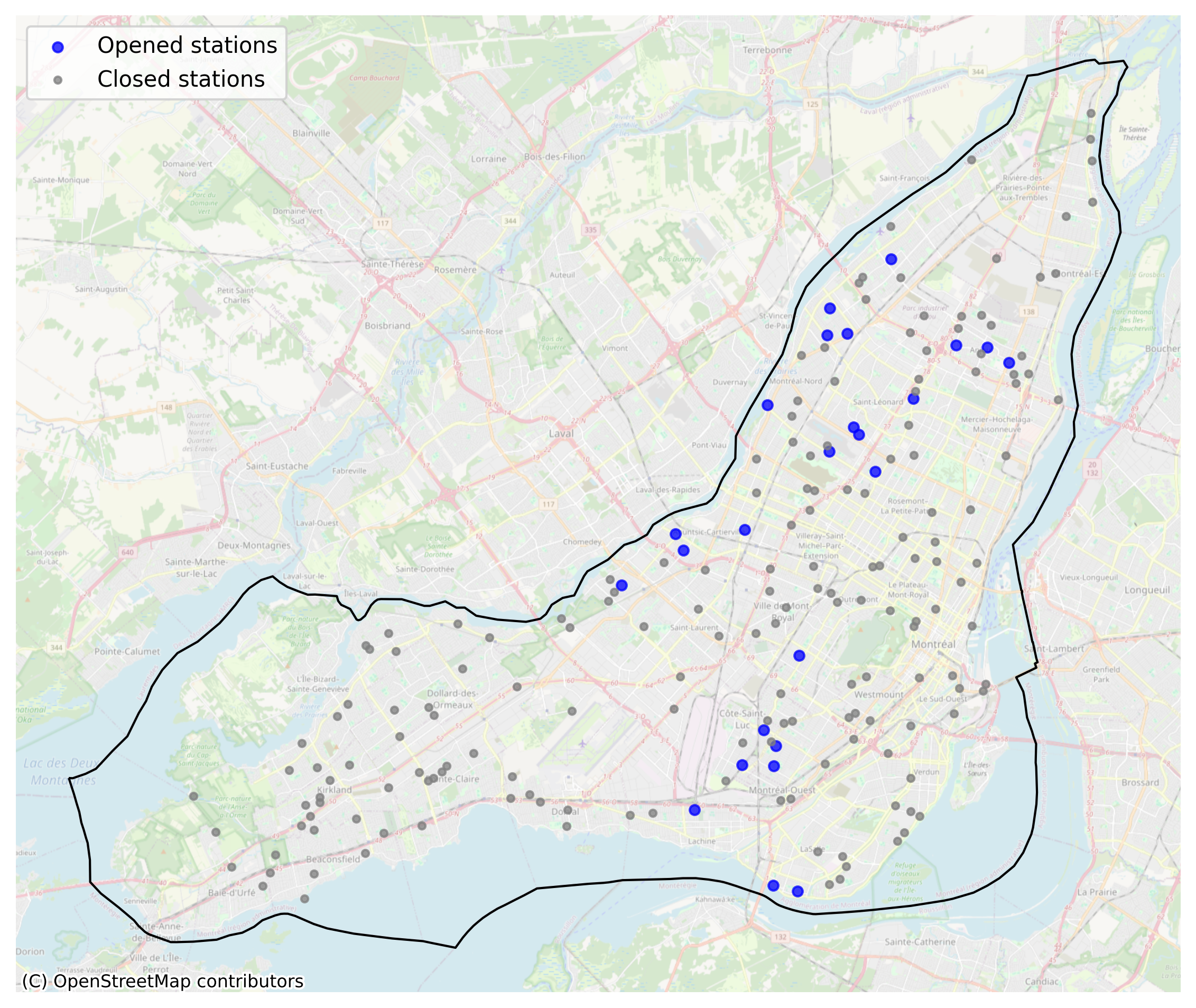}
            \caption
            {{\small \mnl}}    
            \label{FigureMNLSolution}
        \end{subfigure}
        \vskip\baselineskip
        \begin{subfigure}[b]{0.475\textwidth}
            \centering
            \includegraphics[width=\textwidth]{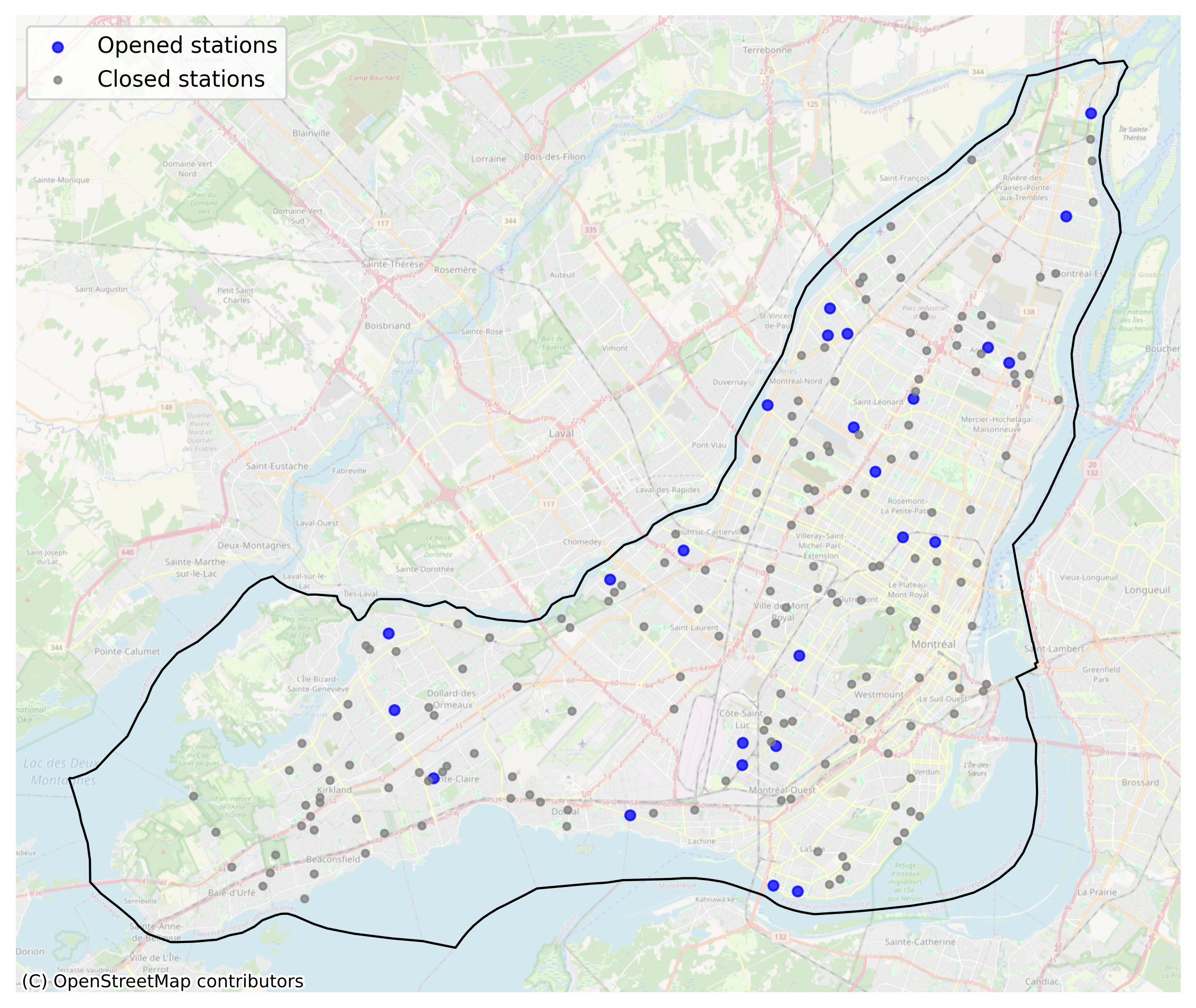}
            \caption
            {{\small \mxlmean}}    
            \label{FigureDistanceSolution}
        \end{subfigure}
        \hfill
        \begin{subfigure}[b]{0.475\textwidth}
            \centering
            \includegraphics[width=\textwidth]{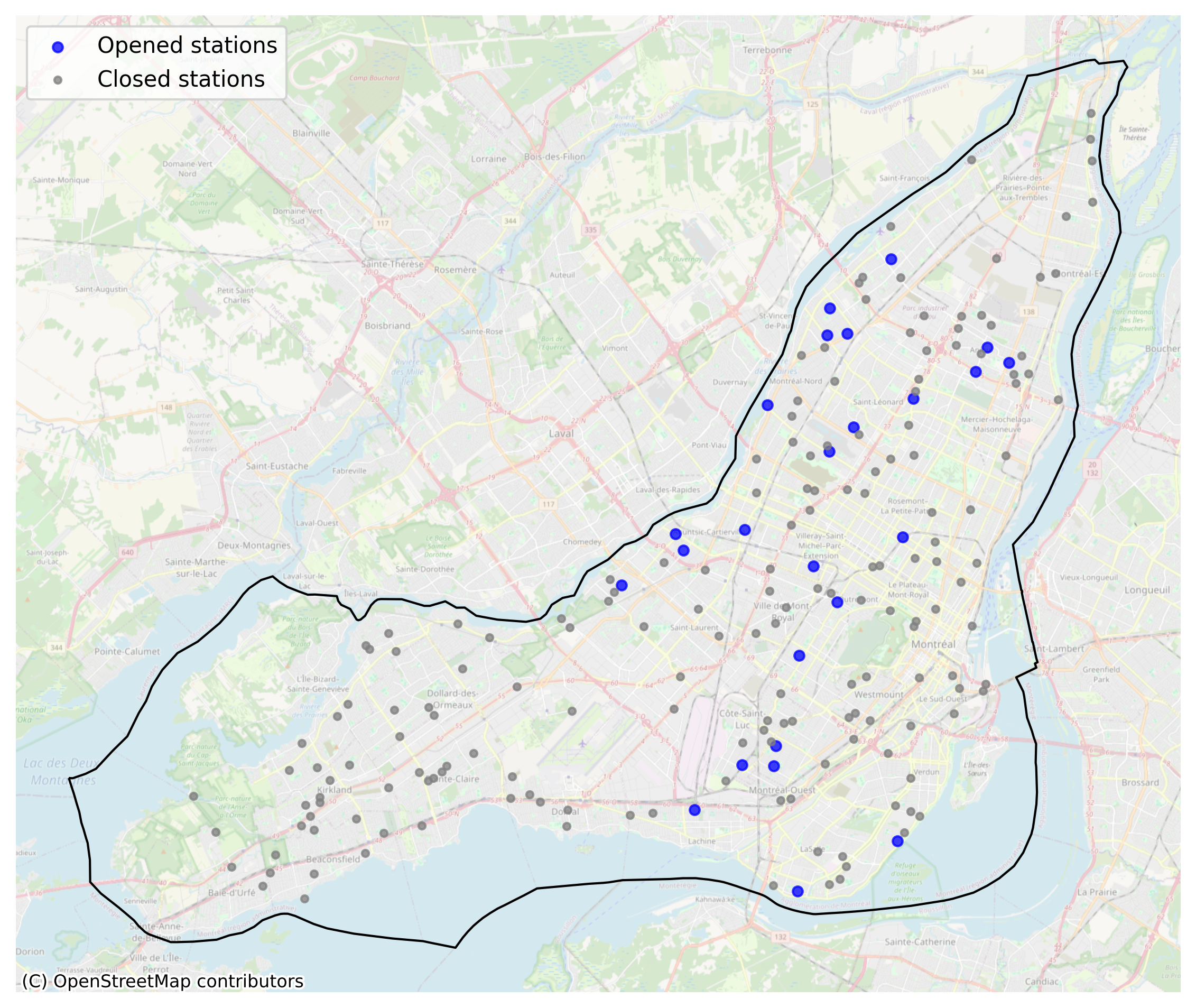}
            \caption
            {{\small \mxltwentyfive}}    
            \label{FigureMNLSolution}
        \end{subfigure}
        \caption{\small Solution (opened and closed stations) for each type of utility.} 
        \label{FigureSolutions}
\end{figure}

On the one hand, the similarity for the MXL specifications under the p-median model is significant, reaching at least approximately 72\%. On the other hand, overall--and especially for the max-min model--the solution similarity is low. This suggests that the utility specifications must differ sufficiently to lead to distinct charging station network decisions.
To illustrate one example, we show in Figure~\ref{FigureSolutions} the stations which were selected to be opened and closed in the optimal solutions of each utility specification. As can be seen, notably with respect to opened stations, there are distinct differences. 

 Additionally, we compare the objective function values obtained under different utility specifications to assess whether these solutions, despite their differences, yield equivalent-quality solutions. In Table~\ref{TableObjectiveComparison}, we examine an example using the p-median model with $p=25$ and considering existing charging stations. 
The rows indicate which set of utilities were used to find the solution. The `Objective' column indicates the objective value while the remaining columns are the optimality gaps when applied to the set of utilities. In other words, entry $(i,j)$ indicates the optimality gap of the solution obtained using utilities $i$ when applied to utilities $j$. We observe that when we evaluate the solutions of {\distance} and {\mnl} in {\mxltwentyfive}, we obtain the largest optimality gaps. In contrast, for the remaining cases, the gaps are below 4.44\% with the solutions of {\mxlmean} and {\mxltwentyfive} being the ones providing high-quality solutions across all the different utility sets. Hence, there seems to be advantages in adopting the utilities of the mixed logit. The complete table across all experiments is presented in~\ref{app:Addtional_opt_results}.

\begin{table}[]
    \centering
    \resizebox{1.0\textwidth}{!}{
\begin{tabular}{lrrr|rrr}
\toprule
 & \multicolumn{3}{c}{P-Median} & \multicolumn{3}{c}{Max-Min} \\
 & Active & Similarity, all & Similarity, MXL & Active & Similarity, all & Similarity, MXL \\
\midrule
mean & 1.0000 & 0.5399 & 0.7834 & 0.5128 & 0.4479 & 0.5744 \\
std & 0.0000 & 0.1035 & 0.0871 & 0.4159 & 0.1262 & 0.1109 \\
min & 1.0000 & 0.3704 & 0.6250 & 0.0000 & 0.2632 & 0.5000 \\
25\% & 1.0000 & 0.4926 & 0.7232 & 0.0367 & 0.3687 & 0.5242 \\
50\% & 1.0000 & 0.5360 & 0.8043 & 0.8100 & 0.4609 & 0.5319 \\
75\% & 1.0000 & 0.6157 & 0.8333 & 0.8600 & 0.4710 & 0.5319 \\
max & 1.0000 & 0.6897 & 0.8929 & 1.0000 & 0.7407 & 0.7937 \\
\bottomrule
\end{tabular}

}
    \caption{Descriptive statistics over the set of experiments.}
    \label{TableExperiments}
\end{table}

\begin{table}[]
    \centering
\begin{tabular}{lrrrrr}
\toprule
 & Objective & \multicolumn{4}{|c}{Optimality Gap} \\
Utilities & \multicolumn{1}{c}{Value} &  \multicolumn{1}{|r}{Distance} & MNL & MXLMean & MXL25 \\
\midrule
Distance & -349.63	& 0.00	& 1.52	&0.52	&17.81\\
MNL 	 & 3806.87	& 4.44	& 0.00	&0.17	&15.70\\
MXLMean  & 4236.95	& 3.15	& 0.45	&0.00	&0.96\\
MXL25 	 & 797.29	& 3.49	& 0.20	&0.21	&0.00\\
\bottomrule
\end{tabular}
    \caption{Optimality gap (in \%) of the objective function for the p-median model with $p=25$ and considering existing charging station. 
    }
    \label{TableObjectiveComparison}
\end{table}

\section{Conclusion}
\label{ThirdArticle:SectionConclusion}

In this work, we estimated multinomial logit and mixed logit models for analysing and predicting electric vehicle charging station usage, for both level 2 and level 3 locations. Our models rely on characteristics of charging stations that are readily available for charging network operators, and thus make them well-suited for optimisation purposes. Internal 5-fold cross validation was performed to evaluate the accuracy and predictive power of our models, demonstrating improved performance over random selection.

By design, these estimation results can be integrated into electric vehicle charging stations network design models such as \cite{Luo2015, Cui2018, Lamontagne2023a, shen2024sequential}. In addition, our findings indicate that the demand models integrated into these optimisation programs should factor panel effects. Notably, users who recharge frequently on the public network significantly influence the average preferences of the MNL model, while the MXL model captures taste heterogeneity. As a consequence, optimisation models which do not include panel effects may integrate unrealistic user behaviour, thus incorrectly evaluating the potential benefits of stations. We provided an application example, which illustrated that the placement of charging stations could be impacted by the model selection. Although we observed significant differences in the charging stations network design solutions depending on the utility model used, solutions generated with the MXL model consistently demonstrated high quality, even when evaluated with using other choice model specifications.   

Overall, users prefer level 2 charging stations closer to their home (especially within walking distance), while preferences for level 3 charging stations were more variable. Far from home, the most significant attribute was the distance between the home area of the users and the charging locations, with the number of outlets having a positive (if more muted) impact. High heterogeneity was observed for nearly all attributes, suggesting the existence of users with higher preference for all attributes (even those for which a general aversion is observed). 

For the amenities, we observe that the preferences varied depending on the charging level. For level 2 charging, users exhibited a preference for charging near restaurants, while avoiding fast food locations and gas stations. For level 3 charging, users favoured charging stations near leisure locations, shopping malls, and supermarkets, while they had an aversion to those with higher densities of shopping facilities and sports facilities, and to those near gas stations. For both charging levels, the densities of the other attributes had an overall neutral effect (which varied by fold), and accompanied by standard deviation terms much higher than the mean terms. This suggests very high levels of preference heterogeneity across the population, with subsets of users who either heavily favour or are heavily averse to these characteristics.

Our work focuses on the intracity case, where both the charging stations and the users are located within the same city. However, future work could extend this approach to the intercity case, where the charging stations are being used during the course of long-distance trips. Additionally, while not relevant in our case study, future work could also incorporate charging price to the selection process.

\section*{Acknowledgements}

The authors gratefully acknowledge the assistance of Jean-Luc Dupre from \emph{Direction Mobilit\'e} of \emph{Hydro-Qu\'ebec} for sharing his expertise on EV charging stations and the network.  We also gratefully acknowledge Miguel F. Anjos, Ismail Sevim, and Nagisa Sugishita for their insights into the project.

\section*{Funding Sources}
This research was supported by Hydro-Québec, NSERC Collaborative Research and Development Grant CRDPJ 536757 - 19, and the FRQ-IVADO Research Chair in Data Science for Combinatorial Game Theory.



\bibliography{References/References_copy, References/Data}

\clearpage
\appendix 

\section{Data Description}
\label{ThirdArticle:AppendixDataShared}

Though this dataset has been used in \cite{Lamontagne2023a, Elhattab2023, Parent2023}, descriptive statistics were not provided for charging behaviour. These are provided  along with a comparison with other reports on charging behaviour (mentioned previously in the literature review).

\begin{figure}
    \centering
\includegraphics[width = 0.75\textwidth]{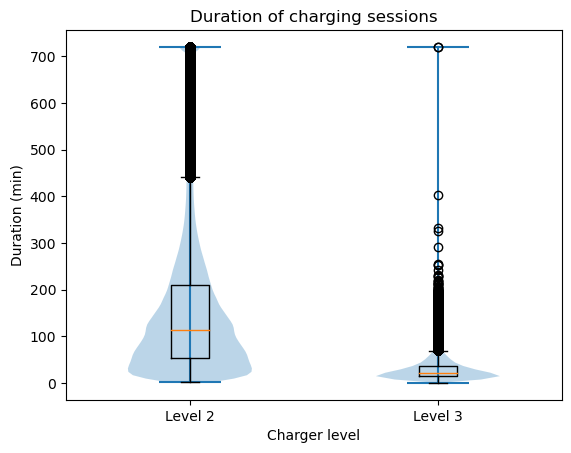}
    \caption{Distribution of duration of charging, by level of charging outlet.}
    \label{ThirdArticle:Data:SessionDuration}
\end{figure}

\begin{figure}
    \centering
\includegraphics[width = 0.75\textwidth]{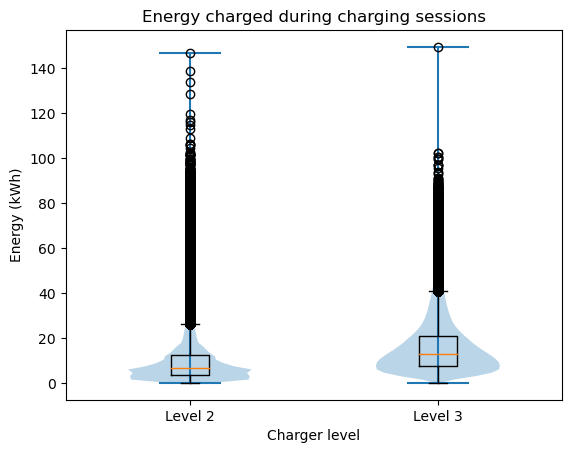}
    \caption{Distribution of energy from charging, by level of charger.}
    \label{ThirdArticle:Data:SessionEnergy}
\end{figure}

\begin{figure}
    \centering
\includegraphics[width = 0.75\textwidth]{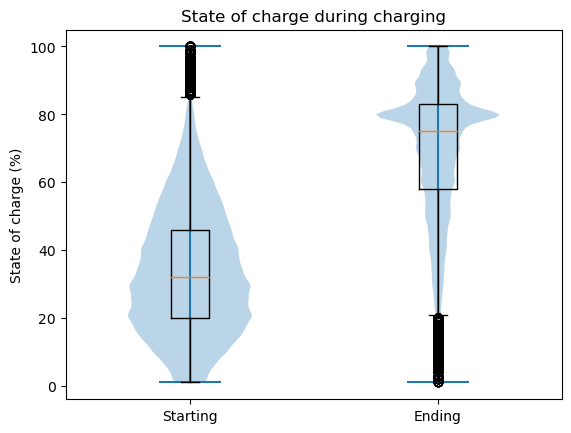}
    \caption{Distribution of starting and ending state of charge at level 3 chargers.}
    \label{ThirdArticle:Data:SessionSoC}
\end{figure}

We present summary statistics of charging session characteristics in Figures~\ref{ThirdArticle:Data:SessionDuration}-\ref{ThirdArticle:Data:SessionSoC}. 
These statistics are separated, as applicable, by level of the charging outlet. However, these statistics are similar among different account types, and thus account types are aggregated. These values are compared with reports of level 2 charging sessions on public charging infrastructure in \cite{vanDenHoed2013, Morrissey2016, Helmus2020} and level 3 sessions on public charging in \cite{Morrissey2016}. While charging sessions information is presented in \cite{Tal2020}, home charging is not separated from public charging, rendering the comparability questionable.  However, as level 3 charging is not available at home, the relevant charging session information in \cite{Tal2020} can be used. 
\begin{itemize}
    \item In Figure~\ref{ThirdArticle:Data:SessionDuration}, we report the distribution of the duration of charging sessions. In this graph, the horizontal lines indicate the 25th, 50th, and 75th percentile, while the width of the shaded background illustrates the distribution of each duration (with larger width corresponding to more sessions).   As vehicles may be connected to chargers much longer than their charging time~\citep[e.g.][]{Morrissey2016, Helmus2020}, we cap the duration at 720 minutes, which applies for around 2.2\% of level 2 sessions and only two level 3 sessions. The median of 116 minutes for level 2 and 23 minutes for level 3 is comparable to the corresponding values of 128.78 minutes and 26.62 minutes for level 2 and level 3 presented in \cite{Morrissey2016}, and within the distribution of around 17 minutes to 37 minutes (depending on vehicle type and day) reported in \cite{Tal2020}. However, this does differ significantly from the average level 2 connection time of 435 minutes presented in \cite{vanDenHoed2013} and median of 537 minutes in \cite{Helmus2020}. In \cite{Helmus2020}, they note that over 75\% of EV users are dependent on public charging facilities. As such, a higher ratio of overnight charging may contribute to this increased duration. 
    \item In Figure~\ref{ThirdArticle:Data:SessionEnergy}, we report the distribution of energy (in kWh) charged during the session. The median energy per level 2 session is 6.9 kWh while the median for level 3 session is 13.4 kWh. The value for level 2 charging matches with those in the literature, with 6.8 kWh per session in \cite{Morrissey2016} and \cite{Helmus2020}, and 8.31 kWh in \cite{vanDenHoed2013}.  The value for level 3 charging is higher than the 8.32 kWh reported in \cite{Morrissey2016}, and within the distribution of 7.9 kWh to 26.3 kWh (depending on vehicle type) reported in \cite{Tal2020}.
    \item In Figure~\ref{ThirdArticle:Data:SessionSoC}, we report the distribution for the starting and ending state of charge for level 3 charging sessions, with median values of 30\% and 73\%, respectively. We recall that state of charge information is not available for level 2 sessions. These values are slightly below the distribution of 30.3 to 48.4\% for starting and 74.7 to 92.5\% for ending state of charge (depending on vehicle type) presented in \cite{Tal2020}. State of charge information is not reported in \cite{Morrissey2016}. 
\end{itemize}

\begin{figure}
    \centering
\includegraphics[width = 0.75\textwidth]{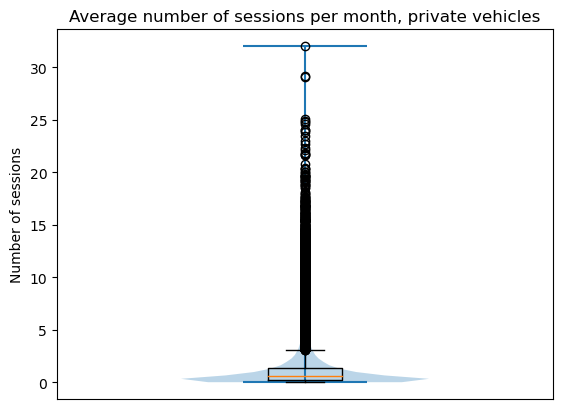}
    \caption{Distribution of average number of sessions per month, private vehicles.}
    \label{ThirdArticle:Data:PrivateSessions}
\end{figure}

\begin{figure}
    \centering
\includegraphics[width = 0.75\textwidth]{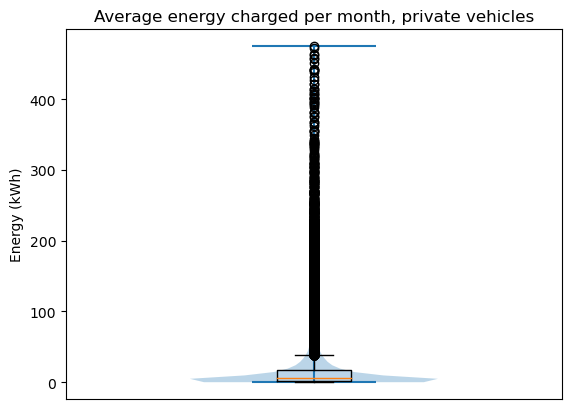}
    \caption{Distribution of average energy charged per month, private vehicles.}
    \label{ThirdArticle:Data:PrivateEnergy}
\end{figure}

\begin{figure}
    \centering
\includegraphics[width = 0.75\textwidth]{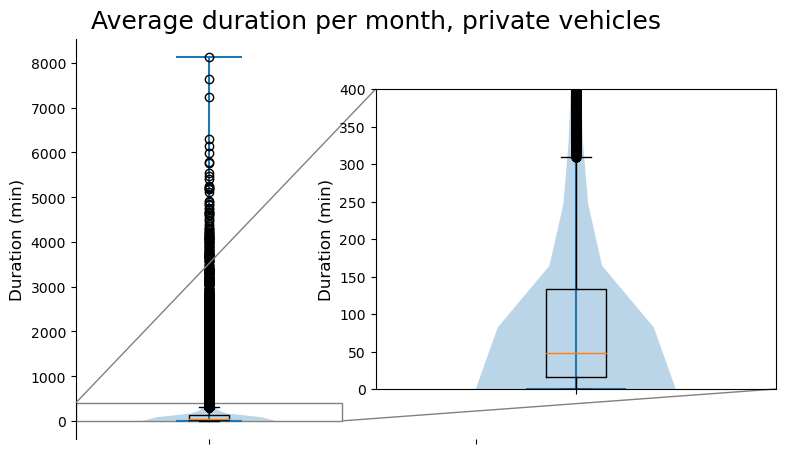}
    \caption{Distribution of average monthly time spent charging, private vehicles.}
    \label{ThirdArticle:Data:PrivateDuration}
\end{figure}

\begin{figure}
    \centering
\includegraphics[width = 0.75\textwidth]{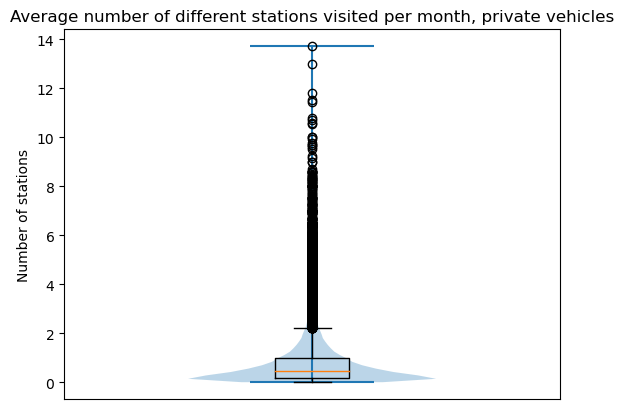}
    \caption{Distribution of average number of different stations per month, private vehicles.}
    \label{ThirdArticle:Data:PrivateStations}
\end{figure}

Rather than reporting on charging sessions, in Figures~\ref{ThirdArticle:Data:PrivateSessions}-\ref{ThirdArticle:Data:PrivateStations}, we present summary statistics for private vehicle accounts. Summary statistics for shared accounts are presented and discussed in \ref{ThirdArticle:AppendixDataShared}, while they are not presented for rental vehicle or unplugged accounts. For rental vehicle accounts this is due to the small sample size, while there is no charging activity to report for  unplugged accounts. 
\begin{itemize}
    \item By construction, the average number of sessions and total energy charged are lower for private vehicles (Figures~\ref{ThirdArticle:Data:PrivateSessions} and \ref{ThirdArticle:Data:PrivateEnergy}) than for shared accounts. Median values for the number of charging sessions and total energy per month for private vehicles are, respectively, 0.57 sessions and 5.89 kWh, while the median values for shared accounts are 2.75 sessions and 34.7 kWh. Comparable values in the literature were discussed at the beginning of this section as part of our user classification process.
    \item In Figures~\ref{ThirdArticle:Data:PrivateDuration}, we present the distribution for the average amount of time spent charging every month, with the duration of each charging session capped at 720 minutes as before. The median value is 47.7 minutes for private vehicles. In terms of comparable values, \cite{Tal2020} reports an average \emph{daily} duration of between 97.33 and 266.96 minutes (depending on vehicle type). However, we note that this charging duration includes home charging and, as such, may not be an accurate proxy for public charging, even for those that lack access to home chargers. For public charging specifically, by combining the average number of 18 to 20 charging sessions per year and 7.25 hour session duration in \cite{vanDenHoed2013}, we obtain a monthly average of between 652.5 and 720 minutes. As with the duration of individual sessions, this higher value for members may be attributed to a higher rate of overnight charging.
    \item  In Figure~\ref{ThirdArticle:Data:PrivateStations}, we present the distribution for the average number of different stations visited each month. The median value is 0.45 stations for private vehicles. Only \cite{vanDenHoed2013} reported the number of stations visited, with an average of 4 and 77 different locations visited per year for private vehicles and car sharing vehicles, respectively. Accordingly, the value for private vehicles is quite comparable to that of \cite{vanDenHoed2013}. 
\end{itemize}

\begin{figure}
    \centering
\includegraphics[width = 0.6\textwidth]{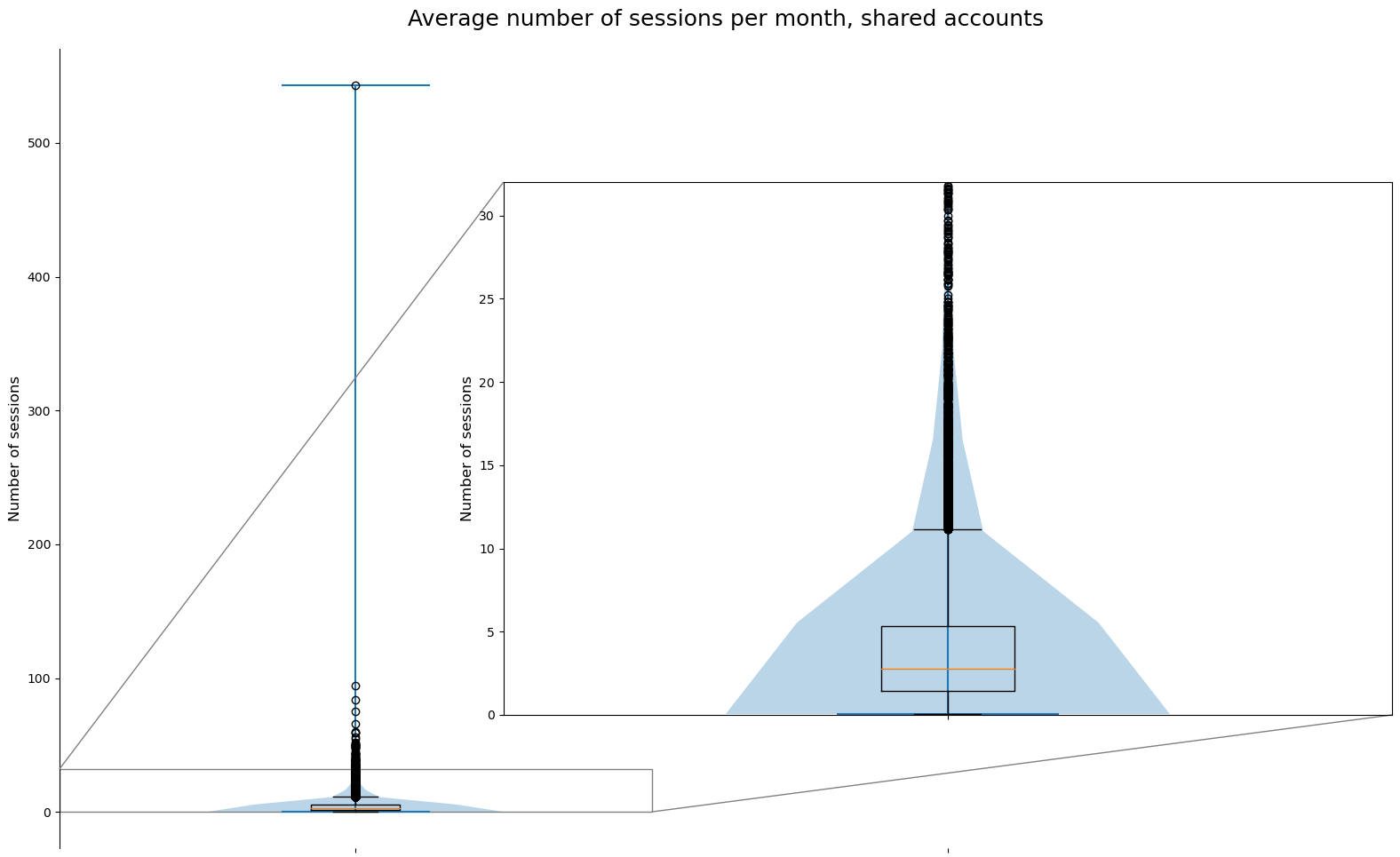}
    \caption{Distribution of average number of sessions per month, shared vehicles.}
    \label{ThirdArticle:Data:SharedSessions}
\end{figure}

\begin{figure}
    \centering
\includegraphics[width = 0.6\textwidth]{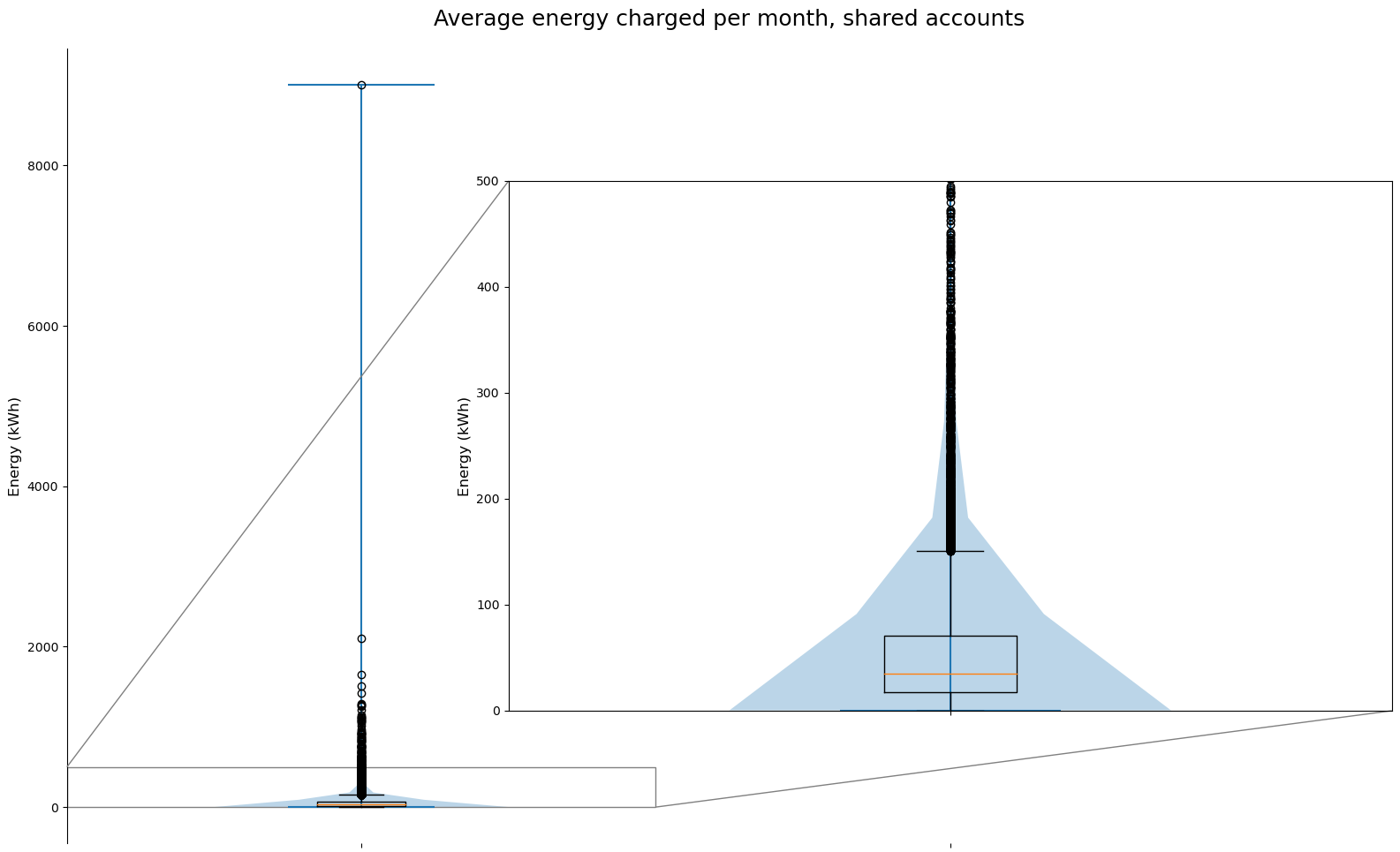}
    \caption{Distribution of average energy charged per month, shared vehicles.}
    \label{ThirdArticle:Data:SharedEnergy}
\end{figure}

\begin{figure}
    \centering
\includegraphics[width = 0.6\textwidth]{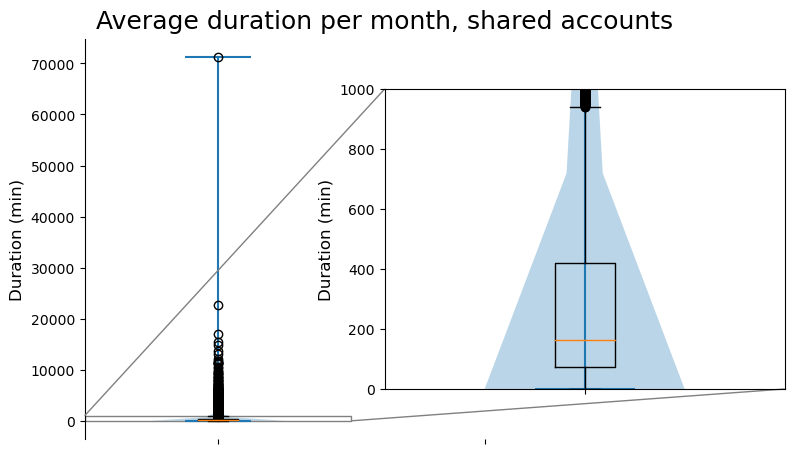}
    \caption{Distribution of average monthly time spent charging, shared vehicles.}
    \label{ThirdArticle:Data:SharedDuration}
\end{figure}

\begin{figure}
    \centering
\includegraphics[width = 0.6\textwidth]{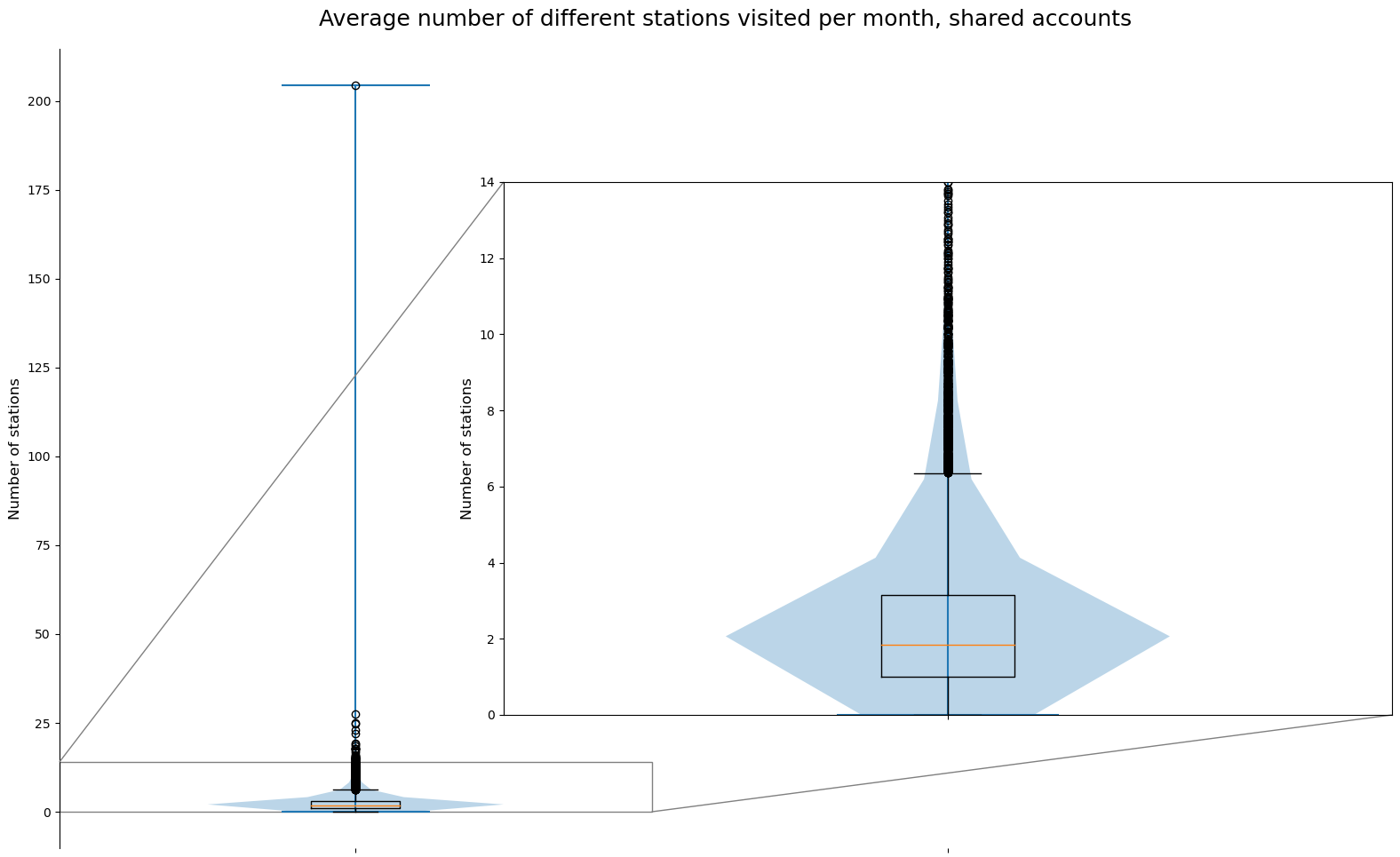}
    \caption{Distribution of average number of different stations per month, shared vehicles.}
    \label{ThirdArticle:Data:SharedStations}
\end{figure}

In general, shared members have a much higher usage of public charging infrastructure across all metrics. This can be seen by comparing the figures for private vehicles with their counterpart in Figures~\ref{ThirdArticle:Data:SharedSessions}-\ref{ThirdArticle:Data:SharedStations}. For ease of reading, the median values for private vehicles are repeated in this section.
\begin{itemize}
    \item By construction, the average number of sessions and total energy charged are presented in  Figures~\ref{ThirdArticle:Data:SharedSessions} and \ref{ThirdArticle:Data:SharedEnergy}. Median values for the number of charging sessions and total energy per month for private vehicles are, respectively, 0.57 sessions and 5.89 kWh, while the median values for shared accounts are 2.75 sessions and 34.7 kWh.
    \item In Figure~\ref{ThirdArticle:Data:SharedDuration}, we present the distribution for the average amount of time spent charging every month, with the duration of each charging session capped at 720 minutes as before. The median values are 47.7 minutes for private vehicles and 162.7 minutes for shared accounts. 
    \item  In Figure~\ref{ThirdArticle:Data:SharedStations}, we present the distribution for the average number of different stations visited each month. The median values are 0.45 and 1.85 stations, respectively, for private vehicles and shared accounts. 
\end{itemize}

\clearpage

\section{Distribution of Attributes}
\label{ThirdArticle:AppendixAttributes}
\begin{figure}[h]
    \centering
    \includegraphics[width = 0.4\textwidth]{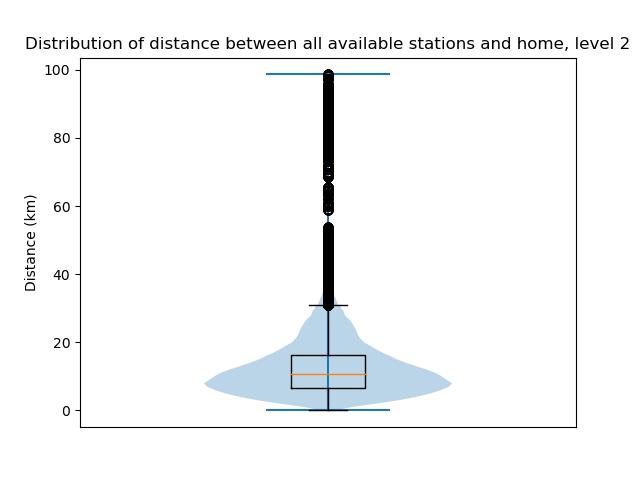}
    \includegraphics[width = 0.4\textwidth]{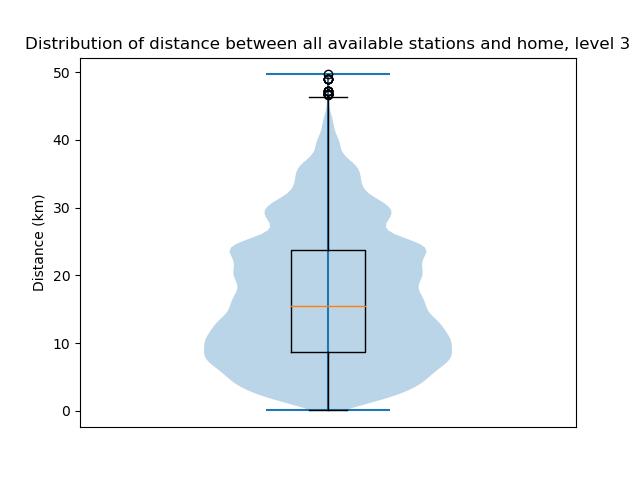}
    \caption{Distance from the user's home to the all available charging stations.}
    \label{ThirdArticle:FigureDistributionDistanceToHomeAll}
\end{figure}

\begin{figure}
    \centering
    \includegraphics[width = 0.4\textwidth]{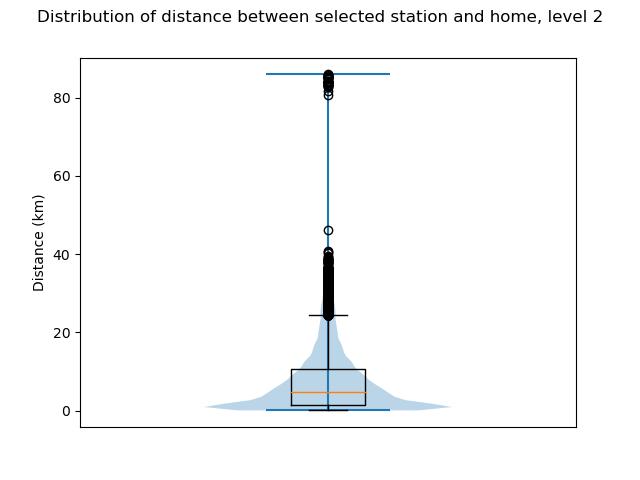}
    \includegraphics[width = 0.4\textwidth]{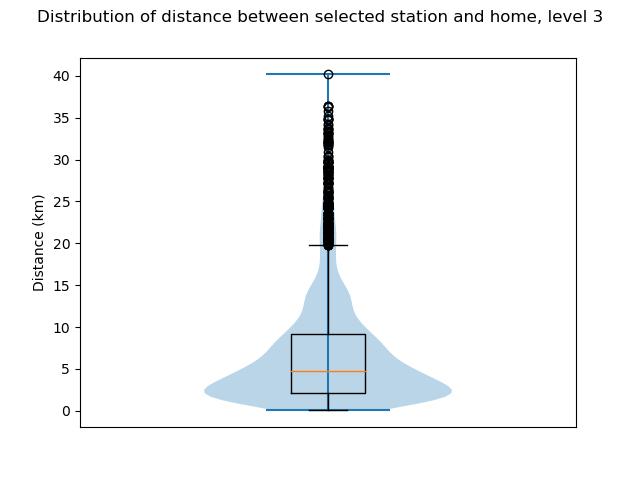}
    \caption{Distance from the user's home to the selected charging station, before removing outliers.}
    \label{ThirdArticle:FigureDistributionDistanceToHomeSelected}
\end{figure}

\begin{figure}
    \centering
    \includegraphics[width = 0.4\textwidth]{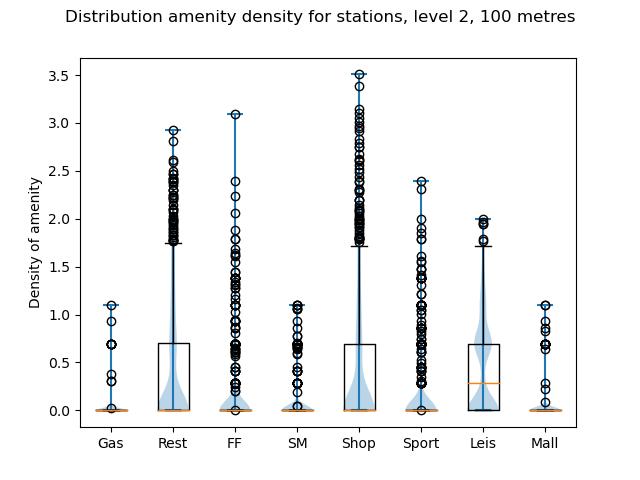}
    \includegraphics[width = 0.4\textwidth]{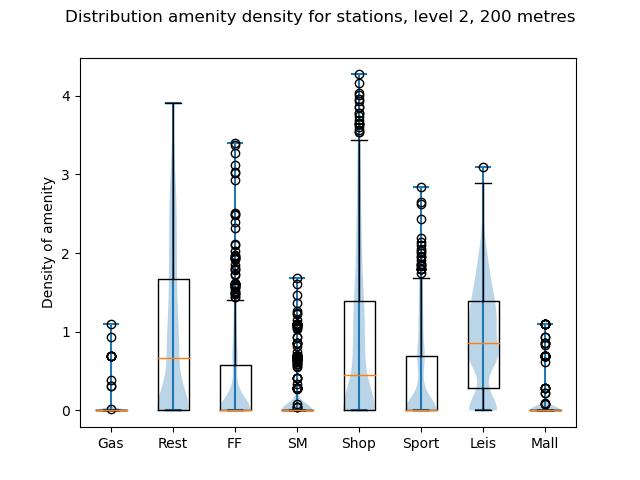}
    \includegraphics[width = 0.4\textwidth]{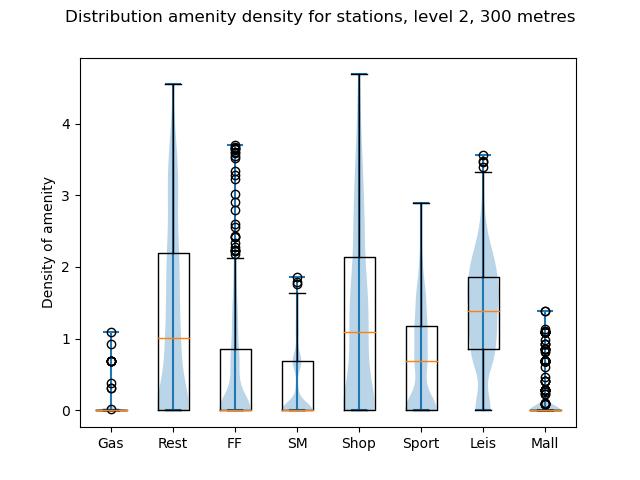}
    \includegraphics[width = 0.4\textwidth]{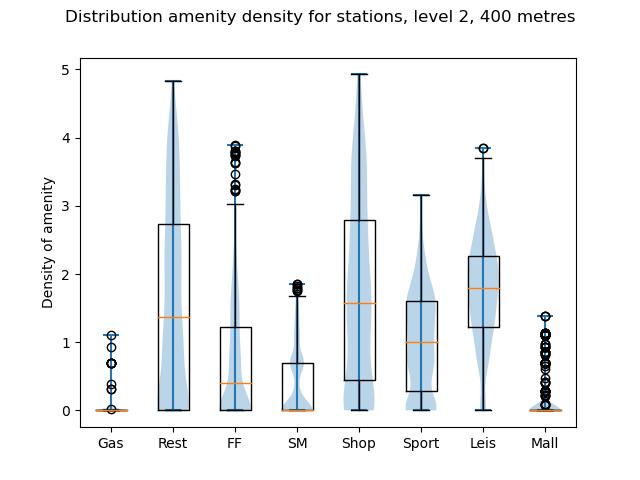}
    \caption{Density of amenities at different distance thresholds, level 2.}
    \label{ThirdArticle:FigureDensityLevel2}
\end{figure}

\begin{figure}
    \centering
    \includegraphics[width = 0.4\textwidth]{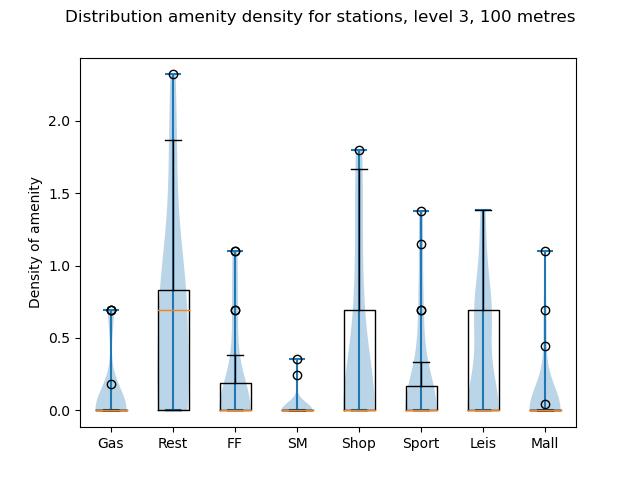}
    \includegraphics[width = 0.4\textwidth]{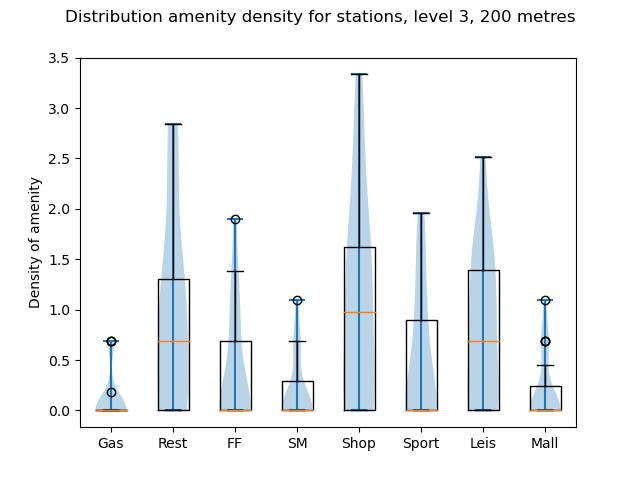}
    \includegraphics[width = 0.4\textwidth]{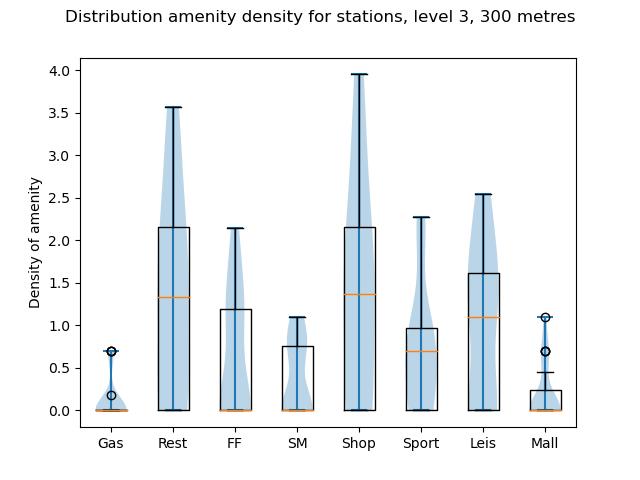}
    \includegraphics[width = 0.4\textwidth]{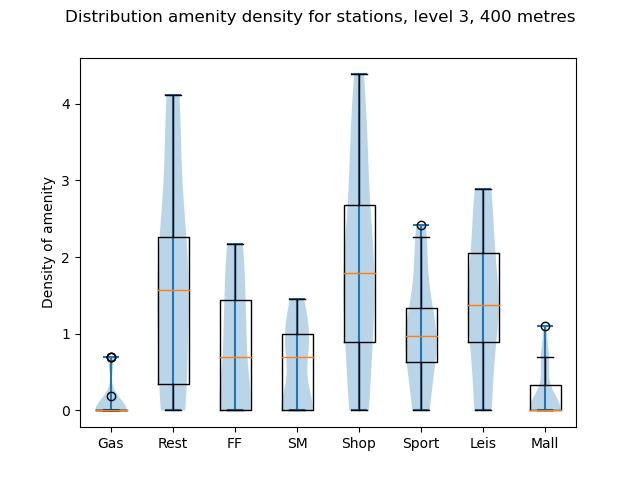}
    \caption{Density of amenities at different distance thresholds, level 3.}
    \label{ThirdArticle:FigureDensityLevel3}
\end{figure}

 See Figures~\ref{ThirdArticle:FigureDistributionDistanceToHomeAll} to~\ref{ThirdArticle:FigureDensityLevel3}, and Tables~\ref{ThirdArticle:TableThresholdSplitL2} and~\ref{ThirdArticle:TableThresholdSplitL3}.

\begin{table}[]
    \centering
\begin{tabular}{lrrrrr}
\toprule
 & 500m & 750m & 1000m & 1250m & 1500m \\
\midrule
mean & 0.8127 & 2.0735 & 3.9996 & 6.5259 & 9.5881 \\
std & 1.1767 & 2.4713 & 4.5060 & 6.9330 & 9.7409 \\
min & 0 & 0 & 0 & 0 & 0 \\
25\% & 0 & 0 & 1 & 1 & 2 \\
50\% & 0 & 1 & 3 & 5 & 7 \\
75\% & 1 & 3 & 6 & 10 & 15 \\
max & 15 & 23 & 42 & 55 & 70 \\
\bottomrule
\end{tabular}

    \caption{Distribution of the number of charging stations with the utility function thresholds for level 2 charging.}
    \label{ThirdArticle:TableThresholdSplitL2}
\end{table}

\begin{table}[]
    \centering
\begin{tabular}{lrrrrr}
\toprule
 & 500m & 750m & 1000m & 1250m & 1500m \\
\midrule
mean & 0.0152 & 0.0658 & 0.1080 & 0.1609 & 0.2223 \\
std & 0.1225 & 0.2480 & 0.3156 & 0.3843 & 0.4442 \\
min & 0 & 0 & 0 & 0 & 0 \\
25\% & 0 & 0 & 0 & 0 & 0 \\
50\% & 0 & 0 & 0 & 0 & 0 \\
75\% & 0 & 0 & 0 & 0 & 0 \\
max & 1 & 1 & 2 & 2 & 2 \\
\bottomrule
\end{tabular}
    \caption{Distribution of the number of charging stations with the utility function thresholds for level 3 charging.}
    \label{ThirdArticle:TableThresholdSplitL3}
\end{table}

\clearpage

\section{Additional Estimation Results}
\label{ThirdArticle:AppendixAdditionalResults}
See Tables~\ref{ThirdArticle:Table:ParametersDistancesL2} to~\ref{ThirdArticle:Table:ScoresDistances200L3}.

\begin{table}
    \centering
    \resizebox{0.6\textwidth}{!}{
\begin{tabular}{lllllll}
\toprule
 &  & 0 & 1 & 2 & 3 & 4 \\
Parameter & Distance &  &  &  &  &  \\
\midrule
\multirow[t]{4}{*}{$\beta_{distFar}^{\mu}$} & 100 & -0.1591*** & -0.1422*** & -0.1283*** & -0.1619*** & -0.1582*** \\
 & 200 & -0.1616*** & -0.145*** & -0.1293*** & -0.1608*** & -0.1608*** \\
 & 300 & -0.163*** & -0.1454*** & -0.1296*** & -0.164*** & -0.1608*** \\
 & 400 & -0.1623*** & -0.1443*** & -0.129*** & -0.1628*** & -0.1596*** \\
\cline{1-7}
\multirow[t]{4}{*}{$\beta_{distNear}^{\mu}$} & 100 & 0.5991*** & 1.0374*** & 0.7813*** & 0.723*** & 0.8911*** \\
 & 200 & 0.6434*** & 1.0774*** & 0.8139*** & 0.9184*** & 0.9184*** \\
 & 300 & 0.6391*** & 1.0422*** & 0.7845*** & 0.6879*** & 0.8845*** \\
 & 400 & 0.5709*** & 0.9572*** & 0.6848*** & 0.6034*** & 0.811*** \\
\cline{1-7}
\multirow[t]{4}{*}{$\beta_{ff}^{\mu}$} & 100 & 0.0005 & 0.0791** & 0.0129 & -0.0465 & -0.0472 \\
 & 200 & -0.1015*** & -0.1728*** & -0.2049*** & -0.1663*** & -0.1663*** \\
 & 300 & -0.0961*** & -0.1161*** & -0.1098*** & -0.1311*** & -0.0985*** \\
 & 400 & -0.1042*** & -0.0907*** & -0.1004*** & -0.1099*** & -0.0513** \\
\cline{1-7}
\multirow[t]{4}{*}{$\beta_{isGas}^{\mu}$} & 100 & -0.324*** & -0.415*** & -0.003 & -0.4791*** & -0.4317*** \\
 & 200 & -0.2123** & -0.3062*** & 0.1074 & -0.3652*** & -0.3652*** \\
 & 300 & -0.2414** & -0.3383*** & 0.1231 & -0.4756*** & -0.4155*** \\
 & 400 & -0.2582*** & -0.3583*** & 0.08 & -0.483*** & -0.4052*** \\
\cline{1-7}
\multirow[t]{4}{*}{$\beta_{isWalkHome}^{\mu}$} & 100 & 3.0644*** & 3.1369*** & 3.4862*** & 3.3313*** & 3.2841*** \\
 & 200 & 3.1021*** & 3.1773*** & 3.5186*** & 3.3121*** & 3.3121*** \\
 & 300 & 3.0992*** & 3.1498*** & 3.4934*** & 3.3087*** & 3.2845*** \\
 & 400 & 3.0426*** & 3.0794*** & 3.407*** & 3.2376*** & 3.2243*** \\
\cline{1-7}
\multirow[t]{4}{*}{$\beta_{leis}^{\mu}$} & 100 & -0.1515*** & 0.0269 & -0.2136*** & 0.2344*** & -0.0001 \\
 & 200 & -0.0884*** & -0.0846*** & -0.1754*** & -0.1332*** & -0.1332*** \\
 & 300 & -0.159*** & -0.1166*** & -0.2607*** & -0.1748*** & -0.1916*** \\
 & 400 & -0.2475*** & -0.1591*** & -0.3302*** & -0.2523*** & -0.2983*** \\
\cline{1-7}
\multirow[t]{4}{*}{$\beta_{mall}^{\mu}$} & 100 & 0.6322*** & 0.6795*** & 0.5981*** & 1.1677*** & 0.7985*** \\
 & 200 & 0.3663*** & 0.3058*** & 0.2705*** & 0.4025*** & 0.4025*** \\
 & 300 & 0.2711*** & 0.2975*** & 0.197*** & 0.6687*** & 0.4962*** \\
 & 400 & 0.1605*** & 0.2353*** & 0.2327*** & 0.5829*** & 0.3509*** \\
\cline{1-7}
\multirow[t]{4}{*}{$\beta_{outletsFar}^{\mu}$} & 100 & 0.1374*** & 0.1284*** & 0.1163*** & 0.1102*** & 0.0954*** \\
 & 200 & 0.1221*** & 0.1146*** & 0.1058*** & 0.0839*** & 0.0839*** \\
 & 300 & 0.1136*** & 0.1067*** & 0.0958*** & 0.0923*** & 0.0724*** \\
 & 400 & 0.1169*** & 0.1068*** & 0.0927*** & 0.078*** & 0.0718*** \\
\cline{1-7}
\multirow[t]{4}{*}{$\beta_{outletsNear}^{\mu}$} & 100 & 0.3219*** & 0.2669*** & 0.2393*** & 0.2475*** & 0.1863*** \\
 & 200 & 0.3258*** & 0.2711*** & 0.2422*** & 0.19*** & 0.19*** \\
 & 300 & 0.3251*** & 0.2681*** & 0.2399*** & 0.2442*** & 0.1857*** \\
 & 400 & 0.3181*** & 0.2602*** & 0.2314*** & 0.2347*** & 0.1766*** \\
\cline{1-7}
\bottomrule
\end{tabular}
}
\caption[Parameter values for level 2 charging with the MNL model, and various amenity thresholds]{Parameter values for level 2 charging with the MNL model, and various amenity thresholds. The utility function threshold is set at 1.5km and the sample size is 5,000. *** indicates significance at 1\% level, ** significance at 5\% level, and * significance at 10\% level}
    \label{ThirdArticle:Table:ParametersDistancesL2}
\end{table}

\begin{table}[]
    \centering
    \resizebox{0.6\textwidth}{!}{
\begin{tabular}{lllllll}
\toprule
 &  & 0 & 1 & 2 & 3 & 4 \\
Parameter & Distance &  &  &  &  &  \\
\midrule
\multirow[t]{4}{*}{$\beta_{rest}^{\mu}$} & 100 & 0.2273*** & 0.3801*** & 0.377*** & 0.2648*** & 0.3425*** \\
 & 200 & 0.3181*** & 0.4165*** & 0.3732*** & 0.3822*** & 0.3822*** \\
 & 300 & 0.3556*** & 0.3912*** & 0.376*** & 0.3777*** & 0.3956*** \\
 & 400 & 0.3469*** & 0.3746*** & 0.3966*** & 0.3989*** & 0.3719*** \\
\cline{1-7}
\multirow[t]{4}{*}{$\beta_{shop}^{\mu}$} & 100 & -0.0223 & -0.0575** & -0.0291 & 0.0453* & -0.0146 \\
 & 200 & -0.0941*** & -0.0504*** & -0.0095 & -0.0692*** & -0.0692*** \\
 & 300 & -0.098*** & -0.0981*** & -0.0496** & -0.1129*** & -0.1494*** \\
 & 400 & -0.0785*** & -0.1419*** & -0.0949*** & -0.1603*** & -0.1428*** \\
\cline{1-7}
\multirow[t]{4}{*}{$\beta_{sm}^{\mu}$} & 100 & 0.1547*** & 0.0994 & -0.1608*** & 0.0734 & 0.0238 \\
 & 200 & 0.1724*** & 0.1505*** & 0.0684 & 0.3255*** & 0.3255*** \\
 & 300 & 0.1006*** & 0.1129*** & 0.0437 & 0.1702*** & 0.2572*** \\
 & 400 & 0.0646* & 0.1595*** & 0.049 & 0.1336*** & 0.2197*** \\
\cline{1-7}
\multirow[t]{4}{*}{$\beta_{sport}^{\mu}$} & 100 & 0.1797*** & 0.0869** & -0.0112 & -0.0255 & -0.1036** \\
 & 200 & 0.0814*** & 0.1227*** & -0.0167 & 0.0039 & 0.0039 \\
 & 300 & 0.0798*** & 0.0694*** & 0.0765*** & 0.1168*** & 0.0413 \\
 & 400 & 0.0895*** & 0.0161 & 0.0566** & 0.1298*** & 0.0859*** \\
\cline{1-7}
\bottomrule
\end{tabular}
}
\caption[Parameter values for level 2 charging with the MNL model, and various amenity thresholds (continued)]{Parameter values for level 2 charging with the MNL model, and various amenity thresholds (continued). The utility function threshold is set at 1.5km and the sample size is 5,000. *** indicates significance at 1\% level, ** significance at 5\% level, and * significance at 10\% level}
    \label{ThirdArticle:Table:ParametersDistancesL2Continued}
\end{table}

\begin{table}[]
    \centering
    \resizebox{\textwidth}{!}{
\begin{tabular}{llrrrrrrrr}
\toprule
 & Distance & \multicolumn{2}{c}{100} & \multicolumn{2}{c}{200} & \multicolumn{2}{c}{300} & \multicolumn{2}{c}{400} \\
 &  & Mean & Std & Mean & Std & Mean & Std & Mean & Std \\
Test set & Score &  &  &  &  &  &  &  &  \\
\midrule
\multirow[t]{6}{*}{Estimation} & Akaike Information Criterion & 52019.7700 & 652.8600 & 51937.1900 & 604.2400 & 51859.3200 & 668.2300 & 51865.2400 & 664.9100 \\
 & Bayesian Information Criterion & 52104.4900 & 652.8600 & 52021.9100 & 604.2400 & 51944.0400 & 668.2300 & 51949.9600 & 664.9100 \\
 & Final log-likelihood & -25996.8853 & 326.4277 & -25955.5947 & 302.1192 & -25916.6605 & 334.1137 & -25919.6182 & 332.4531 \\
 & Null log-likelihood & -30901.9274 & 50.3773 & -30900.1624 & 50.3833 & -30901.9274 & 50.3773 & -30901.9274 & 50.3773 \\
 & $\rho$ & 0.1587 & 0.0104 & 0.1600 & 0.0095 & 0.1613 & 0.0107 & 0.1612 & 0.0106 \\
 & $\bar{\rho}^2$ & 0.1583 & 0.0104 & 0.1596 & 0.0095 & 0.1609 & 0.0107 & 0.1608 & 0.0106 \\
\cline{1-10}
\multirow[t]{12}{*}{Validation} & Akaike Information Criterion & 53369.9400 & 193.0700 & 53115.9900 & 164.1200 & 53083.3200 & 157.7800 & 53144.9600 & 162.9300 \\
 & Bayesian Information Criterion & 53453.9200 & 193.0700 & 53199.9700 & 164.1200 & 53167.3000 & 157.7800 & 53228.9400 & 162.9300 \\
 & DPSA final, max & 0.5792 & 0.0222 & 0.6033 & 0.0277 & 0.5946 & 0.0226 & 0.5816 & 0.0290 \\
 & DPSA final, mean & 0.0212 & 0.0009 & 0.0213 & 0.0010 & 0.0213 & 0.0010 & 0.0214 & 0.0010 \\
 & DPSA final, min & 0.0000 & 0.0000 & 0.0000 & 0.0000 & 0.0000 & 0.0000 & 0.0000 & 0.0000 \\
 & DPSA null, max & 0.0038 & 0.0000 & 0.0038 & 0.0000 & 0.0038 & 0.0000 & 0.0038 & 0.0000 \\
 & DPSA null, mean & 0.0017 & 0.0000 & 0.0017 & 0.0000 & 0.0017 & 0.0000 & 0.0017 & 0.0000 \\
 & DPSA null, min & 0.0014 & 0.0000 & 0.0014 & 0.0000 & 0.0014 & 0.0000 & 0.0014 & 0.0000 \\
 & Final log-likelihood & -26671.9687 & 96.5372 & -26544.9946 & 82.0617 & -26528.6589 & 78.8882 & -26559.4805 & 81.4667 \\
 & Null log-likelihood & -30199.1322 & 9.0394 & -30198.0326 & 7.5853 & -30199.1322 & 9.0394 & -30199.1322 & 9.0394 \\
 & $\rho$ & 0.1168 & 0.0032 & 0.1210 & 0.0028 & 0.1215 & 0.0026 & 0.1205 & 0.0027 \\
 & $\bar{\rho}^2$ & 0.1164 & 0.0032 & 0.1205 & 0.0028 & 0.1211 & 0.0026 & 0.1201 & 0.0027 \\
\cline{1-10}
\bottomrule
\end{tabular}
}
    \caption[Performance indicators for level 2 charging with the MNL model, and various amenity thresholds]{Performance indicators for level 2 charging with the MNL model, and various amenity thresholds. The utility function threshold is set at 1.5km and the sample size is 5,000.}
    \label{ThirdArticle:Table:ScoresDistancesL2}
\end{table}

\begin{table}[]
    \centering
    \resizebox{0.6\textwidth}{!}{
\begin{tabular}{lllllll}
\toprule
 &  & 0 & 1 & 2 & 3 & 4 \\
Parameter & Distance &  &  &  &  &  \\
\midrule
\multirow[t]{5}{*}{$\beta_{distFar}^{\mu}$} & 0.5 & -0.1952*** & -0.1781*** & -0.1632*** & -0.1967*** & -0.2023*** \\
 & 0.75 & -0.1783*** & -0.168*** & -0.1517*** & -0.1869*** & -0.187*** \\
 & 1.0 & -0.17*** & -0.1627*** & -0.1391*** & -0.1801*** & -0.1757*** \\
 & 1.25& -0.1658*** & -0.1556*** & -0.1331*** & -0.1692*** & -0.169*** \\
 & 1.5 & -0.163*** & -0.1454*** & -0.1296*** & -0.164*** & -0.1608*** \\
\cline{1-7}
\multirow[t]{5}{*}{$\beta_{distNear}^{\mu}$} & 0.5 & -3.8529*** & -2.3207*** & -2.7967*** & -2.9456*** & -3.8382*** \\
 & 0.75 & 0.8015*** & -0.1358 & -0.1549 & -0.1444 & -0.0426 \\
 & 1.0 & 1.2063*** & 0.5303*** & 1.4025*** & 0.4651*** & 1.1343*** \\
 & 1.25 & 0.8708*** & 0.8543*** & 1.0689*** & 0.7694*** & 0.8666*** \\
 & 1.5 & 0.6391*** & 1.0422*** & 0.7845*** & 0.6879*** & 0.8845*** \\
\cline{1-7}
\multirow[t]{5}{*}{$\beta_{ff}^{\mu}$} & 0.5 & -0.1442*** & -0.1496*** & -0.1927*** & -0.1917*** & -0.1749*** \\
 & 0.75 & -0.1201*** & -0.149*** & -0.1523*** & -0.1603*** & -0.1421*** \\
 & 1.0& -0.1119*** & -0.1504*** & -0.1281*** & -0.1446*** & -0.1093*** \\
 & 1.25 & -0.0973*** & -0.1314*** & -0.12*** & -0.1404*** & -0.0938*** \\
 & 1.5 & -0.0961*** & -0.1161*** & -0.1098*** & -0.1311*** & -0.0985*** \\
\cline{1-7}
\multirow[t]{5}{*}{$\beta_{isGas}^{\mu}$} & 0.5 & -0.2311** & -0.319*** & 0.1141 & -0.3755*** & -0.3989*** \\
 & 0.75 & -0.2443*** & -0.3863*** & 0.1015 & -0.5085*** & -0.4941*** \\
 & 1.0 & -0.2303** & -0.3795*** & 0.1146 & -0.5092*** & -0.493*** \\
 & 1.25& -0.2598*** & -0.3775*** & 0.1121 & -0.4966*** & -0.4525*** \\
 & 1.5 & -0.2414** & -0.3383*** & 0.1231 & -0.4756*** & -0.4155*** \\
\cline{1-7}
\multirow[t]{5}{*}{$\beta_{isWalkHome}^{\mu}$} & 0.5 & 5.6532*** & 5.3508*** & 5.5593*** & 5.5003*** & 5.691*** \\
 & 0.75 & 2.8536*** & 3.6115*** & 3.7672*** & 3.9158*** & 4.005*** \\
 & 1.0 & 2.7698*** & 3.1226*** & 3.1435*** & 3.2018*** & 3.192*** \\
 & 1.25 & 2.9403*** & 3.1045*** & 3.2977*** & 3.2061*** & 3.2324*** \\
 & 1.5 & 3.0992*** & 3.1498*** & 3.4934*** & 3.3087*** & 3.2845*** \\
\cline{1-7}
\multirow[t]{5}{*}{$\beta_{leis}^{\mu}$} & 0.5 & -0.0241 & 0.0216 & -0.1333*** & -0.0832*** & -0.0384 \\
 & 0.75 & -0.0933*** & -0.0416* & -0.1881*** & -0.1181*** & -0.1134*** \\
 & 1.0 & -0.126*** & -0.0662*** & -0.2648*** & -0.1577*** & -0.1477*** \\
 & 1.25 & -0.146*** & -0.0819*** & -0.2593*** & -0.1631*** & -0.1747*** \\
 & 1.5 & -0.159*** & -0.1166*** & -0.2607*** & -0.1748*** & -0.1916*** \\
\cline{1-7}
\multirow[t]{5}{*}{$\beta_{mall}^{\mu}$} & 0.5 & 0.1741*** & 0.2047*** & 0.0848 & 0.5263*** & 0.385*** \\
 & 0.75 & 0.216*** & 0.2582*** & 0.1505*** & 0.5834*** & 0.4625*** \\
 & 1.0 & 0.2603*** & 0.2759*** & 0.1905*** & 0.6065*** & 0.4995*** \\
 & 1.25 & 0.28*** & 0.2754*** & 0.2038*** & 0.6376*** & 0.5029*** \\
 & 1.5 & 0.2711*** & 0.2975*** & 0.197*** & 0.6687*** & 0.4962*** \\
\cline{1-7}
\multirow[t]{5}{*}{$\beta_{outletsFar}^{\mu}$} & 0.5 & 0.1472*** & 0.1365*** & 0.1259*** & 0.1233*** & 0.1012*** \\
 & 0.75 & 0.1249*** & 0.119*** & 0.1099*** & 0.1131*** & 0.0902*** \\
 & 1.0 & 0.1203*** & 0.1108*** & 0.0998*** & 0.0968*** & 0.0808*** \\
 & 1.25 & 0.1153*** & 0.1102*** & 0.0969*** & 0.0971*** & 0.0744*** \\
 & 1.5 & 0.1136*** & 0.1067*** & 0.0958*** & 0.0923*** & 0.0724*** \\
\cline{1-7}
\multirow[t]{5}{*}{$\beta_{outletsNear}^{\mu}$} & 0.5 & 0.0077 & -0.0314 & -0.0417 & -0.0493 & -0.155*** \\
 & 0.75 & 0.4203*** & 0.3154*** & 0.286*** & 0.167*** & 0.0547* \\
 & 1.0 & 0.3627*** & 0.3354*** & 0.2588*** & 0.2853*** & 0.1637*** \\
 & 1.25 & 0.3354*** & 0.2818*** & 0.2608*** & 0.2475*** & 0.1815*** \\
 & 1.5 & 0.3251*** & 0.2681*** & 0.2399*** & 0.2442*** & 0.1857*** \\
\cline{1-7}
\bottomrule
\end{tabular}
}
\caption[Parameter values for level 2 charging with the MNL model, and various utility function thresholds]{Parameter values for level 2 charging with the MNL model, and various utility function thresholds. The amenity threshold is set at 300m and the sample size is 5,000. *** indicates significance at 1\% level, ** significance at 5\% level, and * significance at 10\% level}
    \label{ThirdArticle:Table:ParametersThresholdL2}
\end{table}

\begin{table}[]
    \centering
    \resizebox{0.6\textwidth}{!}{
\begin{tabular}{lllllll}
\toprule
 &  & 0 & 1 & 2 & 3 & 4 \\
Parameter & Distance &  &  &  &  &  \\
\midrule
\multirow[t]{5}{*}{$\beta_{rest}^{\mu}$} & 0.5 & 0.3192*** & 0.3079*** & 0.3368*** & 0.3564*** & 0.3023*** \\
 & 0.75 & 0.3274*** & 0.337*** & 0.3333*** & 0.3715*** & 0.3531*** \\
 & 1.0 & 0.3458*** & 0.3468*** & 0.3735*** & 0.3875*** & 0.359*** \\
 & 1.25 & 0.3512*** & 0.3554*** & 0.3821*** & 0.3663*** & 0.3635*** \\
 & 1.5 & 0.3556*** & 0.3912*** & 0.376*** & 0.3777*** & 0.3956*** \\
\cline{1-7}
\multirow[t]{5}{*}{$\beta_{shop}^{\mu}$} & 0.5 & -0.0557*** & -0.0587*** & 0.0307 & -0.0657*** & -0.0481** \\
 & 0.75 & -0.0798*** & -0.0722*** & 0.0027 & -0.0991*** & -0.108*** \\
 & 1.0 & -0.0904*** & -0.072*** & -0.0312 & -0.1289*** & -0.1207*** \\
 & 1.25 & -0.1044*** & -0.0846*** & -0.051** & -0.1015*** & -0.1376*** \\
 & 1.5 & -0.098*** & -0.0981*** & -0.0496** & -0.1129*** & -0.1494*** \\
\cline{1-7}
\multirow[t]{5}{*}{$\beta_{sm}^{\mu}$} & 0.5 & 0.1579*** & 0.1612*** & 0.0859*** & 0.3017*** & 0.3057*** \\
 & 0.75 & 0.1599*** & 0.1292*** & 0.0871*** & 0.2691*** & 0.296*** \\
 & 1.0 & 0.1221*** & 0.121*** & 0.076** & 0.2711*** & 0.2793*** \\
 & 1.25 & 0.1173*** & 0.1018*** & 0.056 & 0.1942*** & 0.2834*** \\
 & 1.5 & 0.1006*** & 0.1129*** & 0.0437 & 0.1702*** & 0.2572*** \\
\cline{1-7}
\multirow[t]{5}{*}{$\beta_{sport}^{\mu}$} & 0.5 & 0.1544*** & 0.111*** & 0.1451*** & 0.0881*** & 0.0904*** \\
 & 0.75 & 0.0848*** & 0.1064*** & 0.1171*** & 0.0895*** & 0.051** \\
 & 1.0 & 0.0801*** & 0.0838*** & 0.1129*** & 0.0802*** & 0.044* \\
 & 1.25 & 0.0749*** & 0.0566*** & 0.0953*** & 0.0964*** & 0.0425* \\
 & 1.5 & 0.0798*** & 0.0694*** & 0.0765*** & 0.1168*** & 0.0413 \\
\cline{1-7}
\bottomrule
\end{tabular}
}
\caption[Parameter values for level 2 charging with the MNL model, and various utility function thresholds (continued)]{Parameter values for level 2 charging with the MNL model, and various utility function thresholds (continued). The amenity threshold is set at 300m and the sample size is 5,000. *** indicates significance at 1\% level, ** significance at 5\% level, and * significance at 10\% level}
    \label{ThirdArticle:Table:ParametersThresholdL2Continued}
\end{table}

\begin{table}[]
    \centering
    \resizebox{\textwidth}{!}{
\begin{tabular}{llrrrrrrrrrr}
\toprule
 & Distance & \multicolumn{2}{c}{0.5} & \multicolumn{2}{c}{0.75} & \multicolumn{2}{c}{1.0} & \multicolumn{2}{c}{1.25} & \multicolumn{2}{c}{1.5} \\
 &  & Mean & Std & Mean & Std & Mean & Std & Mean & Std & Mean & Std \\
Test set & Score &  &  &  &  &  &  &  &  &  &  \\
\midrule
\multirow[t]{6}{*}{Estimation} & Akaike Information Criterion & 53092.5700 & 616.9900 & 52498.7400 & 653.5600 & 52118.0600 & 698.4400 & 51993.8400 & 728.4000 & 51859.3200 & 668.2300 \\
 & Bayesian Information Criterion & 53188.3300 & 617.0400 & 52594.5000 & 653.6100 & 52213.8100 & 698.5100 & 52089.5900 & 728.4600 & 51944.0400 & 668.2300 \\
 & Final log-likelihood & -26533.2857 & 308.4955 & -26236.3715 & 326.7782 & -26046.0282 & 349.2219 & -25983.9190 & 364.1997 & -25916.6605 & 334.1137 \\
 & Null log-likelihood & -30901.9274 & 50.3773 & -30901.9274 & 50.3773 & -30901.9274 & 50.3773 & -30901.9274 & 50.3773 & -30901.9274 & 50.3773 \\
 & $\rho$ & 0.1414 & 0.0102 & 0.1510 & 0.0108 & 0.1571 & 0.0113 & 0.1591 & 0.0117 & 0.1613 & 0.0107 \\
 & $\bar{\rho}^2$ & 0.1409 & 0.0102 & 0.1506 & 0.0108 & 0.1567 & 0.0113 & 0.1587 & 0.0117 & 0.1609 & 0.0107 \\
\cline{1-12}
\multirow[t]{12}{*}{Validation} & Akaike Information Criterion & 54012.0549 & 126.7889 & 53563.5112 & 135.7686 & 53267.2156 & 135.0333 & 53126.8781 & 145.8170 & 53083.3178 & 157.7764 \\
 & Bayesian Information Criterion & 54096.0359 & 126.7895 & 53647.4921 & 135.7692 & 53351.1965 & 135.0335 & 53210.8591 & 145.8173 & 53167.2988 & 157.7766 \\
 & DPSA final, max & 0.7335 & 0.0557 & 0.6304 & 0.0568 & 0.6179 & 0.0535 & 0.5945 & 0.0372 & 0.5946 & 0.0226 \\
 & DPSA final, mean & 0.0208 & 0.0006 & 0.0211 & 0.0009 & 0.0214 & 0.0010 & 0.0213 & 0.0009 & 0.0213 & 0.0010 \\
 & DPSA final, min & 0.0000 & 0.0000 & 0.0000 & 0.0000 & 0.0000 & 0.0000 & 0.0000 & 0.0000 & 0.0000 & 0.0000 \\
 & DPSA null, max & 0.0000 & 0.0000 & 0.0000 & 0.0000 & 0.0000 & 0.0000 & 0.0000 & 0.0000 & 0.0000 & 0.0000 \\
 & DPSA null, mean & 0.0000 & 0.0000 & 0.0000 & 0.0000 & 0.0000 & 0.0000 & 0.0000 & 0.0000 & 0.0000 & 0.0000 \\
 & DPSA null, min & 0.0000 & 0.0000 & 0.0000 & 0.0000 & 0.0000 & 0.0000 & 0.0000 & 0.0000 & 0.0000 & 0.0000 \\
 & Final log-likelihood & -26993.0275 & 63.3945 & -26768.7556 & 67.8843 & -26620.6078 & 67.5166 & -26550.4391 & 72.9085 & -26528.6589 & 78.8882 \\
 & Null log-likelihood & -30199.1322 & 9.0394 & -30199.1322 & 9.0394 & -30199.1322 & 9.0394 & -30199.1322 & 9.0394 & -30199.1322 & 9.0394 \\
 & $\rho$ & 0.1062 & 0.0021 & 0.1136 & 0.0022 & 0.1185 & 0.0022 & 0.1208 & 0.0024 & 0.1215 & 0.0026 \\
 & $\bar{\rho}^2$ & 0.1057 & 0.0021 & 0.1132 & 0.0022 & 0.1181 & 0.0022 & 0.1204 & 0.0024 & 0.1211 & 0.0026 \\
\cline{1-12}
\bottomrule
\end{tabular}
}
 \caption[Performance indicators for level 2 charging with the MNL model, and various utility function thresholds]{Performance indicators for level 2 charging with the MNL model, and various utility function thresholds. The amenity threshold is set at 300m and the sample size is 5,000.}
    \label{ThirdArticle:Table:ScoresThresholdsL2}
\end{table}


\begin{table}[]
    \centering
    \resizebox{0.6\textwidth}{!}{
\begin{tabular}{lllllll}
\toprule
 &  & 0 & 1 & 2 & 3 & 4 \\
Parameter & Sample size &  &  &  &  &  \\
\midrule
\multirow[t]{4}{*}{$\beta_{distFar}^{\mu}$} & 1000 & -0.1449*** & -0.1306*** & -0.1213*** & -0.1775*** & -0.1624*** \\
 & 3000 & -0.1615*** & -0.1437*** & -0.1263*** & -0.1686*** & -0.1649*** \\
 & 5000 & -0.163*** & -0.1454*** & -0.1296*** & -0.164*** & -0.1608*** \\
 & 7000 & -0.1592*** & -0.1478*** & -0.1299*** & -0.157*** & -0.1625*** \\
\cline{1-7}
\multirow[t]{4}{*}{$\beta_{distNear}^{\mu}$} & 1000 & 0.6025*** & 0.9028*** & 0.9079*** & 0.7305*** & 0.8192*** \\
 & 3000 & 0.6146*** & 1.0851*** & 0.7936*** & 0.7123*** & 0.8924*** \\
 & 5000 & 0.6391*** & 1.0422*** & 0.7845*** & 0.6879*** & 0.8845*** \\
 & 7000 & 0.6525*** & 1.0282*** & 0.806*** & 0.6908*** & 0.8752*** \\
\cline{1-7}
\multirow[t]{4}{*}{$\beta_{ff}^{\mu}$} & 1000 & -0.0392 & -0.1667*** & -0.0084 & -0.2112*** & -0.0269 \\
 & 3000 & -0.1034*** & -0.1426*** & -0.1069*** & -0.1481*** & -0.134*** \\
 & 5000 & -0.0961*** & -0.1161*** & -0.1098*** & -0.1311*** & -0.0985*** \\
 & 7000 & -0.1018*** & -0.1403*** & -0.1041*** & -0.1368*** & -0.1104*** \\
\cline{1-7}
\multirow[t]{4}{*}{$\beta_{isGas}^{\mu}$} & 1000 & -0.3528 & -0.3691 & 0.0597 & -0.8249*** & -0.6176* \\
 & 3000 & -0.3125** & -0.2669* & 0.0686 & -0.5475*** & -0.4784*** \\
 & 5000 & -0.2414** & -0.3383*** & 0.1231 & -0.4756*** & -0.4155*** \\
 & 7000 & -0.3174*** & -0.3788*** & 0.098 & -0.4722*** & -0.4213*** \\
\cline{1-7}
\multirow[t]{4}{*}{$\beta_{isWalkHome}^{\mu}$} & 1000 & 3.002*** & 3.2534*** & 3.6752*** & 3.1939*** & 3.3787*** \\
 & 3000 & 3.011*** & 3.1172*** & 3.4777*** & 3.2558*** & 3.2743*** \\
 & 5000 & 3.0992*** & 3.1498*** & 3.4934*** & 3.3087*** & 3.2845*** \\
 & 7000 & 3.1229*** & 3.1425*** & 3.561*** & 3.3311*** & 3.2716*** \\
\cline{1-7}
\multirow[t]{4}{*}{$\beta_{leis}^{\mu}$} & 1000 & -0.1944*** & -0.0048 & -0.3333*** & -0.1536*** & -0.1785*** \\
 & 3000 & -0.1494*** & -0.1163*** & -0.2759*** & -0.1645*** & -0.185*** \\
 & 5000 & -0.159*** & -0.1166*** & -0.2607*** & -0.1748*** & -0.1916*** \\
 & 7000 & -0.153*** & -0.1279*** & -0.2453*** & -0.1823*** & -0.1839*** \\
\cline{1-7}
\multirow[t]{4}{*}{$\beta_{mall}^{\mu}$} & 1000 & 0.3014*** & 0.1744 & 0.1894 & 0.7692*** & 0.4778*** \\
 & 3000 & 0.2275*** & 0.3175*** & 0.1614** & 0.6539*** & 0.5013*** \\
 & 5000 & 0.2711*** & 0.2975*** & 0.197*** & 0.6687*** & 0.4962*** \\
 & 7000 & 0.2756*** & 0.2978*** & 0.1829*** & 0.6668*** & 0.4985*** \\
\cline{1-7}
\multirow[t]{4}{*}{$\beta_{outletsFar}^{\mu}$} & 1000 & 0.0991*** & 0.0999*** & 0.0861*** & 0.0845*** & 0.0768*** \\
 & 3000 & 0.1023*** & 0.1147*** & 0.0965*** & 0.087*** & 0.078*** \\
 & 5000 & 0.1136*** & 0.1067*** & 0.0958*** & 0.0923*** & 0.0724*** \\
 & 7000 & 0.1091*** & 0.1027*** & 0.1076*** & 0.0925*** & 0.0764*** \\
\cline{1-7}
\multirow[t]{4}{*}{$\beta_{outletsNear}^{\mu}$} & 1000 & 0.3372*** & 0.3*** & 0.1892*** & 0.2854*** & 0.1459*** \\
 & 3000 & 0.3411*** & 0.2747*** & 0.2403*** & 0.2348*** & 0.1898*** \\
 & 5000 & 0.3251*** & 0.2681*** & 0.2399*** & 0.2442*** & 0.1857*** \\
 & 7000 & 0.3274*** & 0.2725*** & 0.2459*** & 0.2362*** & 0.1894*** \\
\cline{1-7}
\bottomrule
\end{tabular}
}
\caption[Parameter values for level 2 charging with the MNL model, and various sample sizes]{Parameter values for level 2 charging with the MNL model, and various sample sizes. The amenity threshold is set at 300m and the utility function threshold is 1,5km. *** indicates significance at 1\% level, ** significance at 5\% level, and * significance at 10\% level}
    \label{ThirdArticle:Table:ParametersDistancesL2Sample5000}
\end{table}

\begin{table}[]
    \centering
    \resizebox{0.6\textwidth}{!}{
\begin{tabular}{lllllll}
\toprule
 &  & 0 & 1 & 2 & 3 & 4 \\
Parameter & Sample size &  &  &  &  &  \\
\midrule
\multirow[t]{4}{*}{$\beta_{rest}^{\mu}$} & 1000 & 0.3016*** & 0.361*** & 0.4669*** & 0.3992*** & 0.3349*** \\
 & 3000 & 0.3439*** & 0.4012*** & 0.3932*** & 0.3741*** & 0.4106*** \\
 & 5000 & 0.3556*** & 0.3912*** & 0.376*** & 0.3777*** & 0.3956*** \\
 & 7000 & 0.362*** & 0.3917*** & 0.3821*** & 0.3659*** & 0.3892*** \\
\cline{1-7}
\multirow[t]{4}{*}{$\beta_{shop}^{\mu}$} & 1000 & -0.1368*** & -0.0886* & -0.0983* & -0.1189** & -0.1306*** \\
 & 3000 & -0.0952*** & -0.1053*** & -0.0405 & -0.0872*** & -0.138*** \\
 & 5000 & -0.098*** & -0.0981*** & -0.0496** & -0.1129*** & -0.1494*** \\
 & 7000 & -0.0988*** & -0.096*** & -0.0509*** & -0.1044*** & -0.1507*** \\
\cline{1-7}
\multirow[t]{4}{*}{$\beta_{sm}^{\mu}$} & 1000 & 0.2847*** & 0.2098*** & -0.0212 & 0.2968*** & 0.3009*** \\
 & 3000 & 0.1377*** & 0.1048** & 0.056 & 0.1673*** & 0.2358*** \\
 & 5000 & 0.1006*** & 0.1129*** & 0.0437 & 0.1702*** & 0.2572*** \\
 & 7000 & 0.0922*** & 0.1252*** & 0.026 & 0.1392*** & 0.2698*** \\
\cline{1-7}
\multirow[t]{4}{*}{$\beta_{sport}^{\mu}$} & 1000 & 0.0518 & 0.0961* & 0.1639*** & 0.1914*** & 0.0185 \\
 & 3000 & 0.0706** & 0.0936*** & 0.0882*** & 0.0938*** & 0.071** \\
 & 5000 & 0.0798*** & 0.0694*** & 0.0765*** & 0.1168*** & 0.0413 \\
 & 7000 & 0.0589*** & 0.0807*** & 0.0931*** & 0.1125*** & 0.027 \\
\cline{1-7}
\bottomrule
\end{tabular}
}
\caption[Parameter values for level 2 charging with the MNL model, and various sample sizes (continued)]{Parameter values for level 2 charging with the MNL model, and various sample sizes. The amenity threshold is set at 300m and the utility function threshold is 1,5km (continued). *** indicates significance at 1\% level, ** significance at 5\% level, and * significance at 10\% level}
    \label{ThirdArticle:Table:ParametersDistancesL2Sample5000Continued}
\end{table}

\begin{table}[]
    \centering
    \resizebox{\textwidth}{!}{
\begin{tabular}{llrrrrrrrr}
\toprule
 & Sample size & \multicolumn{2}{c}{1000} & \multicolumn{2}{c}{3000} & \multicolumn{2}{c}{5000} & \multicolumn{2}{c}{7000} \\
 &  & Mean & Std & Mean & Std & Mean & Std & Mean & Std \\
Test set & Score &  &  &  &  &  &  &  &  \\
\midrule
\multirow[t]{6}{*}{Estimation} & Akaike Information Criterion & 10480.7700 & 240.7100 & 31224.4800 & 442.5100 & 51859.3200 & 668.2300 & 72541.2300 & 820.6900 \\
 & Bayesian Information Criterion & 10544.5700 & 240.7100 & 31302.5700 & 442.5100 & 51944.0400 & 668.2300 & 72630.3200 & 820.6900 \\
 & Final log-likelihood & -5227.3846 & 120.3561 & -15599.2423 & 221.2528 & -25916.6605 & 334.1137 & -36257.6129 & 410.3451 \\
 & Null log-likelihood & -6185.5388 & 15.4519 & -18539.4485 & 35.4542 & -30901.9274 & 50.3773 & -43262.0048 & 79.5915 \\
 & $\rho$ & 0.1549 & 0.0191 & 0.1586 & 0.0122 & 0.1613 & 0.0107 & 0.1619 & 0.0091 \\
 & $\bar{\rho}^2$ & 0.1528 & 0.0191 & 0.1579 & 0.0122 & 0.1609 & 0.0107 & 0.1616 & 0.0091 \\
\cline{1-10}
\multirow[t]{12}{*}{Validation} & Akaike Information Criterion & 53148.7700 & 179.0100 & 53081.9700 & 161.1900 & 53083.3200 & 157.7800 & 53084.3000 & 137.1300 \\
 & Bayesian Information Criterion & 53232.7500 & 179.0100 & 53165.9500 & 161.1900 & 53167.3000 & 157.7800 & 53168.2800 & 137.1300 \\
 & DPSA final, max & 0.5802 & 0.0412 & 0.5903 & 0.0300 & 0.5946 & 0.0226 & 0.5969 & 0.0204 \\
 & DPSA final, mean & 0.0207 & 0.0014 & 0.0208 & 0.0009 & 0.0213 & 0.0010 & 0.0215 & 0.0009 \\
 & DPSA final, min & 0.0000 & 0.0000 & 0.0000 & 0.0000 & 0.0000 & 0.0000 & 0.0000 & 0.0000 \\
 & DPSA null, max & 0.0038 & 0.0000 & 0.0038 & 0.0000 & 0.0038 & 0.0000 & 0.0038 & 0.0000 \\
 & DPSA null, mean & 0.0017 & 0.0000 & 0.0017 & 0.0000 & 0.0017 & 0.0000 & 0.0017 & 0.0000 \\
 & DPSA null, min & 0.0014 & 0.0000 & 0.0014 & 0.0000 & 0.0014 & 0.0000 & 0.0014 & 0.0000 \\
 & Final log-likelihood & -26561.3830 & 89.5047 & -26527.9849 & 80.5963 & -26528.6589 & 78.8882 & -26529.1505 & 68.5667 \\
 & Null log-likelihood & -30199.1322 & 9.0394 & -30199.1322 & 9.0394 & -30199.1322 & 9.0394 & -30199.1322 & 9.0394 \\
 & $\rho$ & 0.1205 & 0.0029 & 0.1216 & 0.0027 & 0.1215 & 0.0026 & 0.1215 & 0.0023 \\
 & $\bar{\rho}^2$ & 0.1200 & 0.0029 & 0.1211 & 0.0027 & 0.1211 & 0.0026 & 0.1211 & 0.0023 \\
\cline{1-10}
\bottomrule
\end{tabular}
}

\caption[Performance indicators for level 2 charging with the MNL model, and various sample sizes]{Performance indicators for level 2 charging with the MNL model, and various sample sizes. The amenity threshold is set at 300m and the utility function threshold is 1,5km. *** indicates significance at 1\% level, ** significance at 5\% level, and * significance at 10\% level}
    \label{ThirdArticle:Table:ScoresDistancesL2Sample5000}
\end{table}

\begin{table}
\centering
    \resizebox{0.6\textwidth}{!}{
\begin{tabular}{lllllll}
\toprule
 &  & 0 & 1 & 2 & 3 & 4 \\
Parameter & Distance &  &  &  &  &  \\
\midrule
\multirow[t]{4}{*}{$\beta_{distFar}^{\mu}$} & 100 & -0.176*** & -0.1741*** & -0.1832*** & -0.1832*** & -0.1837*** \\
 & 200 & -0.1737*** & -0.1715*** & -0.1798*** & -0.1791*** & -0.1806*** \\
 & 300 & -0.1739*** & -0.1738*** & -0.1807*** & -0.1819*** & -0.1823*** \\
 & 400 & -0.1742*** & -0.1733*** & -0.1817*** & -0.1797*** & -0.1831*** \\
\cline{1-7}
\multirow[t]{4}{*}{$\beta_{distNear}^{\mu}$} & 100 & 1.411*** & 1.743*** & 1.322*** & 1.3636*** & 1.9477*** \\
 & 200 & 1.2063*** & 1.4573*** & 1.1227*** & 1.1926*** & 1.7204*** \\
 & 300 & 1.1971*** & 1.4446*** & 1.0276*** & 1.0765*** & 1.6485*** \\
 & 400 & 1.2011*** & 1.508*** & 1.0773*** & 1.1526*** & 1.6488*** \\
\cline{1-7}
\multirow[t]{4}{*}{$\beta_{ff}^{\mu}$} & 100 & -0.3821*** & -0.6857*** & -0.4282*** & -0.3854*** & -0.5555*** \\
 & 200 & -0.4237*** & -0.5366*** & -0.4259*** & -0.266*** & -0.3712*** \\
 & 300 & -0.0479 & -0.1621*** & 0.1073 & 0.0505 & 0.0409 \\
 & 400 & 0.2099*** & 0.0838 & 0.329*** & 0.3817*** & 0.2778*** \\
\cline{1-7}
\multirow[t]{4}{*}{$\beta_{isGas}^{\mu}$} & 100 & 0.1049 & 0.0454 & -0.1722** & 0.0253 & 0.0037 \\
 & 200 & -0.1824*** & -0.1933*** & -0.4185*** & -0.2357*** & -0.2465*** \\
 & 300 & -0.3423*** & -0.3846*** & -0.5826*** & -0.4111*** & -0.4357*** \\
 & 400 & -0.1634*** & -0.1328** & -0.4461*** & -0.2873*** & -0.2091*** \\
\cline{1-7}
\multirow[t]{4}{*}{$\beta_{isWalkHome}^{\mu}$} & 100 & 2.0464*** & 1.7839*** & 1.8796*** & 1.4404*** & 2.0454*** \\
 & 200 & 1.9622*** & 1.6666*** & 1.7922*** & 1.3722*** & 1.9548*** \\
 & 300 & 1.9035*** & 1.5773*** & 1.6888*** & 1.2823*** & 1.8437*** \\
 & 400 & 1.9187*** & 1.6313*** & 1.7537*** & 1.3088*** & 1.8588*** \\
\cline{1-7}
\multirow[t]{4}{*}{$\beta_{leis}^{\mu}$} & 100 & 0.2584*** & 0.3385*** & 0.1879** & 0.3869*** & 0.2217*** \\
 & 200 & 0.1551*** & 0.2365*** & 0.083* & 0.2364*** & 0.1826*** \\
 & 300 & 0.3323*** & 0.4502*** & 0.2744*** & 0.3302*** & 0.3286*** \\
 & 400 & 0.0084 & 0.1055*** & -0.0594 & -0.0429 & 0.0155 \\
\cline{1-7}
\multirow[t]{4}{*}{$\beta_{mall}^{\mu}$} & 100 & 0.3925*** & 0.4002*** & 0.51*** & 0.4583*** & 0.4416*** \\
 & 200 & 0.2075*** & 0.1304* & 0.4064*** & 0.3938*** & 0.3062*** \\
 & 300 & 0.6607*** & 0.6946*** & 0.8865*** & 0.7802*** & 0.7251*** \\
 & 400 & 0.8997*** & 0.9646*** & 1.0536*** & 0.9838*** & 0.9663*** \\
\cline{1-7}
\multirow[t]{4}{*}{$\beta_{outletsFar}^{\mu}$} & 100 & 0.3965*** & 0.4726*** & 0.4457*** & 0.4166*** & 0.494*** \\
 & 200 & 0.2905*** & 0.2857*** & 0.326*** & 0.2931*** & 0.3479*** \\
 & 300 & 0.2022*** & 0.1823*** & 0.1993*** & 0.1692*** & 0.2174*** \\
 & 400 & 0.3001*** & 0.3195*** & 0.3038*** & 0.297*** & 0.3428*** \\
\cline{1-7}
\multirow[t]{4}{*}{$\beta_{outletsNear}^{\mu}$} & 100 & 0.0506 & -0.0115 & 0.2013*** & 0.2003** & -0.1789** \\
 & 200 & 0.0126 & -0.0773 & 0.1754** & 0.1698* & -0.2252*** \\
 & 300 & -0.0107 & -0.1007 & 0.1431* & 0.1367 & -0.2777*** \\
 & 400 & 0.0371 & -0.0357 & 0.1799** & 0.1862** & -0.2063*** \\
\cline{1-7}
\bottomrule
\end{tabular}
}
\caption[Parameter values for level 3 charging with the MNL model, and various amenity thresholds]{Parameter values for level 3 charging with the MNL model, and various amenity thresholds. The utility function threshold is set to 1,5km. *** indicates significance at 1\% level, ** significance at 5\% level, and * significance at 10\% level}
    \label{ThirdArticle:Table:ParametersDistancesL3}
\end{table}

\begin{table}
\centering
    \resizebox{0.6\textwidth}{!}{
\begin{tabular}{lllllll}
\toprule
 &  & 0 & 1 & 2 & 3 & 4 \\
Parameter & Distance &  &  &  &  &  \\
\midrule
\multirow[t]{4}{*}{$\beta_{rest}^{\mu}$} & 100 & -0.1871*** & -0.1848*** & -0.2134*** & -0.2781*** & -0.2602*** \\
 & 200 & -0.2837*** & -0.2687*** & -0.3107*** & -0.401*** & -0.4175*** \\
 & 300 & -0.1182*** & -0.0878* & -0.1855*** & -0.204*** & -0.2553*** \\
 & 400 & -0.1324** & -0.0559 & -0.1809*** & -0.2868*** & -0.2071*** \\
\cline{1-7}
\multirow[t]{4}{*}{$\beta_{shop}^{\mu}$} & 100 & -0.3348*** & -0.3325*** & -0.3579*** & -0.4193*** & -0.3614*** \\
 & 200 & 0.1204*** & 0.1095*** & 0.1818*** & 0.0996** & 0.1475*** \\
 & 300 & -0.2306*** & -0.3053*** & -0.2609*** & -0.2551*** & -0.2172*** \\
 & 400 & -0.4168*** & -0.4875*** & -0.4565*** & -0.3506*** & -0.4777*** \\
\cline{1-7}
\multirow[t]{4}{*}{$\beta_{sm}^{\mu}$} & 100 & 0.5238 & 1.0688*** & 1.2679*** & 0.9772*** & 0.7517** \\
 & 200 & 0.0131 & 0.0393 & -0.0158 & 0.0178 & -0.0521 \\
 & 300 & 0.2098*** & 0.4397*** & 0.333*** & 0.43*** & 0.306*** \\
 & 400 & 0.8362*** & 0.9826*** & 0.9871*** & 0.9231*** & 1.0172*** \\
\cline{1-7}
\multirow[t]{4}{*}{$\beta_{sport}^{\mu}$} & 100 & -0.4973*** & -0.7591*** & -0.5213*** & -0.6338*** & -0.6354*** \\
 & 200 & -0.3129*** & -0.4919*** & -0.2887*** & -0.3342*** & -0.3847*** \\
 & 300 & -0.2794*** & -0.4374*** & -0.3163*** & -0.2945*** & -0.337*** \\
 & 400 & -0.1467*** & -0.3012*** & -0.0911** & -0.146*** & -0.2263*** \\
\cline{1-7}
\bottomrule
\end{tabular}
}
\caption[Parameter values for level 3 charging with the MNL model, and various amenity thresholds (continued)]{Parameter values for level 3 charging with the MNL model, and various amenity thresholds. The utility function threshold is set to 1,5km (continued). *** indicates significance at 1\% level, ** significance at 5\% level, and * significance at 10\% level}
    \label{ThirdArticle:Table:ParametersDistancesL3Continued}
\end{table}

\begin{table}[]
    \centering
    \resizebox{\textwidth}{!}{
\begin{tabular}{llrrrrrrrr}
\toprule
 & Distance & \multicolumn{2}{c}{100} & \multicolumn{2}{c}{200} & \multicolumn{2}{c}{300} & \multicolumn{2}{c}{400} \\
 &  & Mean & Std & Mean & Std & Mean & Std & Mean & Std \\
Test set & Score &  &  &  &  &  &  &  &  \\
\midrule
\multirow[t]{6}{*}{Estimation} & Akaike Information Criterion & 10693.9800 & 88.9800 & 10606.4300 & 97.6900 & 10735.1200 & 105.4700 & 10655.0700 & 113.7000 \\
 & Bayesian Information Criterion & 10773.7300 & 88.9800 & 10686.1900 & 97.7200 & 10814.8700 & 105.5200 & 10734.8200 & 113.7200 \\
 & Final log-likelihood & -5333.9900 & 44.4899 & -5290.2174 & 48.8472 & -5354.5580 & 52.7369 & -5314.5332 & 56.8513 \\
 & Null log-likelihood & -8479.0167 & 114.4676 & -8479.0167 & 114.4676 & -8479.0167 & 114.4676 & -8479.0167 & 114.4676 \\
 & $\rho$ & 0.3708 & 0.0083 & 0.3760 & 0.0076 & 0.3684 & 0.0069 & 0.3731 & 0.0081 \\
 & $\bar{\rho}^2$ & 0.3693 & 0.0083 & 0.3745 & 0.0076 & 0.3669 & 0.0069 & 0.3716 & 0.0081 \\
\cline{1-10}
\multirow[t]{12}{*}{Validation} & Akaike Information Criterion & 825.2300 & 44.2000 & 812.2000 & 44.7700 & 823.3700 & 45.3500 & 816.9000 & 48.3000 \\
 & Bayesian Information Criterion & 868.9600 & 44.1800 & 855.9400 & 44.7500 & 867.1100 & 45.3400 & 860.6400 & 48.2800 \\
 & DPSA final, max & 0.9549 & 0.0053 & 0.9461 & 0.0065 & 0.9462 & 0.0117 & 0.9497 & 0.0117 \\
 & DPSA final, mean & 0.2884 & 0.0190 & 0.2974 & 0.0188 & 0.2910 & 0.0179 & 0.2938 & 0.0210 \\
 & DPSA final, min & 0.0011 & 0.0006 & 0.0016 & 0.0009 & 0.0015 & 0.0008 & 0.0016 & 0.0009 \\
 & DPSA null, max & 0.2500 & 0.0000 & 0.2500 & 0.0000 & 0.2500 & 0.0000 & 0.2500 & 0.0000 \\
 & DPSA null, mean & 0.0831 & 0.0018 & 0.0831 & 0.0018 & 0.0831 & 0.0018 & 0.0831 & 0.0018 \\
 & DPSA null, min & 0.0556 & 0.0000 & 0.0556 & 0.0000 & 0.0556 & 0.0000 & 0.0556 & 0.0000 \\
 & Final log-likelihood & -399.6134 & 22.0993 & -393.1023 & 22.3839 & -398.6860 & 22.6755 & -395.4509 & 24.1482 \\
 & Null log-likelihood & -548.7335 & 2.7691 & -548.7335 & 2.7691 & -548.7335 & 2.7691 & -548.7335 & 2.7691 \\
 & $\rho$ & 0.2717 & 0.0402 & 0.2836 & 0.0406 & 0.2735 & 0.0409 & 0.2793 & 0.0438 \\
 & $\bar{\rho}^2$ & 0.2481 & 0.0402 & 0.2599 & 0.0406 & 0.2498 & 0.0409 & 0.2556 & 0.0438 \\
\cline{1-10}
\bottomrule
\end{tabular}
}

\caption[Performance indicators for level 3 charging with the MNL model, and various amenity thresholds]{Performance indicators for level 3 charging with the MNL model, and various amenity thresholds. The utility function threshold is set to 1,5km. *** indicates significance at 1\% level, ** significance at 5\% level, and * significance at 10\% level}
    \label{ThirdArticle:Table:ScoresDistances1500L3}
\end{table}


\begin{table}[]
    \centering
    \resizebox{0.6\textwidth}{!}{
\begin{tabular}{lllllll}
\toprule
 &  & 0 & 1 & 2 & 3 & 4 \\
Parameter & Distance &  &  &  &  &  \\
\midrule
\multirow[t]{5}{*}{$\beta_{distFar}^{\mu}$} & 0.5 & -0.1935*** & -0.1915*** & -0.198*** & -0.2002*** & -0.2004*** \\
 & 0.75 & -0.1858*** & -0.184*** & -0.1941*** & -0.1934*** & -0.1928*** \\
 & 1.0 & -0.182*** & -0.1821*** & -0.1897*** & -0.1897*** & -0.1891*** \\
 & 1.25 & -0.1812*** & -0.1787*** & -0.1856*** & -0.1837*** & -0.1866*** \\
 & 1.5 & -0.1737*** & -0.1715*** & -0.1798*** & -0.1791*** & -0.1806*** \\
\cline{1-7}
\multirow[t]{5}{*}{$\beta_{distNear}^{\mu}$} & 0.5 & -6.9284 & 0.2229 & -1.5196 & -0.0655 & 0.927 \\
 & 0.75 & 1.0473 & 1.606*** & 3.2752*** & 1.035 & 0.7309 \\
 & 1.0 & 2.4474*** & 1.9191*** & 2.1093*** & 1.7397*** & 2.4352*** \\
 & 1.25 & 0.676** & 1.3522*** & 1.1083*** & 1.9382*** & 1.5369*** \\
 & 1.5 & 1.2063*** & 1.4573*** & 1.1227*** & 1.1926*** & 1.7204*** \\
\cline{1-7}
\multirow[t]{5}{*}{$\beta_{ff}^{\mu}$} & 0.5 & -0.3137*** & -0.3711*** & -0.3261*** & -0.0649 & -0.2145*** \\
 & 0.75 & -0.409*** & -0.4809*** & -0.3713*** & -0.145* & -0.3271*** \\
 & 1.0 & -0.4871*** & -0.5187*** & -0.4233*** & -0.234*** & -0.4081*** \\
 & 1.25 & -0.4619*** & -0.5708*** & -0.448*** & -0.2716*** & -0.4395*** \\
 & 1.5 & -0.4237*** & -0.5366*** & -0.4259*** & -0.266*** & -0.3712*** \\
\cline{1-7}
\multirow[t]{5}{*}{$\beta_{isGas}^{\mu}$} & 0.5 & -0.1969*** & -0.1714*** & -0.412*** & -0.215*** & -0.2436*** \\
 & 0.75 & -0.1866*** & -0.1533*** & -0.4001*** & -0.2062*** & -0.2324*** \\
 & 1.0 & -0.1954*** & -0.1842*** & -0.4277*** & -0.2400*** & -0.2675*** \\
 & 1.25 & -0.1831*** & -0.191*** & -0.4147*** & -0.2197*** & -0.2533*** \\
 & 1.5 & -0.1824*** & -0.1933*** & -0.4185*** & -0.2357*** & -0.2465*** \\
\cline{1-7}
\multirow[t]{5}{*}{$\beta_{isWalkHome}^{\mu}$} & 0.5 & 5.0925*** & 1.5596 & 1.3657 & 2.2746 & 1.9765 \\
 & 0.75 & 2.0835*** & 2.2644*** & 2.4458*** & 2.1489*** & 2.7706*** \\
 & 1.0 & 2.1762*** & 1.6867*** & 1.7287*** & 1.6861*** & 1.8825*** \\
 & 1.25 & 2.0271*** & 1.7278*** & 1.8199*** & 1.4059*** & 1.9651*** \\
 & 1.5 & 1.9622*** & 1.6666*** & 1.7922*** & 1.3722*** & 1.9548*** \\
\cline{1-7}
\multirow[t]{5}{*}{$\beta_{leis}^{\mu}$} & 0.5 & 0.2254*** & 0.3259*** & 0.0838* & 0.3032*** & 0.2357*** \\
 & 0.75 & 0.1603*** & 0.2512*** & 0.0655 & 0.2482*** & 0.1637*** \\
 & 1.0 & 0.1541*** & 0.2312*** & 0.0688 & 0.2346*** & 0.1586*** \\
 & 1.25 & 0.164*** & 0.2443*** & 0.0873* & 0.2495*** & 0.1824*** \\
 & 1.5 & 0.1551*** & 0.2365*** & 0.083* & 0.2364*** & 0.1826*** \\
\cline{1-7}
\multirow[t]{5}{*}{$\beta_{mall}^{\mu}$} & 0.5 & 0.1487* & 0.1024 & 0.398*** & 0.378*** & 0.2917*** \\
 & 0.75 & 0.1426* & 0.0967 & 0.3991*** & 0.3692*** & 0.2826*** \\
 & 1.0 & 0.1358* & 0.0847 & 0.3637*** & 0.3512*** & 0.2548*** \\
 & 1.25 & 0.1729** & 0.0943 & 0.3861*** & 0.3862*** & 0.2687*** \\
 & 1.5 & 0.2075*** & 0.1304* & 0.4064*** & 0.3938*** & 0.3062*** \\
\cline{1-7}
\multirow[t]{5}{*}{$\beta_{outletsFar}^{\mu}$} & 0.5 & 0.2343*** & 0.2247*** & 0.2872*** & 0.2343*** & 0.2886*** \\
 & 0.75 & 0.2517*** & 0.2514*** & 0.2987*** & 0.2485*** & 0.3119*** \\
 & 1.0 & 0.2695*** & 0.2576*** & 0.2994*** & 0.26*** & 0.3255*** \\
 & 1.25 & 0.2702*** & 0.2648*** & 0.3075*** & 0.2797*** & 0.3354*** \\
 & 1.5 & 0.2905*** & 0.2857*** & 0.326*** & 0.2931*** & 0.3479*** \\
\cline{1-7}
\multirow[t]{5}{*}{$\beta_{outletsNear}^{\mu}$} & 0.5 & 0.1295 & 0.3088 & 1.5361 & 0.2389 & 0.2481 \\
 & 0.75 & 0.0711 & -0.2416 & -0.3532 & -0.0263 & -0.162 \\
 & 1.0 & -0.3237** & -0.2531* & 0.0082 & -0.1406 & -0.3656*** \\
 & 1.25 & 0.0723 & -0.1202 & 0.1297 & -0.0614 & -0.2175** \\
 & 1.5 & 0.0126 & -0.0773 & 0.1754** & 0.1698* & -0.2252*** \\
\cline{1-7}
\bottomrule
\end{tabular}
}
\caption[Parameter values for level 3 charging with the MNL model, and various utility function thresholds]{Parameter values for level 3 charging with the MNL model, and various utility function thresholds. The amenity threshold is set to 200m. *** indicates significance at 1\% level, ** significance at 5\% level, and * significance at 10\% level}
    \label{ThirdArticle:Table:ParametersThresholdsL3}
\end{table}

\begin{table}[]
    \centering
    \resizebox{0.6\textwidth}{!}{
\begin{tabular}{lllllll}
\toprule
 &  & 0 & 1 & 2 & 3 & 4 \\
Parameter & Distance &  &  &  &  &  \\
\midrule
\multirow[t]{5}{*}{$\beta_{rest}^{\mu}$} & 0.5 & -0.3397*** & -0.3278*** & -0.3376*** & -0.4622*** & -0.4658*** \\
 & 0.75 & -0.3234*** & -0.2961*** & -0.3353*** & -0.4532*** & -0.4447*** \\
 & 1.0 & -0.2814*** & -0.2873*** & -0.3237*** & -0.4247*** & -0.4209*** \\
 & 1.25 & -0.2738*** & -0.2495*** & -0.2926*** & -0.3968*** & -0.3898*** \\
 & 1.5 & -0.2837*** & -0.2687*** & -0.3107*** & -0.401*** & -0.4175*** \\
\cline{1-7}
\multirow[t]{5}{*}{$\beta_{shop}^{\mu}$} & 0.5 & 0.1114*** & 0.0837** & 0.2038*** & 0.0527 & 0.1318*** \\
 & 0.75 & 0.1562*** & 0.1185*** & 0.211*** & 0.0999** & 0.1733*** \\
 & 1.0 & 0.1593*** & 0.1298*** & 0.2143*** & 0.1215*** & 0.1841*** \\
 & 1.25 & 0.1278*** & 0.1146*** & 0.1824*** & 0.1038*** & 0.1556*** \\
 & 1.5 & 0.1204*** & 0.1095*** & 0.1818*** & 0.0996** & 0.1475*** \\
\cline{1-7}
\multirow[t]{5}{*}{$\beta_{sm}^{\mu}$} & 0.5 & 0.0736 & 0.1047 & 0.0131 & 0.1807* & -0.0009 \\
 & 0.75 & 0.0795 & 0.1253 & 0.0237 & 0.1745* & -0.0027 \\
 & 1.0 & 0.0137 & 0.0752 & 0.0063 & 0.0998 & -0.0743 \\
 & 1.25 & 0.0097 & 0.0386 & -0.0262 & 0.0294 & -0.0782 \\
 & 1.5 & 0.0131 & 0.0393 & -0.0158 & 0.0178 & -0.0521 \\
\cline{1-7}
\multirow[t]{5}{*}{$\beta_{sport}^{\mu}$} & 0.5 & -0.1829*** & -0.3273*** & -0.1623*** & -0.1981*** & -0.2401*** \\
 & 0.75 & -0.225*** & -0.3769*** & -0.1746*** & -0.2331*** & -0.282*** \\
 & 1.0 & -0.2736*** & -0.3935*** & -0.2045*** & -0.276*** & -0.3285*** \\
 & 1.25 & -0.3075*** & -0.4727*** & -0.2753*** & -0.3279*** & -0.3778*** \\
 & 1.5 & -0.3129*** & -0.4919*** & -0.2887*** & -0.3342*** & -0.3847*** \\
\cline{1-7}
\bottomrule
\end{tabular}
}
\caption[Parameter values for level 3 charging with the MNL model, and various utility function thresholds (continued)]{Parameter values for level 3 charging with the MNL model, and various utility function thresholds (continued). The amenity threshold is set to 200m. *** indicates significance at 1\% level, ** significance at 5\% level, and * significance at 10\% level}
    \label{ThirdArticle:Table:ParametersThresholdsL3Continued}
\end{table}

\begin{table}[]
    \centering
    \resizebox{\textwidth}{!}{
\begin{tabular}{llrrrrrrrrrr}
\toprule
 & Distance & \multicolumn{2}{c}{0.5} & \multicolumn{2}{c}{0.75} & \multicolumn{2}{c}{1.0} & \multicolumn{2}{c}{1.25} & \multicolumn{2}{c}{1.5} \\
 &  & Mean & Std & Mean & Std & Mean & Std & Mean & Std & Mean & Std \\
Test set & Score &  &  &  &  &  &  &  &  &  &  \\
\midrule
\multirow[t]{6}{*}{Estimation} & Akaike Information Criterion & 10909.0800 & 90.5500 & 10791.4700 & 81.5900 & 10758.3700 & 96.0200 & 10704.4000 & 94.7800 & 10606.4300 & 97.6900 \\
 & Bayesian Information Criterion & 10988.8300 & 90.6100 & 10871.2300 & 81.6500 & 10838.1200 & 96.0600 & 10784.1600 & 94.7800 & 10686.1900 & 97.7200 \\
 & Final log-likelihood & -5441.5394 & 45.2774 & -5382.7351 & 40.7970 & -5366.1834 & 48.0075 & -5339.2006 & 47.3914 & -5290.2174 & 48.8472 \\
 & Null log-likelihood & -8479.0167 & 114.4676 & -8479.0167 & 114.4676 & -8479.0167 & 114.4676 & -8479.0167 & 114.4676 & -8479.0167 & 114.4676 \\
 & $\rho$ & 0.3582 & 0.0067 & 0.3651 & 0.0065 & 0.3671 & 0.0070 & 0.3702 & 0.0085 & 0.3760 & 0.0076 \\
 & $\bar{\rho}^2$ & 0.3566 & 0.0068 & 0.3636 & 0.0065 & 0.3655 & 0.0070 & 0.3687 & 0.0085 & 0.3745 & 0.0076 \\
\cline{1-12}
\multirow[t]{12}{*}{Validation} & Akaike Information Criterion & 832.1997 & 38.0500 & 826.7232 & 35.6112 & 820.5784 & 37.9845 & 820.2844 & 42.8279 & 812.2046 & 44.7679 \\
 & Bayesian Information Criterion & 875.9330 & 38.0377 & 870.4565 & 35.5982 & 864.3118 & 37.9693 & 864.0177 & 42.8114 & 855.9379 & 44.7517 \\
 & DPSA final, max & 0.9563 & 0.0241 & 0.9555 & 0.0137 & 0.9551 & 0.0130 & 0.9441 & 0.0092 & 0.9461 & 0.0065 \\
 & DPSA final, mean & 0.2820 & 0.0177 & 0.2842 & 0.0150 & 0.2885 & 0.0173 & 0.2889 & 0.0173 & 0.2974 & 0.0188 \\
 & DPSA final, min & 0.0009 & 0.0004 & 0.0012 & 0.0005 & 0.0013 & 0.0007 & 0.0014 & 0.0007 & 0.0016 & 0.0009 \\
 & DPSA null, max & 0.2500 & 0.0000 & 0.2500 & 0.0000 & 0.2500 & 0.0000 & 0.2500 & 0.0000 & 0.2500 & 0.0000 \\
 & DPSA null, mean & 0.0800 & 0.0000 & 0.0800 & 0.0000 & 0.0800 & 0.0000 & 0.0800 & 0.0000 & 0.0800 & 0.0000 \\
 & DPSA null, min & 0.0600 & 0.0000 & 0.0600 & 0.0000 & 0.0600 & 0.0000 & 0.0600 & 0.0000 & 0.0600 & 0.0000 \\
 & Final log-likelihood & -403.0998 & 19.0250 & -400.3616 & 17.8056 & -397.2892 & 18.9923 & -397.1422 & 21.4140 & -393.1023 & 22.3839 \\
 & Null log-likelihood & -548.7335 & 2.7691 & -548.7335 & 2.7691 & -548.7335 & 2.7691 & -548.7335 & 2.7691 & -548.7335 & 2.7691 \\
 & $\rho$ & 0.2654 & 0.0339 & 0.2704 & 0.0320 & 0.2760 & 0.0343 & 0.2763 & 0.0389 & 0.2836 & 0.0406 \\
 & $\bar{\rho}^2$ & 0.2417 & 0.0339 & 0.2467 & 0.0319 & 0.2523 & 0.0343 & 0.2526 & 0.0389 & 0.2599 & 0.0406 \\
\cline{1-12}
\bottomrule
\end{tabular}
}

\caption[Performance indicators for level 3 charging with the MNL model, and various utility function thresholds]{Performance indicators for level 3 charging with the MNL model, and various utility function thresholds. The amenity threshold is set to 200m. *** indicates significance at 1\% level, ** significance at 5\% level, and * significance at 10\% level}
    \label{ThirdArticle:Table:ScoresDistances200L3}
\end{table}

\clearpage
\section{Simulation of Utilities}
\label{ThirdArticle:AppendixSimulUtilities}

Recall from Section~\ref{ThirdArticle:SectionModelSpecification} that the utility of a customer $i\in N$ for alternative $j \in M$ has the form
$$u_{ij}=\sum_{k=1}^K \beta_k x_{ijk} + \epsilon_{ij},$$
where $k = 1, \dotsc, K$ are the attributes and $\epsilon_{ij}$ follows a Gumbel distribution with location 0 and a scale $\lambda>0$. The true scale $\lambda$ cannot be identified. In practice, the estimation process returns a set of parameter values with an implicit scale $\hat{\lambda}$, i.e. $\hat \beta_k =\frac{\beta_k}{\hat \lambda}$ and $\hat \epsilon_{ij}=\frac{\epsilon_{ij}}{\hat \lambda}$ with 
$$\hat u_{ij}= \frac{u_{ij}}{\hat{\lambda}} = \sum_{k=1}^K \hat \beta_k x_{ijk} + \hat \epsilon_{ij}.$$
In this way, $ \hat\epsilon_{ij}$ has scale 1.

In our application, we recall that our parameter estimations were over five folds (indexed by $n=0,1,2,3,4$) each of which was estimated separately and, as a consequence, with its own implicit scale $\hat{\lambda}^n$. As these implicit scales are different, the utilities are not directly comparable across folds. For this reason we use ratios in the utilities of our optimization programs; More specifically, for each fold $n$, our simulated utilities are given by 
\begin{align*}
    \sum_{k=1}^K \frac{\hat \beta_k^n}{\hat \beta_{ishome}^n}x_{ijk}+\frac{\hat \epsilon_{ij}^n}{\hat \beta_{ishome}^n} &  = \sum_{k=1}^K \frac{ \frac{ \beta_k}{\lambda}}{\frac{\beta_{ishome}}{\lambda}}x_{ijk}+\frac{\frac{\epsilon_{ij}}{\lambda}}{\frac{\beta_{ishome}}{\lambda}} \\
    & = \sum_{k=1}^K \frac{\beta_k}{\beta_{ishome}}x_{ijk}+\frac{\epsilon_{ij}}{\beta_{ishome}},
\end{align*}
which is independent of the fold. As a consequence, the utility values across all folds are comparable.


Next, we specify the derivation of the parameters for the {\mnl}, {\mxlmean}, and {\mxltwentyfive} utilities, similar to the method in~\cite{Paneque2020}. This process derives a single utility value for each pair of census points and charging station, based on the set of simulated utility values. To this end, we denote by $\hat \epsilon_{ij}^{nr} = \frac{\epsilon_{ij}^{r}}{\hat \beta_{ishome}^n}$ the simulated error term in scenario $r$ and fold $n$, where the error term $\epsilon_{ij}^{r}$ is drawn from a Gumbel distribution with location 0 and scale 1.

\paragraph{Multinomial logit} 
For the {\mnl}, we start by generating a set of utility values via 5,000 draws of $\hat \epsilon_{ij}^{nr}$, with 1,000 in each fold $n$.  In other words, we have 
$$\hat u_{ij}^{nr} = \sum_{k=1}^K \frac{\hat \beta_k^n}{\hat \beta_{ishome}^n}x_{ijk}+\hat \epsilon_{ij}^{nr}.$$
We note that the draws of $\epsilon_{ij}^{nr}, r = 1, \ldots, 1000$ are distinct for each census point $i$, charging station $j$, and fold $n$. 

The utility $u_{ij}$ for census point $i$ and charging station $j$--which acts as the weight coefficients in Models~\eqref{ModelPMedian} and~\eqref{ModelMaxMin}--is then given by the average across the simulated utilities:
$$u_{ij} =\frac{1}{5000}\sum_{n=0}^4 \sum_{r=1}^{1000} u_{ij}^{nr}.$$

\paragraph{Mixed logit}

For both mixed logit models, the process is similar as the multinomial logit. We use a set of 10,000 scenarios each for $\hat \epsilon_{ij}^{nr}$ (2,000 in each fold $n$), with $\hat \epsilon_{ij}^{nr}$ as defined for {\mnl}. However, in this case, the parameter values $\hat \beta_k^{n}$ are not constants, and must be simulated as well. For this, we generate a set of 10,000 draws $\hat \beta_k^{nr}$. As all parameters are assumed to follow a normal distribution, we have that $\hat \beta_k^{nr} = \frac{\hat \beta_k^{\mu n} + \beta_k^{\sigma n} \eta_k^{r} }{\hat \beta_{ishome}^{\mu n}}$ with $\eta_k^{r}$ being draws from a normal distribution with location 0 and scale 1. In other words, we have 
$$\hat u_{ij}^{nr} = \sum_{k=1}^K \frac{\hat \beta_k^{nr}}{\hat \beta_{ishome}^{n}}x_{ijk}+\hat \epsilon_{ij}^{nr}.$$ 
We note that taking the parameter ratios with the mean value of $\beta_{ishome}$ in each fold is sufficient, via identical logic as before. Furthermore, we remark that the draws of $\epsilon_{ij}^{nr}$ are distinct for each census point $i$, charging station $j$, and fold $n$, while the draws of $\hat \beta_k^{nr}$ are distinct for each parameter $k$ but shared in each $i,j$ and $n$.

For the mixed logit cases, we then use the set of utilities $\hat u_{ij}^{nr}$ to generate two sets of parameters values for the utility $u_{ij}$ in Models~\eqref{ModelPMedian} and~\eqref{ModelMaxMin}. For \mxlmean, we use the average of this set, i.e. 
$$u_{ij} =\frac{1}{10000}\sum_{n=0}^4 \sum_{r=1}^{2000} u_{ij}^{nr},$$
while for \mxltwentyfive, we take the 25th percentile of this set.

\section{Additional Optimisation Results}\label{app:Addtional_opt_results}

In this section, we present the complete comparison of the objective value from the optimisation process across all experiments. See Tables~\ref{TableObjectiveComparisonFull_min_max} and~\ref{TableObjectiveComparisonFull}.

\newpage
\pagestyle{empty}
\begin{table}[]
\begin{footnotesize}
\begin{tabular}{llllrrrrr} 
\toprule
 &  &  &  & Objective & Distance & MNL & MXLMean & MXL25 \\
Model & Utilities & Consider existing stations & $p$ &  &  &  &  &  \\
\midrule
\multirow[t]{40}{*}{MaxMin} & \multirow[t]{10}{*}{Distance} & \multirow[t]{5}{*}{False} & 10 & -1.15 & -1.15 & 0.33 & 0.33 & -0.23 \\
 &  &  & 25 & -1.15 & -1.15 & 0.33 & 0.42 & -0.11 \\
 &  &  & 50 & -1.15 & -1.15 & 0.37 & 0.42 & -0.11 \\
 &  &  & 75 & -1.15 & -1.15 & 0.37 & 0.42 & -0.11 \\
 &  &  & 100 & -1.15 & -1.15 & 0.37 & 0.42 & -0.11 \\
\cline{3-9}
 &  & \multirow[t]{5}{*}{True} & 10 & -1.15 & -1.15 & 0.74 & 1.07 & -0.06 \\
 &  &  & 25 & -1.15 & -1.15 & 0.74 & 1.07 & -0.06 \\
 &  &  & 50 & -1.15 & -1.15 & 0.74 & 1.07 & -0.06 \\
 &  &  & 75 & -1.15 & -1.15 & 0.74 & 1.07 & -0.06 \\
 &  &  & 100 & -1.15 & -1.15 & 0.74 & 1.07 & -0.06 \\
\cline{2-9} \cline{3-9}
 & \multirow[t]{10}{*}{MNL} & \multirow[t]{5}{*}{False} & 10 & 0.40 & -1.15 & 0.40 & 0.30 & -0.22 \\
 &  &  & 25 & 0.40 & -1.15 & 0.40 & 0.53 & 0.05 \\
 &  &  & 50 & 0.40 & -1.15 & 0.40 & 0.62 & -0.07 \\
 &  &  & 75 & 0.40 & -1.15 & 0.40 & 0.62 & -0.07 \\
 &  &  & 100 & 0.40 & -1.15 & 0.40 & 0.62 & 0.05 \\
\cline{3-9}
 &  & \multirow[t]{5}{*}{True} & 10 & 0.74 & -1.17 & 0.74 & 1.07 & -0.06 \\
 &  &  & 25 & 0.74 & -1.16 & 0.74 & 1.07 & -0.06 \\
 &  &  & 50 & 0.74 & -1.15 & 0.74 & 1.07 & -0.06 \\
 &  &  & 75 & 0.74 & -1.15 & 0.74 & 1.07 & -0.06 \\
 &  &  & 100 & 0.74 & -1.15 & 0.74 & 1.07 & -0.06 \\
\cline{2-9} \cline{3-9}
 & \multirow[t]{10}{*}{MXL25} & \multirow[t]{5}{*}{False} & 10 & 0.05 & -1.16 & 0.38 & 0.53 & 0.05 \\
 &  &  & 25 & 0.05 & -1.16 & 0.38 & 0.53 & 0.05 \\
 &  &  & 50 & 0.05 & -1.15 & 0.38 & 0.53 & 0.05 \\
 &  &  & 75 & 0.05 & -1.15 & 0.38 & 0.53 & 0.05 \\
 &  &  & 100 & 0.05 & -1.15 & 0.40 & 0.62 & 0.05 \\
\cline{3-9}
 &  & \multirow[t]{5}{*}{True} & 10 & 0.05 & -1.16 & 0.74 & 1.07 & 0.05 \\
 &  &  & 25 & 0.05 & -1.15 & 0.74 & 1.07 & 0.05 \\
 &  &  & 50 & 0.05 & -1.15 & 0.74 & 1.07 & 0.05 \\
 &  &  & 75 & 0.05 & -1.15 & 0.74 & 1.07 & 0.05 \\
 &  &  & 100 & 0.05 & -1.15 & 0.74 & 1.07 & 0.05 \\
\cline{2-9} \cline{3-9}
 & \multirow[t]{10}{*}{MXLMean} & \multirow[t]{5}{*}{False} & 10 & 0.62 & -1.17 & 0.35 & 0.62 & -0.25 \\
 &  &  & 25 & 0.62 & -1.16 & 0.37 & 0.62 & -0.21 \\
 &  &  & 50 & 0.62 & -1.15 & 0.40 & 0.62 & -0.07 \\
 &  &  & 75 & 0.62 & -1.15 & 0.40 & 0.62 & -0.07 \\
 &  &  & 100 & 0.62 & -1.15 & 0.37 & 0.62 & 0.05 \\
\cline{3-9}
 &  & \multirow[t]{5}{*}{True} & 10 & 1.07 & -1.17 & 0.74 & 1.07 & -0.06 \\
 &  &  & 25 & 1.07 & -1.16 & 0.74 & 1.07 & -0.06 \\
 &  &  & 50 & 1.07 & -1.15 & 0.74 & 1.07 & -0.06 \\
 &  &  & 75 & 1.07 & -1.16 & 0.74 & 1.07 & -0.06 \\
 &  &  & 100 & 1.07 & -1.15 & 0.74 & 1.07 & -0.06 \\
\cline{1-9} \cline{2-9} \cline{3-9}
\end{tabular}
\end{footnotesize}
 \caption{Value of the objective function when applying the solution from one set of utilities to a different set. }
    \label{TableObjectiveComparisonFull_min_max}
\end{table}

\newpage
\pagestyle{empty}
\begin{table}[]
 \begin{footnotesize}
\begin{tabular}{llllrrrrr} 
\toprule
 &  &  &  & Objective & Distance & MNL & MXLMean & MXL25 \\
Model & Utilities & Consider existing stations & $p$ &  &  &  &  &  \\
\midrule
\multirow[t]{40}{*}{PMedian} & \multirow[t]{10}{*}{Distance} & \multirow[t]{5}{*}{False} & 10 & -369.22 & -369.22 & 1544.80 & 1532.99 & 66.79 \\
 &  &  & 25 & -350.78 & -350.78 & 1478.09 & 1741.19 & -216.93 \\
 &  &  & 50 & -340.56 & -340.56 & 1619.95 & 1868.58 & 140.10 \\
 &  &  & 75 & -336.17 & -336.17 & 1747.33 & 2158.54 & 173.32 \\
 &  &  & 100 & -333.91 & -333.91 & 1974.53 & 2704.87 & 609.96 \\
\cline{3-9}
 &  & \multirow[t]{5}{*}{True} & 10 & -364.37 & -364.37 & 3725.68 & 4209.09 & 644.95 \\
 &  &  & 25 & -349.63 & -349.63 & 3749.03 & 4214.87 & 655.30 \\
 &  &  & 50 & -340.26 & -340.26 & 3792.59 & 4229.58 & 672.23 \\
 &  &  & 75 & -335.95 & -335.95 & 3826.20 & 4243.67 & 803.63 \\
 &  &  & 100 & -333.79 & -333.78 & 3857.18 & 4254.13 & 811.95 \\
\cline{2-9} \cline{3-9}
 & \multirow[t]{10}{*}{MNL} & \multirow[t]{5}{*}{False} & 10 & 1602.80 & -390.73 & 1602.80 & 1412.24 & -252.21 \\
 &  &  & 25 & 1730.90 & -373.41 & 1730.90 & 1538.70 & -189.46 \\
 &  &  & 50 & 1879.86 & -358.45 & 1879.86 & 2433.08 & 594.53 \\
 &  &  & 75 & 1999.78 & -345.33 & 1999.78 & 2529.72 & 611.47 \\
 &  &  & 100 & 2083.81 & -337.56 & 2083.81 & 2596.72 & 621.69 \\
\cline{3-9}
 &  & \multirow[t]{5}{*}{True} & 10 & 3761.35 & -376.56 & 3761.35 & 4215.24 & 655.11 \\
 &  &  & 25 & 3806.87 & -365.16 & 3806.87 & 4229.70 & 672.12 \\
 &  &  & 50 & 3858.94 & -353.65 & 3858.94 & 4249.94 & 691.73 \\
 &  &  & 75 & 3898.48 & -344.60 & 3898.48 & 4263.21 & 707.66 \\
 &  &  & 100 & 3926.28 & -340.38 & 3926.28 & 4272.30 & 719.11 \\
\cline{2-9} \cline{3-9}
 & \multirow[t]{10}{*}{MXL25} & \multirow[t]{5}{*}{False} & 10 & 555.16 & -392.13 & 1434.83 & 2227.39 & 555.16 \\
 &  &  & 25 & 573.77 & -369.12 & 1562.49 & 2317.26 & 573.77 \\
 &  &  & 50 & 597.06 & -356.46 & 1724.50 & 2427.41 & 597.06 \\
 &  &  & 75 & 614.97 & -344.37 & 1983.13 & 2513.44 & 614.97 \\
 &  &  & 100 & 626.86 & -337.53 & 2056.56 & 2568.38 & 626.86 \\
\cline{3-9}
 &  & \multirow[t]{5}{*}{True} & 10 & 786.26 & -374.82 & 3756.47 & 4216.47 & 786.26 \\
 &  &  & 25 & 797.29 & -361.83 & 3799.21 & 4228.15 & 797.29 \\
 &  &  & 50 & 810.95 & -350.66 & 3844.63 & 4241.99 & 810.95 \\
 &  &  & 75 & 820.17 & -343.66 & 3880.80 & 4251.56 & 820.17 \\
 &  &  & 100 & 826.65 & -340.53 & 3906.34 & 4258.63 & 826.65 \\
\cline{2-9} \cline{3-9}
 & \multirow[t]{10}{*}{MXLMean} & \multirow[t]{5}{*}{False} & 10 & 2441.11 & -392.14 & 1428.23 & 2441.11 & -1082.62 \\
 &  &  & 25 & 2528.42 & -372.23 & 1560.43 & 2528.42 & 570.63 \\
 &  &  & 50 & 2632.04 & -351.97 & 1722.99 & 2632.04 & 590.32 \\
 &  &  & 75 & 2717.87 & -344.89 & 1994.12 & 2717.87 & 609.70 \\
 &  &  & 100 & 2781.32 & -337.97 & 2081.32 & 2781.32 & 619.80 \\
\cline{3-9}
 &  & \multirow[t]{5}{*}{True} & 10 & 4221.15 & -376.13 & 3747.91 & 4221.15 & 649.03 \\
 &  &  & 25 & 4236.95 & -360.66 & 3789.83 & 4236.95 & 789.65 \\
 &  &  & 50 & 4255.00 & -350.68 & 3841.58 & 4255.00 & 800.89 \\
 &  &  & 75 & 4268.62 & -344.89 & 3886.52 & 4268.62 & 811.63 \\
 &  &  & 100 & 4277.81 & -340.14 & 3915.48 & 4277.81 & 819.22 \\
\cline{1-9} \cline{2-9} \cline{3-9}
\bottomrule
\end{tabular}
\end{footnotesize}
 \caption{Value of the objective function when applying the solution from one set of utilities to a different set. }
    \label{TableObjectiveComparisonFull}
\end{table}

\end{document}